\documentclass[11pt,a4paper,reqno]{amsproc}
\usepackage{amssymb, amsmath, amscd}

\newtheorem{theorem}{Theorem}
\newtheorem{proposition}{Proposition}[section]
\newtheorem{lemma}{Lemma}[section]
\newtheorem{corollary}{Corollary}[section]
\newtheorem*{corollary*}{Corollary}
\newtheorem*{theorem*}{Theorem}
\newtheorem*{theoremCartan*}{Cartan's theorem}
\newtheorem*{theoremHarm*}{The Harmonicity Theorem}
\newtheorem*{theoremFinit*}{The Finiteness Theorem}
\newtheorem*{theoremeiconal*}{The Eiconal Cubic Theorem}
\newtheorem*{theoremA*}{Theorem A}
\newtheorem*{theoremB*}{Theorem B}
\newtheorem*{theoremC*}{Theorem C}
\newtheorem*{theoremD*}{Theorem D}
\newtheorem*{theoremRepr*}{The Representation Theorem}
\newtheorem*{theoremClass*}{The Classification Theorem}
\theoremstyle{remark}
\newtheorem*{definition*}{Definition}
\newtheorem*{conjecture*}{Finiteness Conjecture}
\newtheorem{remark}{Remark}[section]
\newtheorem*{remark*}{Remark}
\newtheorem{example}{Example}[section]
\numberwithin{equation}{section}

\theoremstyle{plain}

\def\scal#1#2{\langle #1, #2\rangle}
\def\R#1{\mathbb{R}^{#1}}
\def\Com#1{\mathbb{C}^{#1}}
\def\I{\mathrm{i}}
\def\TOP{\top}

\DeclareMathOperator{\re}{\mathrm{Re}}

\DeclareMathOperator{\trace}{\mathrm{tr}}
\def\SO{{SO}}
\def\OO{{O}}

\def\barr{{\bar x}}

\begin{document}

\title{On a classification of minimal cubic cones in $\R{n}$}
\author{Tkachev Vladimir G.}
\begin{flushright}
\end{flushright}
\address{V. G. Tkachev\\
              Royal Institute of Technology}
\email{tkatchev@kth.se}

\begin{abstract}
We establish a  classification of cubic minimal cones in case of the so-called   radial eigencubics. Our principal result states that any radial eigencubic is either a member of the infinite family of eigencubics of Clifford type, or belongs to one of 18 exceptional families. We prove that at least 12 of the 18 families are non-empty and study their algebraic structure. We also establish that any radial eigencubic satisfies  the  trace identity $\det \mathrm{Hess}^3 (f)=\alpha f$ for the Hessian matrix of $f$, where $\alpha\in \R{}$. Another result of the paper is a correspondence between radial eigencubics and isoparametric hypersurfaces with four principal curvatures.

\end{abstract}
\maketitle


\tableofcontents


%
%

%
%

\tableofcontents

\section{Preliminaries and the main results}

\subsection{Introduction}
In 1969, Bombieri, De Giorgi and Giusti \cite{BGG}  found the first non-affine entire solution of the minimal surface equation
\begin{equation}\label{MSE}
(1+|\nabla u|^2)\Delta u- \sum_{i,j=1}^{n-1}u_{x_{i}x_{j}}u_{x_i}u_{x_j}=0
\end{equation}
Because of its geometric significance, the minimal surface equation (\ref{MSE}) and, especially, Bernstein's problem on the existence of non-affine entire solutions of (\ref{MSE}), have historically attracted perhaps more interest than any other quasilinear elliptic equation.  We refer to \cite{Miranda}, \cite{MMM}, \cite{Nitsche}, \cite{Osserman}, \cite{SimonL}, and the references therein for a detailed discussion of the history of the solution of Bernstein's problem.  Although  many non-affine examples of entire solutions of (\ref{MSE}) for $n\ge 9$ were shown to exist (see, for instance, \cite{Simon89}, \cite{SS}), no explicit examples have been constructed.
L.~Simon \cite{Simon89} established  that for $n=9$ all entire solutions of (\ref{MSE}) are of  polynomial growth and it is a long-standing conjecture that   this property holds in general \cite{BG},  \cite{Osserman}. Even a simpler question  \cite{SimonL}, \cite{Miranda}, whether or not there exists a  solution of (\ref{MSE}) which is actually a polynomial in $x_i$, is still unanswered.

These questions  prompt one to study algebraic minimal hypersurfaces, and, in particular, algebraic minimal cones. The latter occure naturally as singular `blow-ups' of entire solutions of (\ref{MSE}) at infinity; for example, the  seven-dimensional Simon's  cone $\{(x,y)\in \R{4}\times \R{4}: |x|^2=|y|^2\}$ played an important role in the solution of Bernstein's problem \cite{Fleming}, \cite{Simons} and  in the constructing of non-affine examples by Bombieri, de Giorgi and Giusti \cite{BGG}.  We mention also a recent appearance of algebraic minimal hypersurfaces as selfsimilar solutions of the mean curvature flow in codimension one \cite{Smoczyk}. Note also that any progress in algebraic minimal cones leads  to a better understanding of algebraic aspects of minimal submanifolds of codimension one in the unit spheres because of the well-known correspondence between this objects.

Minimal cones of lower degrees were classified by Hsiang \cite{Hsiang67}:  the only first degree  minimal cones are hyperplanes in $\R{n}$, and the only (up to a congruence in $\R{n}$)  quadratic minimal cones are given by the zero-locus $g^{-1}(0)$ of the  quadratic forms
\begin{equation}\label{ClSim}
g(x)=(n-p-1)(x_1^2+\ldots+x_p^2)-(p-1)(x_{p+1}^2+\ldots+x_{n}^2), \quad 2\le p\le n-1.
\end{equation}

On the other hand,  a classification (and even construction) of algebraic minimal cones of degree higher than three remains a long-standing difficult problem \cite{SimonL}, \cite{Hsiang67}, \cite{Fomenko}, \cite{Osserman}.
The lack of `canonical' normal forms for higher degree polynomials  makes a classification of algebraic  minimal cones, at first sight, defeating. On the other hand, a close analysis of the available examples of minimal cubic cones, see e.g. \cite{Taka}, \cite{Hsiang66}, \cite{Hsiang-Lawson}, \cite{Lawson}, reveals that these cones  have a rather distinguished algebraic structure which, in some content, resembles that of isoparametric hypersurfaces in the spheres.  In \cite{Hsiang67}, Hsiang began developing a systematic approach  to study real algebraic minimal submanifolds of degree higher than two and by using the geometric invariant theory  constructed new examples of  non-homogeneous minimal cubic cones  in $\R{9}$ and $\R{15}$. According  Hsiang, the study of real algebraic minimal cones is equivalent to a classification  of  polynomial solutions  $f=f(x_1,\ldots,x_n)\in \R{}[x_1,\ldots,x_n]$ of  the following congruence:
\begin{equation}\label{Lgamma1}
L(f) \equiv 0 \mod f,
\end{equation}
where
\begin{equation*}
L(f)= |\nabla f|^2\Delta f-
\sum_{i,j=1}^n f_{x_i} f_{x_j} f_{x_ix_j}
\end{equation*}
is the normalized mean curvature operator and (\ref{Lgamma1}) is understood in the usual sense, i.e. $L(f)$ is divisible by $f$ in the polynomial ring $\R{}[x_1,\ldots, x_n]$.  Observe, that (\ref{Lgamma1}) geometrically means   that the zero-locus $f^{-1}(0)$ has zero mean curvature everywhere where the gradient $\nabla f\ne 0$.

A polynomial solutions $f\not\equiv 0$ of (\ref{Lgamma1}) is called an \textit{eigenfunction} of $L$. The ratio  $L(f)/f$ (which  is obviously  a polynomial in $x$) is called the \textit{weight} of an eigenfunction $f$. An eigenfunction $f$ which  is a cubic homogeneous polynomial  is also called an \textit{eigencubic}.

\begin{remark}
For geometric reasons, we  make no distinction between two eigenfunctions $f_1$ and $f_2$ which give rise to two congruent cones $f_1^{-1}(0)$ and $f_2^{-1}(0)$; such eigenfunctions will also be called \textit{congruent}. It follows from the real Nullstellensatz \cite{Milnor} (see also \cite[Proposition~2.4]{TkCliff}) that two irreducible homogeneous cubic polynomials $f_1$ and $f_2$ are congruent if and only if there exists an orthogonal endomorphism of $U\in \OO(\R{n})$ and a constant $c\in \R{}$, $c\ne0$, such that  $f_1(x)=cf_2(Ux)$.
\end{remark}

In \cite{Hsiang67}, Hisang observes  that all available {cubic} minimal cones arise as solutions of  the following   non-linear equation:
\begin{equation}\label{mainlambda}
L(f)=\lambda  |x|^2 f,\qquad \lambda \in \R{},
\end{equation}
and poses the problem to  determine all solutions of (\ref{mainlambda}) up to congruence  in $\R{n}$.
We call the solutions of (\ref{mainlambda}) \textit{radial eigencubics}.

It is the purpose of the present paper to provide a general framework for a classification of  radial eigencubics. We prove that any radial eigencubic $f$ is a harmonic polynomial and associate to it a pair $(n_1,n_2)$ of non-negative integers, called the type of $f$. We show that the type is a congruence invariant of $f$  and establish that $n_1$ can be recovered by  the following remarkable trace identity:
$$
\trace \mathrm{Hess}^3 (f)=3(n_1-1)\lambda\, f,
$$
where $\mathrm{\mathrm{Hess}} (f)$ is the Hessian matrix of $f$ and the constant factor $\lambda$ is the same as in (\ref{mainlambda}).
The principal result of the paper states that any radial eigencubic is either a member of the infinite family of eigencubics of Clifford type introduced and classified recently in \cite{TkCliff}, or belongs to one of 18 exceptional families which types $(n_1,n_2)$  and ambient dimensions $n$ listed in Table~\ref{tabsBas} below. We also establish that at least 12 of the 18 families are non-empty and  provide examples of radial eigencubics for each realizable family.

\def\MM{6mm}
\begin{table}[ht]
\renewcommand\arraystretch{1.5}
\noindent
\begin{small}
\begin{tabular}{p{3mm}||p{2.5mm}|p{2.5mm}|p{2.5mm}|p{2.5mm}|p{2.5mm}|p{2.5mm}|p{2.5mm}|p{2.5mm}|p{2.5mm}|p{2.5mm}|p{2.5mm}|p{2.5mm}|p{2.5mm}|p{2.5mm}|p{2.5mm}|p{2.5mm}|p{2.5mm}|p{2.5mm}|p{2.5mm}}
$n_1$& $2$ & $3$  & $5$  & $9$  & $0$ & $1$  & $2$  & $4$  & $0$  & $1$  &    $5$ & $9$  & $0$  & $1$   &   $3$   &  $1$  &  $3$ & $7$ \\\hline
$n_2$&  $0$ & $0$  & $0$  & $0$  &$5$ & $5$  & $5$  & $5$  & $8$  & $8$  &    $8$ & $8$  &$14$  & $14$  &  $14$  &  $26$ &  $26$ & $26$      \\\hline
$n$  &  $5$ & $8$  & $14$  & $26$  &$9$ & $12$ & $15$ & $21$ & $15$ & $18$ &   $30$&  $42$& $27$ & $30$  &  $36$  & $54$ &  $60$ & $72$       \\\hline
 &   &   &   &   & &  & ? &  &  &   &  $ ?$  & $ ?$  &   &    &  $ ?$   &   &  $ ?$  &  $ ?$     \\
\end{tabular}
\end{small}
\bigskip
\caption{\small Exceptional eigencubics:  $?$ stands for the unsettled cases.}
\label{tabsBas}
\end{table}

As was already mentioned, a classification of general radial eigencubics closely resembles that of isoparametric hypersurfaces with four principal curvatures. The isoparametric hypersurfaces have been intensively studied for several decades now and, at the present, a complete classification is available for all but for four exceptional isoparametric families, see  \cite{Mun1}, \cite{Mun2}, \cite{OT1}, \cite{OT2}, \cite{Abresh}, \cite{Stoltz}, \cite{Cecil}, \cite{Cecil1}. It would be interesting to work out an explicit correspondence between these theories. We mention that, in one direction, a  theorem of Nomizu \cite{Nomi} states that each focal variety of an isoparametric hypersurface is a minimal submanifold of the ambient unit sphere. On the other hand, in the present paper we show that, in the other direction, to any {non-isoparametric} radial eigencubic one can associate an isoparametric hypersurface with four principal curvatures. Combining the latter correspondence  with a  deep characterization of isoparametric quartics obtained recently by T.~Cecil, Q.S.Chi and G.~Jensen \cite{CJ}, and by S.~Immerwoll \cite{Immer}, we obtain an obstruction to the existence of some exceptional families of radial eigencubics.

\begin{remark}
We would like to emphasize that, in general, real algebraic minimal cones have a much more rich structure than isoparametric hypersurfaces.  Indeed, in the former case, the examples constructed recently in \cite{TkCliff} show that there exist irreducible minimal cones of arbitrary high degree, while  the well-known theorem of to M\"unzner  \cite{Mun1} allows the defining polynomials of isoparametric hypersurfaces to be  only of degrees  $g=1, 2, 3, 4$ and $6$.
\end{remark}

In section~\ref{sec:main} below  we consider our  results in more detail. First we recall some basic facts about eigencubics of Clifford type and Cartan's isoparametric polynomials.

\subsection{Eigencubics of Clifford type}
A system of symmetric endomorphisms $\mathcal{A}=\{A_{i}\}_{0\le i\le q}$ of $\R{2m}$ is called a {symmetric Clifford system} \cite{Huse}, \cite{Cecil1}, \cite{Baird}, equivalently $\mathcal{A}\in \mathrm{Cliff}(\R{2m},q)$ if
$$
A_i A_j+A_jA_i=2\delta_{ij}\cdot \mathbf{1}_{\R{2m}},
$$
where $\mathbf{1}_{V}$ stands for  the identity operator in a linear space $V$.
To any symmetric Clifford system  $\mathcal{A}\in \mathrm{Cliff}(\R{2m},q)$ with two distinguished elements $A_0, A_1\in \mathcal{A}$ one can associate an orthogonal eigen-decomposition $\R{2m}=\R{m}\oplus \R{m}$ such that
\begin{equation}\label{correspond00}
A_0=\left(
        \begin{array}{cc}
          \mathbf{1}_{\R{m}} & 0 \\
          0 & -\mathbf{1}_{\R{m}}\\
        \end{array}
      \right),\qquad
A_1=\left(
        \begin{array}{cc}
          0 & \mathbf{1}_{\R{m}} \\
          \mathbf{1}_{\R{m}} & 0\\
        \end{array}
      \right),\qquad
A_i=\left(
        \begin{array}{cc}
          0 & P_i \\
          P_i^\TOP & 0\\
        \end{array}
      \right).
\end{equation}
where the skew-symmetric  transformations $P_1$, \ldots, $P_{q-1}$  satisfy $P_iP_j+P_jP_i=-2\delta_{ij}$. This  determines a  representation $\{P_1, \ldots, P_{q-1}\}$ of the Clifford algebra $\mathrm{Cl}_{q-1}$ on $\R{m}$ \cite{Huse}, \cite{Baird}. Conversely, any representation of the Clifford algebra $\mathrm{Cl}_{q-1}$ on $\R{m}$ induces a symmetric Clifford system by virtue of (\ref{correspond00}). It follows from  the representation  theory of Clifford algebras that the class $\mathrm{Cliff}(\R{2m},q)$ is non-empty if an only if  \begin{equation}\label{rhommm}
q\le \rho(m),
\end{equation}
where  the Hurwitz-Radon function $\rho$ is defined by
\begin{equation}\label{foll}
\rho(m)=8a+2^b, \qquad \text{if} \;m=2^{4a+b}\cdot \mathrm{odd} , \;\; 0\leq b\le 3.
\end{equation}

Two symmetric Clifford systems $\mathcal{A}\in \mathrm{Cliff}(\R{2m},q)$ and $\mathcal{B}\in \mathrm{Cliff}(\R{2m'},q')$ are called \textit{geometrically equivalent}, if $q=q'$, $m=m'$, and there exist  orthogonal endomorphisms $U\in \OO(\R{2m})$ and $u\in \OO(\R{q+1})$ such that
\begin{equation*}\label{geom1}
A_{uz}=U^\TOP B_z U,\quad \forall z\in \R{q+1},
\end{equation*}
where $A_z=\sum_{i=0}^q z_iA_i$, $B_z=\sum_{i=0}^q z_iB_i$.
Then the cardinality $\kappa(m,q)$ of the quotient set of $\mathrm{Cliff}(\R{2m},q)$ with respect to the geometric equiavlence is equal to 1 for $q=0$ and is determined for $q\ge 1$ by the following  formula (see also \cite[\S~4.7]{Cecil1}):
\begin{equation}\label{kappa}
\kappa(m,q)
=\left\{
  \begin{array}{ll}
    0, & \hbox{if $\delta(q)\nmid m$;} \\
    1, & \hbox{if $\delta(q)\mid m$ and $q\not\equiv 0 \mod 4$;} \\
    \lfloor\frac{m}{2\delta(q)}\rfloor+1, & \hbox{if $\delta(q)\mid m$ and $q\equiv 0 \mod 4$,}
      \end{array}
\right.
\end{equation}
where $\lfloor x\rfloor$ is the  integer part of  $x$ and  $\delta(q)=\min\{2^k: \rho(2^k)\ge q\},$
or equivalently by the following table \cite[p.~156]{Baird}:
\begin{center}
\begin{tabular}{c|r|r|r|r|r|r|r|r|r|c}
  $q$ & 1 & 2 & 3 & 4 & 5 & 6 & 7 & 8 &\ldots& $k$ \\\hline
  $\delta(q)$ & 1 & 2 & 4 & 4 & 8 &8 & 8 & 8 &\ldots& $16\delta(k-8)$ \\
\end{tabular}
\end{center}

In \cite{TkCliff} we associated to a {Clifford symmetric system}  $\mathcal{A}\in \mathrm{Cliff}(\R{2m},q)$ the cubic form
\begin{equation}\label{Clifffff}
C_\mathcal{A}(x):=\sum_{i=0}^{q}  \scal{y}{A_iy}\,x_{i+1}, \qquad  y=(x_{q+2},\ldots,x_{q+1+2m})\in \R{2m},
\end{equation}
and proved that $C_\mathcal{A}$ is a radial eigencubic  in $\R{2m+q+1}$.
\begin{definition*}
An arbitrary radial eigencubic is said to be of \textit{ Clifford type} if it is congruent to some $C_{\mathcal{A}}$. Otherwise it is called an \textit{exceptional} radial eigencubics.
\end{definition*}

We also proved in \cite{TkCliff} that the congruence classes of eigencubics of Clifford type are in one-to-one correspondence with the equivalence classes of geometrically equivalent Clifford systems which, in view of the remarks made above, yields a complete classification of eigencubics of Clifford type.

\subsection{The Cartan isoparametric polynomials}\label{sec:Cartan}
In \cite{Cartan}, E.~Cartan  proved  that, up to congruence, the only \textit{irreducible} cubic polynomial solutions  of the isoparametric system
\begin{equation}\label{Mun2}
|\nabla f|^2=9x^{4}, \qquad \Delta f=0
\end{equation}
are the following four polynomials:
\begin{equation}\label{CartanFormula0}
\begin{split}
\theta_\ell =&x_{1}^3+\frac{3x_1}{2}(|z_1|^2+|z_2|^2-2|z_3|^2-2x_{2}^2)
+\frac{3\sqrt{3}}{2}[x_{2}(|z_1|^2-|z_2|^2)+\re z_1z_2z_3],
\end{split}
\end{equation}
where $z_k=(x_{k\ell-\ell+3},\ldots,x_{k\ell+2})\in \R{\ell }=\mathbb{F}_\ell$,  and $\mathbb{F}_\ell$ denote the division algebra of dimension $\ell$ (over reals): $\mathbb{F}_1=\mathbb{R}$ (reals), $\mathbb{F}_2=\mathbb{C}$ (complexes), $\mathbb{F}_4=\mathbb{H}$ (quaternions) and $\mathbb{F}_8=\mathbb{O}$ (octonions). The real part in (\ref{CartanFormula0}) should be understood for a general $\mathbb{F}_\ell$ as
\begin{equation}\label{realpart}
\re z_1z_2z_3=\frac{1}{2}((z_1z_2)z_3+\bar z_3(\bar z_2\bar z_1))=\frac{1}{2}(z_1(z_2z_3)+(\bar z_3\bar z_2)\bar z_1),
\end{equation}
(observe that the real part is associative, see also Lemma~15.12 in \cite{Adams}).

It follows from (\ref{CartanFormula0}) that the Cartan isoparametric cubics are well-defined only if the ambient dimension $n\in \{5,8,14,26\}$. Moreover, in virtue of (\ref{Mun2})
$$
L(f)=-\frac{1}{2}\scal{\nabla |\nabla f|^2}{\nabla f}=-18\scal{x}{\nabla f}=-54 x^2f,
$$
hence any Cartan polynomial is also a radial eigencubic. It is easily seen that any  $\theta_\ell$ is in fact an exceptional eigencubic. Indeed, we note  that the squared norm of the gradient  is a congruence invariant and  $|\nabla \theta_\ell|^2=9x^4$, while (\ref{Clifffff}) yields that   the squared norm of the gradient of an eigencubic of Clifford type is at most quadratic in some variables.

A crucial role in our further analysis will play the  following generalization of the Cartan theorem obtained by the author in \cite{TkCartan}.

\begin{theoremeiconal*}
Let $f(x)$ be a  cubic polynomial solution of the first equation in $(\ref{Mun2})$ alone. Then $f$ is either reducible and  congruent to  $x_n(x_n^2-3x_1^2-\ldots-3x_{n-1}^2)$, or irreducible and congruent to some Cartan polynomial $\theta_\ell(x)$.
\end{theoremeiconal*}

\subsection{Main results}\label{sec:main}
As the  first step we  obtain the following characterization of general radial eigencubics.
\begin{theorem}\label{thA}
Any radial eigencubic in $\R{n}$ is a harmonic function.
\end{theorem}

\begin{remark}
Observe, however, that there are (non-radial) eigencubics which are  non-harmonic, for example, $f=x_1 g(x)$ with $g$ given by (\ref{ClSim}). All such non-harmonic eigenfunctions are reducible, so it would be interesting  to know whether there exist \textit{irreducible} non-harmonic eigencubics.
\end{remark}

Our next step is to establish (Proposition~\ref{pro:main})  that given a radial eigencubic $f$ in $\R{n}$ one can associate the orthogonal coordinates $\R{n}=\mathrm{span}(e_n)\oplus V_1\oplus V_2\oplus V_3$ in which $f$ takes  the so-called \textit{normal form}
\begin{equation}\label{given0}
f=x_n^3+\phi x_n+\psi\equiv x_n^3-\frac{3}{2}x_n(2\xi^2+\eta^2-\zeta^2)+\psi_{111}+\psi_{102}+\psi_{012}+\psi_{030},
\end{equation}
where  $x=(\xi, \eta, \zeta, x_n)\in \R{n}$,
$\xi\in V_1$, $\eta\in V_2$,   $\zeta\in V_3$ and   $\psi_{ijk}$ denotes a cubic form of homogeneous class $\xi^i\otimes \eta^j\otimes\zeta^k$. (Here and in what follows, if no ambiguity possible, we abuse the norm notation by writing, e.g.,  $\xi^2$ for $|\xi|^2$). In addition, the harmonicity of $f$ yields the following restrictions:
\begin{equation}\label{n1n2n3}
n_3=2n_1+n_2-2,\quad n=3n_1+2n_2-1.
\end{equation}

\begin{definition*}
The pair $(n_1,n_2)$ is called the \textit{type} of the normal form.
\end{definition*}

Thus a very natural question appears from the very beginning: whether the type $(n_1,n_2)$ has an invariant meaning? We answer this question in positive, but  what is more important, we establish the following remarkable trace identity for determining of the dimension $n_1$.

\begin{theorem}\label{thB}
Let $f$ be a radial eigencubic in the normal form $(\ref{given0})$. Then the associated dimensions $n_i=\dim V_i$, $i=1,2,3$, do not depend on a particular choice of the normal form of $f$ and can be recovered by virtue of the cubic trace formula
\begin{equation}\label{tracecubicformula}
n_1=1+\frac{x^2\trace \mathrm{Hess}^3(f)}{3L(f)},
\end{equation}
and relations
$n_2=\frac{1}{2}(n-3n_1+1)$, $n_3=2n_1+n_2-2,$
where $\mathrm{Hess}(f)$ is the Hessian matrix of $f$.
\end{theorem}

We emphasize that the ratio in (\ref{tracecubicformula}) is an integer number.

In  \cite{TkCliff}, we established the cubic trace identity for eigencubics of Clifford type. In the present paper, we extend the cubic trace  identity to the general radial eigencubics. Our argument  is heavily based on the characterization of exceptional eigencubics by means of the $\psi_{030}$-term in (\ref{given0}) which we now describe. Let us assume that $f$ be a radial eigencubic written in the normal form. First we show in Proposition~\ref{cor:n1n2constraint} that the combination  $\psi_{111}+\frac{1}{\sqrt{3}}\psi_{102}$ induces a symmetric Clifford system in $\mathrm{Cliff}(\R{2(n_1+n_2-1)},n_1-1)$ which immediately yields   by virtue of (\ref{rhommm}) the inequality
\begin{equation}\label{obstr}
n_1-1\le \rho(n_2+n_1-1).
\end{equation}
Next we prove (Proposition~\ref{pr:hidden}) that  the cubic form $\psi_{030}$ in (\ref{given0}) satisfies an eiconal type equation $|\nabla \psi_{030}(\eta)|^2=\frac{9}{2}\eta^2$, $\eta\in \R{n_2}$. Combining  these observations and some further properties of the normal form, we are able to prove the following important characterization of the Clifford eigencubics.

\begin{theorem}\label{thC}
A radial eigencubic $f$  is of Clifford type if and only if  for any particular choice of its normal form  $(\ref{given0})$, the component $\psi_{030}\not\equiv 0$ and reducible.
\end{theorem}

In particular, by combining Theorem~\ref{thC} with the Eiconal Cubic Theorem above, one obtains that for any exceptional eigencubic there holds  $n_2\in \{0,5,8,14,26\}$. Then by using some special properties of the Hurwitz-Radon function $\rho$ one is able to show that there exists only finitely many pairs $(n_1,n_2)$ satisfying the above inclusion and  (\ref{obstr}). This yields the finiteness of the number of types of exceptional eigencubics. In fact, we have the the following criterion.

\begin{theorem}\label{thD}
Let $f$ be a radial eigencubic in $\R{n}$. Then the following statements are equivalent:
\begin{itemize}
\item[(a)]
$f$ is an exceptional radial eigencubic;
\item[(b)]
for any choice of the normal form $(\ref{given0})$, the form $\psi_{030}$ is either irreducible  or identically zero;
\item[(c)]
$n_2\in \{0,5,8,14,26\}$ and the quadratic form
$$
\sigma_2(f):=-\frac{1}{3\lambda}\mathrm{Hess}^2(f), \quad \text{where} \,\, L(f)=\lambda x^2 f,
$$
has a single eigenvalue.
\end{itemize}
The only possible types of exceptional eigencubics are those displayed in Table~\ref{tabs} below.
\end{theorem}

In the remaining part of the paper we investigate which of the 23 pairs $(n_1,n_2)$ in Table~\ref{tabs} are indeed realizable as the types of exceptional eigencubics. Below we summarize the corresponding results.

\begin{itemize}
\item[(i)]
 For $n_2=0$ all  types $(\ell+1,0)$, $\ell=1,2,4,8$, are realizable. For each $\ell$, there is exactly one congruence class of exceptional eigencubics of type $(\ell+1,0)$ represented by the Cartan polynomial $\theta_\ell$.
\item[(ii)]
For $n_1=0$  the only three types $(0,5)$, $(0,8)$ and $(0,14)$ are realizable.
\item[(iii)]
For  $n_1=1$ then  four types $(1,5)$, $(1,8)$, $(1,14)$ and $(1,26)$ are realizable  and in each case there is exactly one congruence class of exceptional eigencubics.
\item[(iv)]
There is an exceptional eigencubic of type of type $(4,5)$.
\item[(v)]
The types $(2,8)$, $(2,14)$, $(2,26)$ and $(3,8)$ are not realizable.
\end{itemize}

Thus, it remains unsettled  the six exceptional pairs: $(2,5)$, $(5,8)$, $(9,8)$, $(3,14)$, $(3,26)$, $(7,26)$.

To obtain the non-existence result (v) we develop a correspondence between general radial eigencubics with $n_2\ne 0$  and isoparametric quartic polynomials which can be described as follows. Recall that a hypersurface in the unit sphere in $\R{n}$ is called \textit{isoparametric} if it has constant principal curvatures \cite{Thorb}, \cite{Cecil}. A celebrated theorem of M\"{u}nzner \cite{Mun1} states that any isoparametric hypersurface is  algebraic and its defining polynomial $h$ is homogeneous of degree $g=1,2,3,4$ or $6$, where $g$ is the number of distinct principal curvatures. Moreover, if $g=4$ then, suitably normalized, $h$  satisfies the system M\"{u}nzner-Cartan differential equations (cf. with (\ref{Mun2} above)
\begin{equation}\label{Muntzer4}
|\nabla h|^2=16x^{6}, \qquad \Delta h=8(m_2-m_1)x^{2}, \quad x\in \R{n},
\end{equation}
where $m_i$ are the multiplicities of the maximal and minimal principal curvature of $M$, $m_1+m_2=\frac{n-2}{2}$. Let $\mathrm{Isop}(m_1,m_2)$ denote the class of all quartic polynomials satisfying (\ref{Muntzer4}). Then  each $h\in \mathrm{Isop}(m_1, m_2)$ with $m_1,m_2\ge 1$ gives rise to a family of isoparametric hypersurfaces
$$
M_c=\{x\in \mathbb{{S}}^{n-1}\subset\R{n}|\,h(x)=c\}, \qquad c\in (-1,1),
$$
see for instance \cite[p.~96-97]{Cecil1}. In this case $m_1$ and $m_2$ are, up to a permutation, the multiplicities of the maximal and minimal principal curvature of the hypersurface $M_c$ and
\begin{equation}\label{dimension1}
n-2=2(m_1+m_2).
\end{equation}

\begin{theorem}\label{thE}
Let $f$ be any radial eigencubic in $\R{n}$ of type $(n_1,n_2)$, $n_2\ne 0$, and normalized by $\lambda =-8$ in $(\ref{mainlambda})$. Then $f$ can be written in some orthogonal coordinates in the degenerate form
\begin{equation}
f=(u^2-v^2)x_n+a(u,w)+b(y,w)+c(u,y,w),
\end{equation}
where $u=(x_1,\ldots, x_m)$, $v=(x_{m+1},\ldots, x_{2m})$, $w=(x_{2m+1},\ldots, x_{n-1})$, and the cubic forms
$a\in u\otimes w^2$, $b\in v\otimes w^2$, $c\in u\otimes v\otimes w$. Moreover, the quartic polynomials
\begin{equation*}
\begin{split}
h_0(u,v)&:= (u^2+v^2)^2-2 c_w^2\in \mathrm{Isop}(n_1-1,m-n_1),\\
h_1(u,v)&:= -u^4+6u^2v^2-v^4-2 c_w^2\in \mathrm{Isop}(n_1,m-n_1-1).
\end{split}
\end{equation*}
If $f$ is in addition an exceptional eigencubic then $n_2=3\ell+2$, $\ell\in \{1,2,4,8\}$ and $m=\ell+n_1+1$.
\end{theorem}

By using  a recent classification result of T.~Cecil, Q.S.Chi and G.~Jensen \cite{CJ}, and S.~Immerwoll \cite{Immer}, and Theorem~\ref{thE}, one can deduce the nonexistence of types mentioned in (v).

The paper is organized as follows.  In section~\ref{sec:harm} we prove Theorem~\ref{thA} and in section~\ref{sec:def} we establish the normal representation (\ref{given0}). In Proposition~\ref{cor:n1n2constraint} we exhibit a hidden Clifford structure associated with any radial eigencubic and prove (\ref{obstr}). In Proposition~\ref{pro:n20} we obtain a complete classification of  radial eigencubics with $n_2=0$  mentioned in the item (i) above. The proofs of  Theorem~\ref{thB}, Theorem~\ref{thC} and Theorem~\ref{thD} will be given in sections~\ref{sec:ifpart} and ~\ref{sec:trace}. In section~\ref{sec:someexamples} we establish the classification results (ii)--(iv) and also  review some examples of exceptional radial eigencubics and  outline some their aspects. In section~\ref{sec:isop} we prove Theorem~\ref{thE} and the non-existence result (v).

\bigskip
\textbf{{Notation.}} We use the standard convention that $f_\xi$ denotes the vector-column of partial derivatives $f_{\xi_i}$ and  $f_{\xi\eta}$ stands for  the Jacobian matrix with entries $f_{\xi_i\eta_j}$, etc. By $\Delta_\xi f=\trace f_{\xi\xi}$ we denote the Laplacian with respect to $\xi$. We suppress the variable notation for the full gradient gradient $\nabla f$, the Hessian matrix $\mathrm{Hess} f$ and the full Laplacian $\Delta f$.
In what follows,  if no ambiguity possible, we abuse the norm notation by writing, e.g.,  $\xi^2$ for $|\xi|^2$.
The bar notation $\bar x$ is usually  used

\section{The harmonicity of radial eigencubics}\label{sec:harm}

We begin with treating the normal form of a radial eigencubic. To this end, let us consider an arbitrary radial eigencubic $f$ and let $x_0\in \mathbb{{S}}^{n-1}$ be a maximum point of $f$ on the unite sphere $\mathbb{{S}}^{n-1}$. It is well known and easily verified that in any orthogonal  coordinates with $x_0$ chosen to be the $n$-th basis vector, $f$ expands as follows:
\begin{equation}\label{reduced}
f(x)=cx_n^3+x_n\phi(\bar x)+\psi(\bar x), \qquad \bar x=(x_1,\ldots, x_{n-1}).
\end{equation}
We shall refer to (\ref{reduced}) as the \textit{normal form} of $f$.
Using the freedom to scale $f$, we can ensure that $c=1$.  Rewrite the definition of radial eigencubic as follows:
\begin{equation}\label{mainproblem}
L(f)=18\alpha x^2 f,
\end{equation}
where the factor $\lambda(f)=18\alpha$ is chosen for the further convenience. Then, identifying the coefficients of $x_n^i$, $0\le 1\le 5$, in (\ref{mainproblem})  we arrive at the following system
\begin{eqnarray}\label{e1}
&&\Delta_{\bar x} \phi=2\alpha,\\
\label{e2}
&&\Delta_{\bar x} \psi=0,\\
\label{e3}
&&\phi_{\bar x}^2\Delta_{\bar x}\phi-\phi_{\bar x}\phi_{\bar x\bar  x}\phi_{\bar x}=6\alpha (\phi+3{\bar x}^2),\\
\label{e4}
&&2\phi_{\bar x}^{\TOP} \phi_{{\bar x}{\bar x}}\psi_{\bar x}+\phi_{\bar x}^{\TOP} \psi_{{\bar x}{\bar x}}\phi_{\bar x}-(6+4\alpha)\phi_{\bar x}^{\TOP} \psi_{{\bar x}}+18\alpha \psi=0,\\
\label{e5}
&&2(3+\alpha)\psi_{\bar x}^2 -\psi_{\bar x}^{\TOP} \phi_{{\bar x}{\bar x}}\psi_{\bar x}-2\phi_{\bar x}^{\TOP} \psi_{{\bar x}{\bar x}}\psi_{\bar x}=2\phi(9\alpha {\bar x}^2 +\phi_{\bar x}^2-\alpha \phi),\\
\label{e6}
&&2\phi\, \phi_{\bar x}^{\TOP} \psi_{{\bar x}}+18\alpha   {\bar x}^2\psi+\psi_{\bar x}^{\TOP} \psi_{{\bar x}{\bar x}}\psi_{\bar x}=0.
\end{eqnarray}
We may choose orthogonal coordinates in $\R{n}$ such that the quadratic form  $\phi$ becomes diagonal, say $\phi( x)=\sum_{i=1}^{n-1}\phi_ix_i^2$. Then  (\ref{e1}) yields
\begin{equation}\label{a1}
\sum_{i=1}^{n-1}\phi_i=\alpha.
\end{equation}
By expanding (\ref{e3}), we see that each eigenvalue $\phi_i$ satisfies the equation
\begin{equation}\label{char1}
\chi_\alpha(t):=4t^3-4\alpha t^2+3\alpha t+ 9\alpha=0.
\end{equation}

First notice that we can always assume that $\alpha\ne0$ because otherwise (\ref{char1}) yields $\chi_0\equiv 4t^4$, hence $\phi_i=0$ for all $i$ and thus  $\phi\equiv 0$. But the latter implies  $\psi_{\bar x}^2=0$ by virtue of  (\ref{e5}), hence $\psi\equiv0$. This yields $f=x_n^3$, i.e. $n=1$, a contradiction.

Thus, assuming $\alpha\ne 0$, we denote by $ t_i$, $1\le i\le \nu(\alpha)$, all distinct \textit{real} roots of (\ref{char1}). Since  (\ref{char1}) is  a cubic equation with real coefficients, one has  $1\le \nu(\alpha)\le 3$. Regarding $t_i$ as an eigenvalue of $\phi$, let  $V_i$ denote the corresponding eigenspace ($V_i$ may be null-dimensional). Then
\begin{equation}\label{decomp}
\R{n}=\mathrm{span}(e_n)\oplus V, \qquad V=\bigoplus_{i=1}^{\nu(\alpha)}V_i.
\end{equation}
From (\ref{a1}) we infer the following constraints on the dimensions $n_i=\dim V_i$:
\begin{equation}\label{a11}
\sum_{i=1}^{\nu(\alpha)} t_in_i=\alpha, \qquad \sum_{i=1}^{\nu(\alpha)}n_i=n-1.
\end{equation}

Now we are going to specify the algebraic structure of the cubic form $\psi$. Note that the eigen decomposition (\ref{decomp}) extends to the tensor products, thus we have for the cubic forms:
$$
{V^*}^{\otimes 3}=\bigoplus_{|q|=3}{V^*}^{\otimes q}, \quad {V^*}^{\otimes q}:=\bigotimes_{i=1}^{\nu(\alpha)}{V_{i}^*}^{q_i},
$$
where $q=(q_1,\ldots,q_{\nu(\alpha)})$ and  $|q|=q_1+\ldots+q_{\nu(\alpha)}=3$. Write $\psi=\sum_{|q|=3}\psi_q$ according to the above decomposition.

\begin{lemma}\label{prop1}
In the above notation, let
\begin{equation}\label{rdef}
R_q:=\frac{9\alpha}{2}+(\sum_{k=1}^{\nu(\alpha)} t_kq_k)^2+\sum_{k=1}^{\nu(\alpha)} t_k^2 q_k
-(3+2\alpha)\sum_{k=1}^{\nu(\alpha)} t_kq_k.
\end{equation}
If $R_q\ne 0$ for some $q$, $|q|=3$, then the corresponding homogeneous component $\psi_q$ is zero. In other words, $\psi$ is completely determined by the homogeneous components  $\psi_q$ whose indices $q$ satisfy $R_q= 0$.
\end{lemma}

\begin{proof}
By virtue of the Euler homogeneous function theorem,
\begin{equation*}\label{cphi}
\begin{split}
\phi_{\bar x}^{\TOP} \phi_{{\bar x}{\bar x}}(\psi_q)_{\bar x}
&=4\psi_q\sum_{k=1}^{\nu(\alpha)} t_k^2q_k,\\
\phi_{\bar x}^{\TOP}(\psi_q)_{{\bar x}{\bar x}}\phi_{\bar x}
&=4\psi_q\biggl((\sum_{k=1}^{\nu(\alpha)} t_kq_k)^2-\sum_{k=1}^{\nu(\alpha)} t_k^2 q_k\biggr),\\
\phi_{\bar x}^{\TOP} (\psi_q)_{\bar x}&=2\psi_q \sum_{k=1}^3 t_kq_k,
\end{split}
\end{equation*}
hence (\ref{e4}) yields
\begin{equation}\label{b3}
\sum_{q}R_q \psi_q=0.
\end{equation}
Since the non-zero components  $\psi_q$ are linear independent we get the required conclusion.
\end{proof}

\begin{lemma}
\label{pr-disc}
If $f$ is a radial eigencubic of dimension $n\ge 2$ then equation (\ref{char1}) must have three distinct real roots, i.e. $\nu(\alpha)=3$. In particular, the discriminant of $\chi_\alpha$ is nonzero.
\end{lemma}

\begin{proof}
To prove the theorem we shall argue by contradiction and  assume that  $\nu(\alpha)\le 2$. This holds only if either (i) all roots $\phi_i$ are real but the discriminant of $\chi_\alpha(t)$ is zero, or (ii) $\chi_\alpha(t)$ has a  pair of conjugate complex roots.

First consider (i). We have for the discriminant (see, for example, \cite{Waerden})
$$
\mathcal{D}(\chi_\alpha)=144\alpha^2(17\alpha^2-57\alpha-243).
$$
Except for the trivial case $\alpha=0$, the discriminant vanishes only for $\alpha_{\pm}:=\frac{57\pm 39\sqrt{13}}{34}$. Since analysis of the two numbers is similar we treat only $\alpha^{+}$. For this value, (\ref{char1}) has two distinct roots
$t_1=\frac{3-6\sqrt{13}}{17}$ and $t_2=\frac{3+3\sqrt{13}}{4}$, the latter of multiplicity two. Thus $\nu(\alpha^{+})=2$ and by virtue of (\ref{e1}), $t_1n_1+t_2n_2=\alpha^{+}$. A unique integer solution of the latter equation is easily found to be $(n_1,n_2)=(1,2)$, hence, in view of (\ref{a11}), the total dimension  $n=4$. Choose $V_1=\mathrm{span}(e_1)$ and $V_2=\mathrm{span}(e_2,e_3)$ so that
\begin{equation}\label{IF}
\phi=t_1x_1^2+t_2(x_2^2+x_3^2).
\end{equation}
In order to determine $\psi$, we apply Lemma~\ref{prop1}. A direct examination of (\ref{rdef}) shows that among the $R_q$-coefficients with $q=(i,3-i)$, $0\le i \le 3$, there is only one zero coefficient, namely $R_{1,2}=0$. Thus, $\psi\equiv \psi_{12}$, i.e. $\psi$ is linear in $x_1$ and bilinear in $(x_2,x_3)$. By (\ref{e2}), $\Delta_x \psi=0$, hence $ \psi$ is congruent to the form
\begin{equation}\label{bbb}
\psi=bx_1 x_2x_3, \quad b\in \R{}.
\end{equation}
Applying the explicit form of $\phi$, we get
$$
2\phi\, \phi_{\bar x}^{\TOP}\psi_{\bar x}=4\phi b(t_1+2t_2)x_1 x_2x_3= 4b\,\alpha^{+}\,\phi\,x_1 x_2x_3,
$$
and
$$
18\alpha^{+}   {\bar x}^2\psi+\psi_{\bar x}^{\TOP} \psi_{{\bar x}{\bar x}}\psi_{\bar x}=b(18\alpha^{+}+2b^2)(x_1^2+x_2^2+x_3^2)x_1 x_2x_3.
$$
Therefore, (\ref{e6}) yields that either $b=0$ or
$$
4\alpha^{+}\,\phi=(18\alpha^{+}+2b^2)(x_1^2+x_2^2+x_3^2).
$$
The latter relation impossible because  (\ref{IF}) and  $t_1\ne t_2$. Thus $b=0$. But this implies by virtue of (\ref{bbb}) that $\psi\equiv 0$ and by (\ref{e5}), $\phi=0$, so the contradiction follows.

Now we consider the alternative (ii). This implies $\nu(\alpha)=1$ because a cubic polynomial with real coefficients must have at least one real root. We have $V\equiv V_1$, hence $n_2=n_3=0$. As a corollary, we have $\phi( \bar x)=t_1 x^2$, $x\in \R{n_1}$, where $t_1$ is the unique real root of $\chi_\alpha(t)$. In this case,  $\psi\equiv \psi_{3}$ and (\ref{b3}) reduces to a single equation
$
\psi_{3}R_{3}=0.
$
If $\psi_{3}\not\equiv0$ then in view of  (\ref{rdef}) we have $R_3(t_1)\equiv \frac{9\alpha}{2}+12t_1^2-3(3+2\alpha)t_1=0$. Therefore $R_2(t)$  and $\chi_\alpha(t)$ have a common root $t=t_1$, which implies that their resultant must be zero:
$$
\mathcal{R}(g,\chi_\alpha)\equiv -486\alpha(16\alpha+3)(\alpha-6)(\alpha+3)=0.
$$
In the cases $\alpha=-3$ and $\alpha=6$ the characteristic equation $\chi_\alpha(t)=0$ has three real roots. Thus, $\alpha=-\frac{16}{3}$, in which case the unique real root is $t_1=\frac{3}{4}$. But by virtue of (\ref{a11}), $n_1=\frac{\alpha}{t_1}=-\frac{1}{4}$, a contradiction.

It remains to consider $\psi_{3}\equiv0$. In this case we have $\psi\equiv 0$ and from (\ref{e5}) we obtain $9\alpha {\bar x}^2 +\phi_{\bar x}^2-\alpha \phi=0$, which is possible only if $\alpha\in \{0,-3,-\frac{1}{2}\}$. By the assumption, $\nu(\alpha)=1$, hence $\alpha=-\frac{1}{2}$ and the unique real root in this case is $t_1=1$. A contradiction follows because  $n_1=\frac{\alpha}{t_1}=-\frac{1}{2}$, so lemma is proved completely.
\end{proof}

Now we are ready to give a prof of the main result of this section.

\begin{proof}[Proof of Theorem~\ref{thA}]
By virtue of (\ref{reduced}),
\begin{equation}\label{HARM}
\Delta f=2x_n(\alpha+3),
\end{equation}
hence it suffices to show that for any radial eigencubic given by (\ref{reduced}) and normalized by $c=1$, there holds $\alpha=-3$. Let $f$ be an arbitrary such eigencubic.  By  Lemma~\ref{pr-disc} we have $\nu(\alpha)=3$, i.e. the characteristic polynomial (\ref{char1}) has three distinct real roots  $t_1<t_2<t_3$. Since by (\ref{char1}) $t_1t_2t_3=-9\alpha\ne 0$, we have $t_i\ne 0$.

Now we proceed by contradiction and suppose that $\alpha\ne -3$. Then $\psi\not\equiv 0$ because otherwise (\ref{e6}) would imply three alternatives $\alpha\in \{0,-3,-\frac{1}{2}\}$, of which only $\alpha=-3$ yields $\nu(\alpha)=3$, a contradiction. Let $\psi=\sum_{|q|=3} \psi_q$ be the decomposition of $\psi$ into homogeneous parts $\psi_q\in {V^*}^{\otimes q}$. We claim  that is there exists $q\ne (1,1,1)$ such that $\psi_q\not\equiv 0$. Indeed, let us suppose the contrary, i.e. that $\psi\equiv \psi_{111}$. Since $\psi\not\equiv 0$, we have $n_k=\dim V_k>0$ for $k=1,2,3$. Furthermore, $\psi\in V_1^*\otimes V_2^*\otimes V_3^*$ implies
$$
\psi_{\bar x}^2\equiv \sum_{i=1}^{n-1}\psi_{x_i}^2\in \sum_{j\ne k}{V_j^*}^{\otimes 2}\otimes {V_k^*}^{\otimes 2}=: W.
$$
We also have from the diagonal form of $\phi$
$$
\psi_{\bar x}^{\TOP} \phi_{{\bar x}{\bar x}}\psi_{\bar x}=2\sum_{i=1}^{n-1}\phi_i\psi_{x_i}^2\in W.
$$
Similarly,  $\phi_{\bar x}^{\TOP} \psi_{{\bar x}{\bar x}}\psi_{\bar x}\in  W$ because if $x_i\in V_k^*$ and $x_j\in V_l^*$ for  $k\ne l$ then
$$
\phi_{x_i}\psi_{x_ix_j}\psi_{x_j}\in {V_k^*}\otimes {V_m^*}\otimes({V_k^*}\otimes V_m^*)\subset W,
$$
where $\{k,l,m\}=\{1,2,3\}$, and if $k=l$ then $\psi=\psi_{111}$ yields $\psi_{x_ix_j}=0$. This shows that the left hand side of (\ref{e5}) belongs to $W$.

On the other hand, combining terms in the right hand side of (\ref{e5}), we get
\begin{equation*}
\begin{split}
2\sum_{i=1}^{n-1}\phi_i x_i^2\cdot \sum_{i=1}^{n-1} (4 \phi_i^2-\alpha  \phi_i+9\alpha)x_i^2
&=2\sum_{j=1}^{3}t_j u_i^2\cdot \sum_{j=1}^{3}(4 t_j^2-\alpha  t_j+9\alpha)u_i^2 \\
&=2\sum_{j=1}^{3}t_j (4 t_j^2-\alpha  t_j+9\alpha)u_i^4+h,
\end{split}
\end{equation*}
where $h\in W$ and $u_i$ is the projection of $x$ onto $V_i$. This yields $\sum_{j=1}^{3}c_iu_i^4\in W$, where $c_i=t_j (4 t_j^2-\alpha  t_j+9\alpha)$, hence  $c_1=c_2=c_3=0$. Since $t_j\ne 0$, we conclude that $4 t_j^2-\alpha  t_j+9\alpha=0$ for all $j=1,2,3$. But this yields that the quadratic polynomial $4 t^2-\alpha  t+9\alpha$ has three distinct real roots, a contradiction.
This proves that there exists $q\ne (1,1,1)$ such that $\psi_q\not\equiv0$. Applying Lemma~\ref{prop1}, we see that the corresponding $R$-coefficient must be zero. This gives
\begin{equation}\label{rpto}
0=\prod_{q\ne (1,1,1)}R_q.
\end{equation}
Write the latter product as $\rho_1\rho_2$, where
$$
\rho_{1}=R_{300}R_{030}R_{003},\qquad \rho_{2}=R_{210}R_{201}R_{120}R_{021}R_{102}R_{012}.
$$
Then $\rho_1$ and $\rho_2$ are symmetric functions of $t_i$, $i=1,2,3$, hence can be expressed as polynomials in $\alpha$. For instance, in order to find $\rho_2$ we note that
\begin{equation*}
\begin{split}
R_{210}&=2 t_2^2+6 t_1^2+4  t_1t_2-(3+2\alpha)(2t_1+t_2)+\frac{9\alpha}{2},\\
R_{201}&=2 t_3^2+6 t_1^2+4  t_1t_3-(3+2\alpha)(2t_1+t_3)+\frac{9\alpha}{2},
\end{split}
\end{equation*}
hence eliminating  $t_2$ and $t_3$ in $R_{210}R_{201}$ by virtue of Vi\`{e}te's formulas, we get
$$
R_{210}R_{201}=4( 2 t_{1}-3)^2( 12 t_{1}^2-8 \alpha  t_{1}+3\alpha)\equiv 16(  t_{1}-\frac{3}{2})^2\chi_\alpha'(t_1),
$$
which yields
\begin{equation*}\label{R120}
\begin{split}
\rho_{2}&=16^3\prod_{i=1}^3(  t_{i}-\frac{3}{2})^2 \prod_{i=1}^3\chi_\alpha'(t_i)=-4^3\chi_\alpha( 3/2)^2\,\mathcal{D}(\chi_\alpha)\equiv -2^4 3^4(\alpha+3)^2\,\mathcal{D}(\chi_\alpha),
\end{split}
\end{equation*}
where $\mathcal{D}(\chi_\alpha)=-4\chi'(\alpha, t_1)\chi'(\alpha, t_2)\chi'(\alpha, t_3)$ is the discriminant of $\chi_\alpha$.
By our assumption, the characteristic polynomial $\chi_\alpha$ has exactly three distinct real roots, hence $\mathcal{D}(\chi_\alpha)\ne 0$. Thus, in view of $\alpha\ne -3$ we conclude that $\rho_2\ne 0$. This yields by virtue of Lemma~\ref{prop1},  that $\phi_{q}\equiv 0$ for any $q$ obtained from $(1,2,0)$ by permutations. In particular,
\begin{equation}\label{psop}
\psi=\psi_{111}+\psi_{300}+\psi_{030}+\psi_{003}.
\end{equation}

On the other hand, $\rho_2\ne 0$ yields by virtue of (\ref{rpto}) that $\rho_1=0$. We have from (\ref{rdef})
\begin{equation*}\label{RQ1}
\begin{split}
0=\rho_{1}&=\prod_{k=1}^3(\frac{9\alpha}{2}+12 t_k^2-3(3+2\alpha) t_k)=12^3\prod_{k=1}^3( t_k-\frac{3}{4})( t_k-\frac{\alpha}{2}),\\
&=\frac{12^3}{4^6}\cdot \chi_\alpha(3/4) \chi_\alpha(\alpha/2)=-\frac{3^5}{2^{11}}\cdot\alpha(\alpha+3)(\alpha-6)(16\alpha+3).
\end{split}
\end{equation*}
It is easily verified that, except for $\alpha=-3$, only for $\alpha=6$ the characteristic polynomial has three real roots.
By solving the corresponding characteristic equation $\chi_{6}(t)\equiv 2(t-3)(2t^2-6t-9)=0$, we obtain $t_1=3$ and $t_{2,3}=\frac{3\pm 3\sqrt{3}}{2}$. Then (\ref{a11}) yields the relation between the dimensions $n_i=\dim V_i$ of the corresponding eigen spaces $V_i$:
$$
\frac{6n_1+3n_2+3n_3}{2}+\frac{\sqrt{3}}{2}(n_2-n_3)=6,
$$
Since $n_i$ are nonnegative integers, we find $n_2=n_3$ and $n_1+n_2=2$, which gives the following admissible triples
$$
(n_1,n_2,n_3)\in \{(2,0,0), (1,1,1), (0,2,2)\}.
$$
Substituting the found $t_i$ into (\ref{rdef}), we obtain additionally that $R_{030}$ and $R_{003}$ are non-zero, which by Lemma~\ref{prop1} and (\ref{psop}) yields
$
\psi=\psi_{111}+\psi_{300}.
$
Note that $\psi_{111}$ is harmonic because it is linear in each variable  and $\psi$ is harmonic by virtue of (\ref{e2}). Thus, $\psi_{300}$ is harmonic. On the other hand,  $\psi_{300}\not\equiv 0$ because $\psi\not \equiv \psi_{111}$. This yields $n_1\ge 2$. This strikes the triples $(1,1,1)$ and $(0,2,2)$ from the list.

Consider the only remaining triple $(n_1,n_2,n_3)=(2,0,0)$. In this case $V_2$ and $V_3$ are trivial, and $V_1$ is two-dimensional, hence $\psi\equiv \psi_{300}$. Since $t_1=3$, we have $\phi=3(x_1^2+x_2^2)$. This implies for the left hand side of (\ref{e5})
$$
2(3+\alpha)\psi_{\bar x}^2 -\psi_{\bar x}^{\TOP} \phi_{{\bar x}{\bar x}}\psi_{\bar x}-2\phi_{\bar x}^{\TOP} \psi_{{\bar x}{\bar x}}\psi_{\bar x}18\psi_{\bar x}^2-6\psi_{\bar x}^2-24\psi_{\bar x}^2=-12\psi_{\bar x}^2.
$$
On the other hand, the right hand side is strictly positive:
$$
2\phi(9\alpha {\bar x}^2 +\phi_{\bar x}^2-\alpha \phi)=2^4\cdot 3^3\, (x_1^2+x_2^2)^2.
$$
The contradiction shows that $\alpha\ne 6$, the theorem is proved completely.
\end{proof}

\section{A hidden Clifford structure}\label{sec:def}
Let us consider an arbitrary radial eigencubic $f$  given in the normal form (\ref{reduced}) normalized by $c=1$. Then by Theorem~\ref{thA}, any radial eigencubic is harmonic, hence (\ref{HARM}) yields $\alpha=-3$ in (\ref{mainproblem}) (equivalently, $\lambda(f)=-54$ in  (\ref{mainlambda})) and we have for  the characteristic polynomial (\ref{char1})
$$
\chi_{-3}(t)=4t^3+12t^2-9t-27=4(t+3)(t-\frac{3}{2})(t+\frac{3}{2}),
$$
which yields  $t_1=-3$, $t_{2}=-\frac{3}{2}$ and $t_{3}=\frac{3}{2}$. Write $\R{n}=\mathrm{span}(e_n)\oplus V$ and denote by
$
V=V_1\oplus V_2\oplus V_3
$
the eigen decomposition of $V$ associated with $\phi$. Then
\begin{equation}\label{Axi}
\phi(\bar x)=-3\xi^2-\frac{3}{2}\eta^2+\frac{3}{2}\zeta^2,\qquad \text{where $\bar x=\xi\oplus \eta\oplus \zeta\in V$}
\end{equation}
If $n_i=\dim V_i$ then in virtue of (\ref{a11})
\begin{equation}\label{a111}
\begin{split}
n_3&=2n_1+n_2-2, \qquad n=3n_1+2n_2-1.
\end{split}
\end{equation}
A close examination of (\ref{rdef}) shows  that $R_q=0$ vanish only if $q$ is one of the following: $(111), (102), (012), (030)$. This yields by virtue of Lemma ~\ref{prop1}
\begin{equation}\label{Psi}
\psi=\psi_{111}+\psi_{102}+\psi_{012}+\psi_{030},\quad \psi_q\in {V^*}^{\otimes q},
\end{equation}
and
\begin{equation}\label{Psi1}
\Delta\psi_{102}=0, \quad \Delta\psi_{012}=-\Delta\psi_{030}
\end{equation}
by virtue of (\ref{e2}).
The remaining equations (\ref{e5}) and (\ref{e6}) are read as follows:
\begin{equation}\label{e5new}
\psi_{\bar x}^{\TOP} \phi_{{\bar x}{\bar x}}\psi_{\bar x}+2\phi_{\bar x}^{\TOP} \psi_{{\bar x}{\bar x}}\psi_{\bar x}=\frac{27}{2}(\zeta^2-2\xi^2-\eta^2)(5\eta^2+3\zeta^2),
\end{equation}
\begin{equation}\label{e6new}
2\phi\, \phi_{\bar x}^{\TOP} \psi_{{\bar x}}+\psi_{\bar x}^{\TOP} \psi_{{\bar x}{\bar x}}\psi_{\bar x}=54(\xi^2+\eta^2+\zeta^2)\psi.
\end{equation}

In summary, we have

\begin{proposition}\label{pro:main}
Given a radial eigencubic $f$ in $\R{n}$, there is an orthogonal decomposition $\R{n}=\mathrm{Span}[e_n]\oplus V_1\oplus V_2 \oplus V_3,$ such that
\begin{equation}\label{given}
f=x_n^3-\frac{3}{2}x_n(2\xi^2+\eta^2-\zeta^2)+\psi_{111}+\psi_{102}+\psi_{012}+\psi_{030},
\end{equation}
where $x=(\xi,\eta,\zeta,x_n)$, $\xi\in V_1$, $\eta \in V_2$, $\zeta \in V_3$, and $\dim V_i=n_i$ satisfy $(\ref{a111})$.  Moreover,
$\psi$ satisfy $(\ref{Psi})$ and $(\ref{e5new})$--$(\ref{e6new})$. Conversely, if the cubic polynomial $(\ref{given})$ satisfies $(\ref{Psi})$ and $(\ref{e5new})$--$(\ref{e6new})$ then $f$ is a radial eigencubic.
\end{proposition}

\begin{definition*}
Suppose a radial eigencubic $f$ admits the normal form (\ref{given}). Then the pair $(n_1,n_2)$ is called the \textit{type} of the normal form, where $\dim V_1=n_1$ and $\dim V_2=n_2$.
\end{definition*}


\begin{proposition}
\label{pr:hidden}
Let $f$ be a radial eigencubic given in the normal form (\ref{given}). Then
\begin{eqnarray}
  3(\psi_{111})_\eta^2+(\psi_{102})_\zeta^2&=& 27\zeta^2\xi^2,\label{zz202}\\
   (\psi_{111})_\zeta^2 &=& 9\eta^2\xi^2, \label{zz220}\\
   (\psi_{111})_\zeta^\TOP (\psi_{102})_\zeta &=& 0,\label{zz211}\\
6(\psi_{012})_\eta^2+4(\psi_{102})_\xi^2&=&27\,\zeta^4,\label{zz004}\\
(\psi_{030})_\eta^2&=&\frac{9}{2}\,\eta^4,\label{zz040}\\
2(\psi_{111})_\eta^\TOP (\psi_{030})_\eta + (\psi_{111})_\zeta^\TOP (\psi_{012})_\zeta &=& 0,\label{zz121}\\
2(\psi_{111})_\xi^2+2(\psi_{030})_\eta^\TOP (\psi_{012})_\eta-(\psi_{012})_\zeta^2&=&-9\,\zeta^2\eta^2.\label{zz022}
\end{eqnarray}
\end{proposition}

\begin{proof}
We consider (\ref{e5new}) as an identity in ${V^*}^{\otimes q}$, $|q|=4$. Let $\pi_q$ denote the projection of ${V^*}^{\otimes 4}$ onto ${V^*}^{\otimes q}$ and let $S=\frac{27}{2}(\zeta^2-2\xi^2-\eta^2)(5\eta^2+3\zeta^2)$ denote the right hand side of (\ref{e5new}). Since
$\psi_{\bar x}^{\TOP} \phi_{{\bar x}{\bar x}}\psi_{\bar x}=-6\psi_{\xi}^2-3\psi_{\eta}^2+3\psi_{\zeta}^2
$ and
\begin{equation*}
\begin{split}
2\phi_{\bar x}^{\TOP} \psi_{{\bar x}{\bar x}}\psi_{\bar x}(\psi_{\bar x}^2)_{\bar x}^{\TOP} \phi_{{\bar x}}=\sum_{|q|=4}(-6q_1-3q_2+3q_3)\pi_q(\psi_{\bar x}^2),
\end{split}
\end{equation*}
 we obtain from (\ref{e5new})
\begin{equation}\label{beta}
\begin{split}
\pi_q(S)&=\pi_q\bigl(\psi_{\bar x}^{\TOP} \phi_{{\bar x}{\bar x}}\psi_{\bar x}+2\phi_{\bar x}^{\TOP} \psi_{{\bar x}{\bar x}}\psi_{\bar x}\bigr)\\
&=\pi_q(-6\psi_{\xi}^2-3\psi_{\eta}^2+3\psi_{\zeta}^2)-3(2q_1+q_2-q_3)\pi_q(\psi_{\bar x}^2)\\
&=-3(\beta_q+2)\pi_q(\psi_{\xi}^2)-3(\beta_q+1)\pi_q(\psi_{\eta}^2)-3(\beta_q-1)\pi_q(\psi_{\zeta}^2),
\end{split}
\end{equation}
where $\beta_q=2q_1+q_2-q_3$.

We have for $q=(202)$: $\beta_{202}=2$ and $\pi_{202}(\psi_{\xi}^2)=0$, whereas $\pi_{202}(\psi_{\eta}^2)=(\psi_{111})_\eta^2$ and $\pi_{202}(\psi_{\zeta}^2)=(\psi_{102})_\zeta^2$. This yields by (\ref{beta})
\begin{equation*}
\begin{split}
\pi_{202}(S)&=-3(\beta_{202}+1)(\psi_{111})_\eta^2-3(\beta_{202}-1)(\psi_{102})_\zeta^2\\
&=-9(\psi_{111})_\eta^2-3(\psi_{102})_\zeta^2,
\end{split}
\end{equation*}
which proves (\ref{zz202}) because $\pi_{202}(S)=-81\xi^2\zeta^2$  .

Arguing similarly for $q=(220)$ one obtains
\begin{equation*}
\begin{split}
-135\eta^2\xi^2\equiv \pi_{220}(S)&=-3(\beta_{220}-1)(\psi_{111})_\zeta^2=-15(\psi_{111})_\zeta^2,
\end{split}
\end{equation*}
hence  (\ref{zz220}) follows. The remaining  identities are established similarly.
\end{proof}

Let us rewrite $\psi_{111}$ and $\psi_{102}$ in matrix form as follows:
\begin{equation}\label{PQmatrix}
\psi_{111}=3\eta^\TOP P_\xi \zeta, \qquad \psi_{102}=\frac{3\sqrt{3}}{2}\, \zeta^\TOP Q_\xi \zeta,
\end{equation}
where
$$
P_\xi=\sum_{i=1}^{n_1}\xi_i P_i, \quad Q_\xi=\sum_{i=1}^{n_1}\xi_i Q_i, \qquad P_i\in \R{n_2\times n_3},\quad Q_i\in \R{n_3\times n_3}_{\mathrm{symm}}.
$$
Here and in what follows, by $\R{k\times m}$ we denote the vector space of matrices of the corresponding size and by $\mathbb{R}_{\mathrm{symm}}^{m\times m}$ denote the space of symmetric matrices of size $m$.
In addition to  (\ref{PQmatrix}) it is convenient also to introduce the following matrix notation:
\begin{equation}\label{PQmatrix1}
\psi_{111} =3\eta^\TOP P_\eta \zeta  \equiv3\xi^\TOP N_\eta \zeta, \qquad \psi_{012}=\frac{3\sqrt{2}}{2}\zeta^\TOP R_\eta \zeta,
\end{equation}
and
\begin{equation}\label{Heta}
H(\eta)=\frac{\sqrt{2}}{3}\nabla_\eta \psi_{030}(\eta),
\end{equation}
where  $N_\eta\in \R{n_1\times n_3}$, $R_\xi\in \R{n_3\times n_3}_{\mathrm{symm}}$.
Then the equations (\ref{zz202})--(\ref{zz022}) are rewritten in matrix notation as follows:
\begin{eqnarray}
\sum_{i=1}^{n_2}(\xi^\TOP N_i\zeta)^2+\zeta^\TOP Q_\xi^2 \zeta&=& \zeta^2\xi^2\label{mm202}\\
N_\eta N_\eta^\TOP &=& \eta^2 \mathbf{1}_{ V_1}, \label{mm220}\\
Q_\xi N_\eta^\TOP \xi &=& 0,\label{mm211}\\
\sum_{i=1}^{n_1}(\zeta^\TOP Q_i\zeta)^2+\sum_{j=1}^{n_2}(\zeta^\TOP R_j \zeta)^2&=&\zeta^4,\label{mm004}\\
H^\TOP H&=&\eta^4,\label{mm040}\\
N_H +N_\eta R_\eta &=& 0,\label{mm121}\\
2N_\eta^\TOP N_\eta+R_H-2R_\eta^2+\eta^2\mathbf{1}_{V_3}&=&0.\label{mm022}
\end{eqnarray}

\begin{proposition}[Hidden Clifford structure]\label{cor:n1n2constraint}
Let $f$  be a radial eigencubic with the normal form  $(\ref{given})$ of type $(n_1,n_2)$. Then the cubic form
$$
C(z):=\frac{2}{3}(\psi_{111}+\frac{1}{\sqrt{3}}\psi_{102}),\qquad z=(\eta,\zeta)\in \R{2n_1+2n_2-2},
$$
is an eigencubic of Clifford type.
In particular,
\begin{equation}\label{n1n2constraint}
n_1-1\le \rho(n_2+n_1-1),
\end{equation}
where $\rho$ is the Hurwitz-Radon function $(\ref{foll})$.
\end{proposition}

\begin{proof}
The case $n_1=0$ is trivial. Let us  suppose that $f$ be a radial eigencubic with the normal form $(\ref{given})$ of type $(n_1,n_2)$, where $n_1=\dim V_1\ge1 $. By using (\ref{PQmatrix}),  $(\psi_{111})_\eta=3P_\xi \zeta$, $(\psi_{111})_\zeta=3P_\xi^\TOP \eta$, and $(\psi_{102})_\zeta=6\sqrt{3}\,Q_\xi \zeta$, so that (\ref{zz202}), (\ref{zz220}) and (\ref{zz211}) become the following matrix identities:
\begin{equation}\label{nummm1}
\begin{split}
  P_\xi^\TOP P_\xi+Q_\xi ^2= \xi^2\, \mathbf{1}_{V_3}, \qquad
  P_\xi P_\xi^\TOP = \xi^2\,\mathbf{1}_{V_2},
  \qquad   P_\xi Q_\xi &= 0.
  \end{split}
\end{equation}
The latter is equivalent to that the symmetric matrices
\begin{equation}\label{Ez}
E_i=\left(
        \begin{array}{cc}
          0 & P_i \\
          P_i^\TOP & Q_i \\
        \end{array}
      \right)
\end{equation}
satisfy
$$
E_iE_j+E_jE_i=2\delta_{ij}\,1_{V_2\oplus V_3}\qquad  1\le i,j\le n_1,
$$
which implies that $\{E_i\}_{1\le i\le n_1}$ is a symmetric Clifford system in $V_2\oplus V_3=\R{n_2+n_3}=\R{2(n_1+n_2-1)}.$ In particular, this yields that $C(z)=\frac{2}{3}(\psi_{111}+\frac{1}{\sqrt{3}}\psi_{102})$   is indeed a Clifford eigencubic, cf. (\ref{Clifffff}). By using (\ref{rhommm}) we get (\ref{n1n2constraint}).
\end{proof}

The above results naturally yields a classification of radial eigencubics with $n_2=0$.

\begin{proposition}
\label{pro:n20}
A radial eigencubic  $f$ admits the normal form having the property $n_2=0$ if and only if $f$ is congruent to either the Cartan polynomial $\theta_{\ell}$, $\ell\in\{1,2,4,8\}$ or $\theta_0:=x_2^3-3x_2x_1^2$. Furthermore, in that case  $n_1=\ell+1$.
\end{proposition}

\begin{proof}
First remark by virtue of (\ref{CartanFormula0}) that for the Cartan polynomials $\theta_\ell$, $\ell=0,1,2,4,8$, there holds  $n_1=\ell+1$ and $n_2=0$.

Now suppose that $f$ is an arbitrary radial eigencubic  which admits the normal form (\ref{given}) with the property that $n_2=0$. Then $V_2=\{0\}$ and (\ref{Psi}) yields  $\psi\equiv \psi_{102}$. Moreover, (\ref{a111}) yields $n_3=2n_1-2$. If $n_1=1$  then $n_3=0$ and $\psi_{102}\equiv 0$, so that (\ref{given}) implies the trivial case, $f=x_2^3-3\xi_1^2x_2\equiv \theta_0$. Hence we can assume that $n_1\ge 2$. Then $n_3=\dim V_3\ge 1$ and
$$
f=x_n^3-\frac{3}{2}x_n(2\xi^2-\zeta^2)+\psi_{102}.
$$
In particular, by (\ref{zz004}) and (\ref{zz202}) $\psi_{\bar x}^2\equiv (\psi_{012})_\xi^2+(\psi_{102})_\zeta^2=\frac{27}{4}(\zeta^2+4\xi^2)\zeta^2,$
which implies
\begin{equation*}\label{nabla}
|\nabla  f|^2=9(x_{n}^2+\xi^2+\zeta^2)^2\equiv 9x^4.
\end{equation*}
Taking into account that $f$ is harmonic, we conclude  that $f$ satisfies the M\"untzer-Cartan equations (\ref{Mun2}), thus by the Cartan theorem $f$ must be congruent to $\theta_\ell$ for some $\ell\in \{1,2,4,8\}$. This yields  $3\ell+2=n\equiv n_1+n_3+1=3n_1-1$, hence $n_1=\ell+1$ as required.
\end{proof}

\section{Proof of Theorem~\ref{thC}}\label{sec:ifpart}

We split the proof of Theorem~\ref{thC} into two steps: the  `if'-part will be established  in Proposition~\ref{pro:exceptional} below, and the `only if'-part will be given in  Corollary~\ref{cor:reversal}.

\begin{proposition}
\label{pro:exceptional}
Let $f$ be a radial eigencubic in the normal form $(\ref{given})$. If $\psi_{030}$ is reducible and  not identically zero
 then $f$ is of Clifford type.

\end{proposition}

\begin{proof}[Proof of Proposition~\ref{pro:exceptional}]
By the assumption, $\psi_{030}$ is reducible and not identically zero, hence by the Eiconal Cubic Theorem there exist orthogonal coordinates $(\eta_1,\ldots,\eta_{n_2})$ in $V_2=\R{n_2}$ such that
$$
\psi_{030}=\frac{1}{\sqrt{2}}(\eta_{n_2}^3-3\eta_{n_2}\bar \eta^2), \qquad \bar\eta=(\eta_1,\ldots,\eta_{n_2-1})\in V_2'.
$$
Then we have for the vector field  (\ref{Heta}): $H=(\eta_{n_2}^2-\bar\eta^2, -2\eta_{n_2}\bar\eta)$. This yields
$$
N_H=-2\eta_{n_2}N_{\bar \eta}+(\eta_{n_2}^2-\bar\eta^2)N_{n_2}, \qquad
R_H=-2\eta_{n_2}R_{\bar \eta}+(\eta_{n_2}^2-\bar\eta^2)R_{n_2},
$$
where $N_\eta\in \R{n_1\times n_3}$ and $R_\eta\in \R{n_3\times n_3}_{\mathrm{symm}}$ are defined by (\ref{PQmatrix1}). Thus (\ref{mm121}) and (\ref{mm022}) becomes respectively
\begin{equation}\label{NH1}
-2\eta_{n_2}N_{\bar \eta}+(\eta_{n_2}^2-\bar\eta^2)N_{n_2}+N_\eta R_\eta=0
\end{equation}
and
\begin{equation}\label{NH10}
2N_\eta^\TOP N_\eta-2\eta_{n_2}R_{\bar \eta}+(\eta_{n_2}^2-\bar\eta^2)R_{n_2}-2R_\eta^2-\eta^2\mathbf{1}_{V_3}=0.
\end{equation}
By identifying the coefficients of $\eta_{n_2}^2$ in the latter relations one finds
\begin{equation}\label{NN100}
\begin{split}
&N_{n_2}(\mathbf{1}_{V_3}+R_{n_2})=0,\\
&2N_{n_2}^\TOP N_{n_2}=2R_{n_2}^2-R_{n_2}-\mathbf{1}_{V_3},
\end{split}
\end{equation}
which yields $(2R_{n_2}^2-R_{n_2}-\mathbf{1}_{V_3})(\mathbf{1}_{V_3}+R_{n_2})=0$. The latter equation shows that $R_{n_2}$ has eigenvalues $\pm 1$ and $-\frac{1}{2}$. Let $V_3=Y\oplus Z\oplus W$ be the corresponding eigen decomposition of $R_{n_2}$ and let $\zeta=(y,z,w)$ denote the associated decomposition of a typical vector $\zeta\in V_3$. We have
\begin{equation}\label{dim1}
\dim Y+\dim Z+\dim W=n_3\equiv 2n_1+n_2-2
\end{equation}
and
\begin{equation}\label{dim20}
\trace R_{n_2}=\dim Y-\dim Z-\frac{1}{2}\dim W.
\end{equation}
On the other hand, (\ref{mm220}) yields
\begin{equation}\label{V1V3}
N_{n_2}N_{n_2}^\TOP=\mathbf{1}_{V_1},
\end{equation}
hence we find from the second equation in (\ref{NN100})
$$
2n_1=2\trace N_{n_2} N_{n_2}^\TOP=2\trace N_{n_2}^\TOP N_{n_2}=\trace (2R_{n_2}^2-R_{n_2}-\mathbf{1}_{V_3})=2\dim Z.
$$
Thus $\dim Z=n_1$. Also, in view of the second relation in (\ref{Psi1}),
$$
\trace R_\eta=\frac{1}{3\sqrt{2} }\,\Delta_\zeta\psi_{012}= -\frac{1}{3\sqrt{2} }\,\Delta_\eta\psi_{003}=(n_2-2)\eta_{n_2},
$$
hence $\trace R_{n_2}=n_2-2$ and $\trace R_i=0$ for $1\le i\le n_2-1.$ Combining this with (\ref{dim1}) and (\ref{dim20}) we obtain
$$
\dim W=0, \qquad  \dim Y=n_1+n_2-2, \qquad \dim Z=n_1.
$$ In particular,   $V_3=Y\oplus Z$ and we can write $R_{n_2}$ in the block form
\begin{equation}\label{RN2}
R_{n_2}=\left(
        \begin{array}{cc}
          \mathbf{1}_{Y} & 0 \\
          0 & -\mathbf{1}_{Z} \\
        \end{array}
      \right),
\end{equation}
which yields in view of (\ref{NN100})
\begin{equation}\label{NNN}
N_{n_2}^\TOP N_{n_2}=2R_{n_2}^2-R_{n_2}-\mathbf{1}_{V_3}=\left(
        \begin{array}{cc}
          0 & 0 \\
          0 & \mathbf{1}_{Z} \\
        \end{array}
      \right).
\end{equation}
Regarding $N_{n_2}$ as a mapping from $V_3=Y\oplus Z$ into $V_1$, we obtain from (\ref{NNN}) the corresponding block representation: $N_{n_2}=(\mathbf{0},\, U)$, where $U:Z\to V_1$ satisfies, by virtue of (\ref{NNN}),  $U^\TOP U=1_{Z}$ and also, in view of (\ref{V1V3}), $UU^\TOP=1_{V_1}$, hence $U$ is an isometry. Thus, we may assume without loss of generality that the orthogonal coordinates in $V_1$ and $Z$ are agreed so that $U=\mathbf{1}$ is the unit matrix (of size $\dim V_1=\dim Z=n_1$). This yields $N_{n_2}=(\mathbf{0},\, \mathbf{1})$.

On substituting (\ref{RN2}) into (\ref{mm004}), we obtain
$$
\sum_{i=1}^{n_1}(\zeta^\TOP Q_i\zeta)^2+\sum_{j=1}^{n_2-1}(\zeta^\TOP R_j \zeta)^2=(y^2+z^2)^2-(y^2-z^2)^2=4y^2z^2,
$$
where
\begin{equation}\label{ZetA}
\zeta=(y,z)\in Y\oplus Z,
\end{equation}
which implies
\begin{equation}\label{QR}
Q_j=\left(
        \begin{array}{cc}
          0 & \alpha_j \\
          \alpha_j^\TOP & 0 \\
        \end{array}
      \right),
      \qquad
      R_i=\left(
        \begin{array}{cc}
          0 & \beta_i \\
          \beta_i^\TOP & 0 \\
        \end{array}
      \right), \quad 1\le i\le n_2-1.
\end{equation}

Coming back to (\ref{mm220}), it can be seen that
\begin{equation}\label{inff}
N_{n_2}N_{\bar\eta}^\TOP+N_{\bar\eta} N_{n_2}^\TOP={0}  \quad\text{and}\quad
N_{\bar\eta}N_{\bar\eta}^\TOP=\bar\eta^2\mathbf{1}_{V_1}.
\end{equation}
Regarding $N_{\bar\eta}\equiv \sum_{i=1}^{n_2-1}N_i\eta_i$ as a mapping from $V_3=Y\oplus Z$ into $V_1$ and rewriting it in the block form as $(a_{\bar\eta}, b_{\bar\eta})$, we deduce from (\ref{inff}) that
\begin{equation*}
b_{\bar\eta}+b_{\bar\eta}^\TOP={0}, \qquad a_{\bar\eta}a_{\bar\eta}^\TOP+b_{\bar\eta}b_{\bar\eta}^\TOP=\bar\eta^2\mathbf{1}_{V_1}.
\end{equation*}
Also, by identifying the coefficients of $\eta_{n_2}$ in (\ref{NH1}), one finds
$2N_{\bar \eta}=N_{n_2} R_{\bar\eta}+N_{\bar\eta} R_{n_2},$
which  yields
$a_{\bar\eta}=\beta_{\bar\eta}^\TOP$ and  $b_{\bar\eta}=0$.
Thus,
\begin{equation}\label{Neta}
N_\eta=(\beta^\TOP_{\bar\eta}, \, \eta_{n_2}\mathbf{1})
\end{equation}

Now we consider the normal form (\ref{given}) and rewrite it in the Clifford form (\ref{fx0}). To this end, let us first introduce the new orthogonal coordinates
$$
x_n=\frac{x_0+\sqrt{2}y_0}{\sqrt{3}}, \quad \eta_{n_2}=\frac{\sqrt{2}x_0-y_0}{\sqrt{3}}
$$
in the $(x_n, \eta_{n_2})$-plane. Then  the following identity is verified by a straightforward calculation:
$$
x_n^3+\frac{3}{2}x_n(-\eta_{n_2}^2-\bar\eta^2+y^2)+\frac{\sqrt{2}}{2}(\eta_{n_2}^3-3\eta_{n_2}\bar\eta^2)+\frac{3\sqrt{2}}{2}\eta_{n_2}y^2
=\frac{3\sqrt{3}}{2}x_0(y_0^2+y^2-\bar\eta^2).
$$
Taking into account (\ref{ZetA}), (\ref{Neta}) and (\ref{QR}), we rewrite (\ref{given}) in the new coordinates as folows
\begin{equation}\label{fbecome}
f=\frac{3\sqrt{3}}{2}\,x_0(y_0^2+y^2-\bar\eta^2)+\frac{3x_n}{2}(z^2-2\xi^2)-\frac{3\sqrt{3}}{2}\,\eta_{n_2}z^2+3\eta_{n_2} \,\xi^\TOP z+\Omega,
\end{equation}
where
$$
\Omega=3y^\TOP \beta_{\bar\eta}(\xi+z\sqrt{2})+3\sqrt{3}y^\TOP\alpha_{\xi}z.
$$
Let us introduce the new orthogonal coordinates in $V_1\oplus Z$ by virtue of
\begin{equation*}\label{xitea}
s=\frac{\xi+\sqrt{2}z}{\sqrt{3}}, \qquad t=\frac{z-\sqrt{2}\xi}{\sqrt{3}}.
\end{equation*}
Then
$$
\frac{3x_n}{2}(z^2-2\xi^2)-\frac{3\sqrt{3}}{2}\,\eta_{n_2}z^2+3\eta_{n_2} \,\xi^\TOP z=-\frac{3\sqrt{3}}{2}(\,x_0t^2-2y_0\, s^{\TOP}t),
$$
and (\ref{fbecome}) becomes
\begin{equation}\label{fbecome1}
f_1=x_0(y_0^2+y^2-\bar\eta^2-t^2)+2y_0\, s^{\TOP}t+\Omega,
\end{equation}
where $f_1:=\frac{2\sqrt{3}}{9}f$ and
\begin{equation}\label{OMEGA}
\Omega=2y^\TOP \beta_{\bar\eta}s+2y^\TOP\alpha_{\xi}z.
\end{equation}
Next, note that in view of (\ref{QR}) and  (\ref{mm211})
\begin{equation}\label{NQ}
\begin{split}
\alpha_\xi \xi&=0,\\
\end{split}
\end{equation}
hence setting $\xi=\xi'+\xi''$ in (\ref{NQ})  we get $\alpha_{\xi'}\xi''=-\alpha_{\xi''}\xi'$, which yields for the last term in (\ref{OMEGA})
$$
y^\TOP \alpha_{\xi}z=\frac{1}{3}y^\TOP \alpha_{s-\sqrt{2}t}(t+\sqrt{2}s)=y^\TOP \alpha_{s}t.
$$
Thus $\Omega=2y^\TOP \beta_{\bar\eta}s+2y^\TOP \alpha_{s}t,$
hence we arrive at the following expression
\begin{equation}\label{fbecome2}
f_1=x_0(y_0^2+y^2-\bar\eta^2-t^2)+2y_0\, s^{\TOP}t+2y^\TOP \beta_{\bar\eta}s+2y^\TOP \alpha_{s}t.
\end{equation}
In order to show that the last expression is indeed a Clifford representation, we rewrite $f_1$ in matrix notation as follows. We combine the coordinates as follows:
$$
{\tilde y}=(y_0,y)\in\widetilde Y,\qquad
{\tilde z}=(\bar\eta,t)\in\widetilde Z,
$$
where $Y\cong Z\cong\R{n_1+n_2-1}$, and set  $U=\widetilde Y\oplus \widetilde Z $ and $u=(\tilde y,\, \tilde z)\in U$. Then (\ref{fbecome2}) becomes
\begin{equation}\label{fbecome3}
f_1=x_0u^\TOP A_0 u+\sum_{i=1}^{n_1}s_i\,u^\TOP A_i u\equiv x_0u^\TOP A_0 u + u^\TOP A_s u
\end{equation}
where the matrices
$$
A_0=\left(
        \begin{array}{cc}
          \mathbf{1}_{\widetilde Y} & 0 \\
          0 &  -\mathbf{1}_{\widetilde Z}\\
        \end{array}
      \right),
\qquad A_i=\left(
        \begin{array}{cc}
          0 & D_i \\
          D_i^\TOP & 0 \\
        \end{array}
      \right),
$$
are written in the block form with respect to the orthogonal decomposition $U=\widetilde Y\oplus \widetilde Z$. Here the matrix   $D_s=\sum_{i=1}^{n_1}s_iD_i$ is defined by virtue of
$$
2y_0\, s^{\TOP}t+2y^\TOP \beta_{\bar\eta}s+2y^\TOP \alpha_{s}t=\tilde y^\TOP D_s \tilde z,
$$
in other words,
$$
D_s=\left(
        \begin{array}{cc}
          0 & s^\TOP \\
          G_s & \alpha_s \\
        \end{array}
      \right),
$$
where the latter block-form is associated with the vector decompositions ${\tilde y}=(y_0,y)$ and ${\tilde z}=(\bar\eta,t)$, and the matrix $G_s$ is determined by dualizing
\begin{equation}\label{dualz}
y^\TOP \beta_{\bar\eta} s=y^\TOP G_{s}\bar\eta.
\end{equation}

Coming back to (\ref{fbecome3}), we see that it suffices to show that $\{A_0, A_1,\ldots, A_{n_1}\}$ is a symmetric Clifford system in $\R{2(n_1+n_2-1)}$. This will be done, in view of the explicit form of $A_0$ and the equality $\dim \widetilde Y=\dim \widetilde Z$, if we show that the subsystem $\{A_1,\ldots, A_{n_1}\}$ is also a symmetric Clifford system, which is, in its turn, is equivalent to the following two identities: $D_sD_s^\TOP=s^2 \mathbf{1}_{\widetilde Y}$ and $D_s^\TOP D_s=s^2\mathbf{1}_{\widetilde Z}$. Since the matrix $D_s$ is square, it suffices to prove, e.g., the first of the last two relations. To this end, we write by virtue of (\ref{NQ})
$$
D_sD_s^\TOP=\left(
        \begin{array}{cc}
          s^2 & s^\TOP \alpha_s^\TOP \\
          \alpha_s s & G_s G_s^\TOP+\alpha_s\alpha_s^\TOP \\
        \end{array}
      \right)
=\left(
        \begin{array}{cc}
          s^2  & 0 \\
          0 & G_s G_s^\TOP+\alpha_s\alpha_s^\TOP \\
        \end{array}
      \right).
$$
From (\ref{dualz}) we have
\begin{equation}\label{sy}
y^\TOP G_s G_s^\TOP y=|G_s^\TOP y|^2=\sum_{i=1}^{n_2-1}(y^\TOP G_s e_i)^2=\sum_{i=1}^{n_2-1}(y^\TOP \beta_i s)^2,
\end{equation}
where $\{e_i\}$ is an orthonormal basis in $V_2'=\{(\bar \eta,0)\in V_2\}$.
On the other hand, by rewriting (\ref{mm202}) by virtue of (\ref{Neta}), we obtain
$$
\xi^\TOP z+\sum_{i=1}^{n_2-1}(y^\TOP \beta_i \xi)^2+y^\TOP \alpha_\xi\alpha_\xi^\TOP y+z^\TOP \alpha_\xi^\TOP\alpha_\xi z= (y^2+z^2)\xi^2,
$$
and setting  $z=0$ and $\xi=s$ in the latter identity we get in view of (\ref{sy})
$$
y^\TOP G_s G_s^\TOP y+y^\TOP \alpha_s\alpha_s^\TOP y= y^2s^2,
$$ i.e.
$G_s G_s^\TOP+\alpha_s\alpha_s^\TOP=s^2\mathbf{1}_{\widetilde Y}$. This yields $D_sD_s^\TOP=s^2 \mathbf{1}_{\widetilde Y}$ as required. Thus $f_1$ is a Clifford eigencubic and the proposition is proved completely.

\end{proof}


\begin{proposition}
\label{pro:finiten2}
If $f$ is an exceptional eigencubic in the normal form (\ref{given}) then $\psi_{030}$ is either reducible or identically zero. Furthermore for an exceptional eigencubic $n_2\in \{0,5,8,14,26\}$ and the triple $(n_1,n_2,n_3)$ can take only the values presented in  Table~\ref{tabs}.
\end{proposition}

\begin{small}
\def\MM{6mm}
\begin{table}[ht]
\renewcommand\arraystretch{1.5}
\noindent
\begin{flushright}
\begin{tabular}{|p{4mm}||p{2.8mm}|p{2.8mm}|p{2.8mm}|p{2.8mm}||
p{2.8mm}|p{2.8mm}|p{2.8mm}|p{2.8mm}||p{2.8mm}|p{2.8mm}|p{2.8mm}|p{2.8mm}|p{2.8mm}|p{2.8mm}||p{2.8mm}|p{2.8mm}|p{2.8mm}|p{2.8mm}||p{2.8mm}|
p{2.8mm}|p{2.8mm}|p{2.8mm}|p{2.8mm}|p{2.8mm}|p{2.8mm}|p{2.8mm}|}
\hline
$n_1$& $2$ & $3$  & $5$  & $9$  & $0$ & $1$  & $2$  & $4$  & $0$  & $1$  &  $2$ & $3$ &  $5$ & $9$  & $0$  & $1$   &  $2$  & $3$   & $0$  & $1$  & $2$ & $3$ & $7$ \\\hline
$n_2$&  $0$ & $0$  & $0$  & $0$  &$5$ & $5$  & $5$  & $5$  & $8$  & $8$  &  $8$ &  $8$ &  $8$ & $8$  &$14$  & $14$  & $14$  & $14$  & $26$ & $26$ & $26$ & $26$ & $26$      \\\hline
$n_3$&  $2$ & $4$  & $8$  & $16$  &$3$ & $5$  & $7$  & $11$ & $6$  & $8$  & $10$ & $12$ & $16$ & $24$ &$12$  & $14$  & $16$  & $18$  & $24$ & $26$ & $28$ & $30$ & $38$      \\\hline
$n$  &  $5$ & $8$  & $14$  & $26$  &$9$ & $12$ & $15$ & $21$ & $15$ & $18$ & $21$ & $24$ &  $30$&  $42$& $27$ & $30$  & $33$  & $36$  & $51$ & $54$ & $57$ & $60$ & $72$       \\\hline
\end{tabular}
\end{flushright}
\bigskip
\caption{Possible type of exceptional eigencubics}\label{tabs}
\end{table}
\end{small}

\begin{proof}
Let $f$ be an exceptional eigencubic given in the normal form (\ref{given}). Then by Proposition~\ref{pro:exceptional}, $\psi_{030}$ is either irreducible or identically zero. This implies by virtue of the Eiconal Cubic Theorem  that $n_2\in \{0,5,8,14,26\}$. On the other hand, by Proposition~\ref{cor:n1n2constraint} we have
\begin{equation}\label{n111}
n_1-1\le \rho(n_2+n_1-1).
\end{equation}
Note that the following elementary inequality for the Hurwitz-Radon function holds:
\begin{equation}\label{rhomm}
\rho(m)\le \frac{1}{2}(m+8), \qquad m\ge1.
\end{equation}
Indeed, by the definition (\ref{foll}) $\rho(m)=8a+2^b$, where $m=2^sn_1$ with $n_1$ an odd number and $s=4a+b$, $b=0,1,2,3$. For these $b$, $2^b\le 2b+2$, hence
$
\rho(m)\le 8a+2b+2.
$
On the other hand,
$$
8a+2b+2=2s+2\le 2^{s-1}+4\le \frac{1}{2}(2^sn_1+8),
$$
which proves (\ref{rhomm}).

By (\ref{n111}) and (\ref{rhomm}), $n_1-1\le 2n_2+16$. Thus, given $n_2\in\{0,5,8,14,26\}$ it suffices to examine (\ref{n111}) only for the numbers of $n_1$ for which $n_1\le 2n_2+17$. This easily yields the numbers presented in Table~\ref{tabs}.
\end{proof}

\section{Proof of Theorem~\ref{thB}}\label{sec:trace}

We shall establish the cubic trace formula for radial eigencubic of Clifford type and exceptional eigencubics, separately in Proposition~\ref{pro:norm} and Proposition~\ref{pro:tau} below. We start with some definitions and lemma.

\begin{definition*}
We  say that an arbitrary cubic polynomial $f$ in $\R{n}$ possess the \textit{quadratic trace identity} (with the constant $\beta\in \R{}$) if \begin{equation}\label{someQ}
\trace \mathrm{Hess}^2 f=\beta x^2,\qquad \beta=\beta(f)\in \R{}.
\end{equation}
Similarly, we say $f$ possess the \textit{cubic trace identity}  if it satisfies
\begin{equation}\label{somec}
\trace \mathrm{Hess}^3 f=\alpha f,
\end{equation}
for some $\alpha\in \R{}$.

\end{definition*}

We shall also make use the following quadratic form
\begin{equation}\label{ssigg}
\sigma_2(f)=-\frac{\trace \mathrm{Hess}^2(f)}{\lambda}, \qquad \text{where}\,\, \lambda=\frac{L(f)}{x^2f},
\end{equation}
which spectrum is a congruence invariant of $f$  in view of the invariant properties of the operator $L$.

\begin{lemma}
\label{pro:trace2}
For any  radial eigencubic  $f$ in $\R{n}$ the following holds.
\begin{itemize}
\item[(i)]
If $f$ possess the quadratic trace identity  with $\beta$ then $f$ posses the cubic trace identity with the constant $\alpha=- (n+6)\lambda-6\beta$.
\item[(ii)]
If $f$ possess the cubic trace identity with  $\alpha$ then any normal form of $f$ has type  $(n_1,n_2)$, where $n_1=\frac{\alpha}{3\lambda}+1$ and $n_2=\frac{n-3n_1+1}{2}$.
\end{itemize}
In particular, if a radial eigencubic $f$ possesses the cubic trace identity then its type is uniquely determined by the congruence class of $f$.
\end{lemma}

\begin{proof}
(i) Since $f$ is harmonic, we find from (\ref{mainlambda}) that $\sum_{i,j=1}^{n}f_{x_i}f_{x_ix_j}f_{x_j}=-\lambda x^2f$, hence applying the Laplacian to the latter identity and using the homogeneity of $f$ we obtain
$$
2\trace \mathrm{Hess}^3 f+4\sum_{i,j,k=1}^nf_{x_ix_j}f_{x_ix_j x_k}f_{x_k}=-\lambda (2n+12)f,
$$
On the other hand, by our assumption $\trace \mathrm{Hess}^2 f=\beta x^2$, hence
$$
\sum_{i,j,k=1}^nf_{x_ix_j}f_{x_ix_j x_k}f_{x_k}\equiv \frac{1}{2}\scal{\nabla f}{\nabla \trace \mathrm{Hess}^2 f}=\beta \sum_{k=1}^n x_kf_{x_k}=3\beta f,
$$
which yields $\trace \mathrm{Hess}^3 f=-((n+6)\lambda+6\beta)f$ and thereby proves the first claim of the proposition.

(ii) Now suppose that $f$ possess the cubic trace identity (\ref{somec}). It follows from the invariant properties of the operator $L$, see for instance \cite{Hsiang67}, that the ratio
\begin{equation}\label{alphalambda}
f\to \frac{\alpha}{\lambda}\equiv \frac{\trace \mathrm{Hess}^3 f}{f\cdot \lambda}\in \R{}
\end{equation}
is invariant under orthogonal substitutions and dilatations, hence it suffices to establish (ii) in assumption  that $f$ is given in the normal form, e.g. by virtue of (\ref{reduced}) and normalized by $c=1$. In that case, the Hessian matrix of $f$ has the following block form associated with the orthogonal decomposition $\R{n}=\mathrm{span}(e_n)\oplus \mathrm{span}(e_n)^\bot$:
$$
\mathrm{Hess} f=\left(
        \begin{array}{ll}
          6x_n\,\, & \phi_x^{\TOP} \\
          \phi_x\,\, & x_n\phi_{xx}+\psi_{xx} \\
        \end{array}
      \right),
$$
hence
\begin{equation*}
\begin{split}
\trace \mathrm{Hess}^3 f&=216 x_n^3+18x_n\phi_x^2+3\trace \phi_x\phi_x^\TOP(x_n\phi_{xx}+\psi_{xx})+\trace (x_n\phi_{xx}+\psi_{xx})^3\\
&=(216+\trace \phi_{xx}^3)x_n^3+\ldots
\end{split}
\end{equation*}
where the dots stands for the  degrees of $x_n$ lower than 3.
On the other hand,  the coefficient of $x_n^3$ in $\trace \mathrm{Hess}^3 f$  can be read out from the cubic trace identity $\trace \mathrm{Hess}^3 f=\alpha f$ which yields $\alpha=216+\trace \phi_{xx}^3$. Using (\ref{Axi}) and (\ref{a111}),
$$
\alpha=-216+(-216n_1-27n_2+27n_3)=162(1-n_1)\equiv 3\lambda (n_1-1),
$$
which yields $n_1=\frac{\alpha}{3\lambda}+1$, as required.

Finally, note that the dimension $n_1$ is uniquely determined from the cubic trace identity by virtue of $\lambda$ and $\alpha$, and in view of the remark made above the ratio $\frac{\alpha}{3\lambda}$ is a congruence invariant. Since $n_2$ is determined by $n_2=\frac{n+1-3n_1}{2}$,we conclude  that it is also a congruence invariant. The proposition is proved completely.
\end{proof}

\subsection{The cubic trace identity in the Clifford case}

\begin{proposition}\label{pro:tau}
Let $\mathcal{A}=(A_0,\ldots, A_q)\in \mathrm{Cliff}(\R{2m},q)$ and
\begin{equation}\label{Clifffff0}
C_\mathcal{A}(x)=\sum_{i=0}^{q}  \scal{y}{A_iy}\,x_{i+1}, \qquad  y=(x_{q+2},\ldots,x_{q+1+2m})\in \R{2m}.
\end{equation}
be the Clifford eigencubic  associated with $\mathcal{A}$. Then
\begin{equation}\label{taU2}
\begin{split}
\sigma_2(C_\mathcal{A})&=mz^2+(q+1)y^2,\qquad z=(x_1,\ldots,x_{1+q}),
\end{split}
\end{equation}
and $C_\mathcal{A}$ possess the cubic trace identity with $\alpha=3(q-1)\lambda$. In particular,  $C_\mathcal{A}$ has the type
\begin{equation}\label{Cliffordtype}
(n_1,n_2)=(q,m+1-q).
\end{equation}

\end{proposition}

\begin{proof}
Write the Hessian matrix of $f\equiv C_\mathcal{A}$ in the block form
$$
\mathrm{Hess} f\equiv
\left(
        \begin{array}{ll}
          f_{yy} & f_{yz} \\
          f_{zy} & f_{zz}\\
        \end{array}
      \right)
=\left(
        \begin{array}{ll}
          2A_z\,\, & B \\
          B^\TOP& 0 \\
        \end{array}
      \right),
$$
where $B$ is the matrix with entries $B_{ij}= f_{y_iz_j}=2e_i^{\TOP}A_j y$ and $\{e_i\}_{i=1}^{2m}$ is the standard basis in $\R{2m}$.
Then
\begin{equation*}\label{trt1}
\begin{split}
\trace \mathrm{Hess}^2 f&=4\trace A_z^2+2\trace BB^{\TOP}=8mz^2+8\sum_{i=1}^{2m}\sum_{j=0}^{q}y^\TOP A_j e_ie_i^{\TOP}A_j y\\
&=8(mz^2+(q+1)y^2).
\end{split}
\end{equation*}
By Theorem~3.2 in \cite{TkCliff}  we have $\lambda(C_\mathcal{A}) =-8$ which yields  (\ref{taU2}).

In order to establish  the cubic trace identity, we note that $A_z^2=z^2 \mathbf{1}_{\R{2m}}$ and  $\trace A_i=0$, so that $\trace A_z^3=\trace z^2A_z=0$ and
\begin{equation}\label{combb}
\trace \mathrm{Hess}^3 f=8\trace A_z^3+6\trace A_zBB^{\TOP}=6\trace A_zBB^{\TOP}.
\end{equation}
From (\ref{Cliffsys}) $A_zA_k+A_kA_z=2z_k\mathbf{1}_{\R{2m}}$, hence
\begin{equation}\label{combb1}
\begin{split}
\trace A_zBB^{\TOP}&=\sum_{i,j=1}^{2m}(A_z)_{ij}(BB^{\TOP})_{ij}=4\sum_{k=0}^{q}\sum_{i,j=1}^{2m}(A_z)_{ij}e_i^{\TOP}A_k y\cdot e_j^{\TOP}A_k y\\
&
=4\sum_{k=0}^{q}y^{\TOP}A_k(\sum_{i,j=1}^{2m}(A_z)_{ij}\,e_i\cdot e_j^{\TOP})A_ky=4\sum_{k=0}^{q}y^{\TOP}A_kA_zA_ky\\
&=4\sum_{k=0}^{q}y^{\TOP}(2z_k\mathbf{1}_{\R{2m}}-A_zA_k)A_ky\\
&=4(1-q)y^{\TOP}A_zy.
\end{split}
\end{equation}
Combining (\ref{combb}) and (\ref{combb1}) yields
$$
\trace \mathrm{Hess}^3 f=24(1-q)f,
$$
hence $f$ possess the cubic trace identity with $\alpha=24(1-q)\equiv 3(q-1)\lambda$. Then by using  (ii) in Lemma~\ref{pro:trace2} we find from $n=3n_1+2n_2-1=q+1+2m$ that $n_1=q$ and $n_2=m+1-q$.

\end{proof}

\begin{corollary}\label{cor:qm1}
A pair $(n_1,n_2)$ of non-negative integers is the type of some eigencubic of Clifford type if and only if
\begin{equation}\label{inequality}
n_1\le \rho(n_2+n_1-1).
\end{equation}

\end{corollary}

\begin{proof}
Let us first suppose that  (\ref{inequality}) holds. We may assume without loss of generality that $n_1\ge0$ and $n_2+n_1-1\ge 1$. Then setting $q=n_1$ and $m=n_1+n_2-1$, the inequality  $n_1-1< \rho(n_2+n_1-1)$ becomes equivalent to  $q\le \rho(m)$ which implies that there exists a symmetric Clifford system $\mathcal{A}=\{A_0,\ldots, A_q\}\in \mathrm{Cliff}(\R{2m},q)$. Let $C_\mathcal{A}$ the associated with $\mathcal{A}$  Clifford eigencubic  (\ref{Clifffff}). Then  Proposition~\ref{pro:tau} yields that  $C_\mathcal{A}$ is a Clifford eigencubic of type $(q,m+1-q)\equiv (n_1,n_2)$.

In the converse direction, let us assume  there exists a Clifford eigencubic $f$ of the type $(n_1,n_2)$. Since the case is trivial, we may assume that $n_1\ge1$. By Proposition~\ref{cor:n1n2constraint} we have $n_1-1\le  \rho(n_2+n_1-1)$, hence it suffices to show that the equality in the latter inequality is impossible. We assume the contrary, i.e. that there exists a Clifford eigencubic $f$ of the type $(n_1,n_2)$ with $n_1\ge 1$ and $n_1-1= \rho(n_2+n_1-1)$. By Proposition~\ref{pro:tau}, $f$  is of type $(q,m+1-q)=(n_1,n_2)$, hence $n_1=q$ and $n_2=m+1-q$. On the other hand, the existence of the Clifford eigencubic $f$ implies $q\le \rho(m)$, hence $n_1\le \rho(n_2+n_1-1)$, a contradiction.
\end{proof}

\subsection{The trace identities for exceptional eigencubics}

\begin{proposition}
\label{pro:norm}
Let $f$ be an arbitrary radial eigencubic having the normal form of type $(n_1,n_2)$, $n_2\in \{0,5,8,14,26\}$,  and such that  $\Delta \psi_{030}=0$. Then
\begin{equation}\label{eq:norm}
\trace \mathrm{Hess}^2(f)=-\frac{1}{3}(3n_1+n_2+1)\,\lambda x^2,\end{equation}
\begin{equation}\label{trace_iden}
\trace \mathrm{Hess}^3(f)=3(n_1-1)\, \lambda f,
\end{equation}
where $L(f)=\lambda x^2 f$. In particular, the trace identities $(\ref{eq:norm})$ and $(\ref{trace_iden})$ are valid for  any exceptional eigencubic.
\end{proposition}

\begin{remark}
Observe, that the statement of Proposition~\ref{pro:norm} establishes also the implications (a)$\Rightarrow$(c) and (b)$\Rightarrow$(c) in Theorem~\ref{thD}.
\end{remark}

\begin{proof}
Observe that it suffices  to prove the first identity. Indeed, if $(\ref{eq:norm})$ holds then by Lemma~\ref{pro:trace2} we have   $\trace \mathrm{Hess}^3(f)=\alpha f$ with
$$
\alpha=-(n+6)\lambda+2(3n_1+n_2+1)\lambda=3(1-n_1)\lambda,
$$
so that (\ref{trace_iden}) follows.

We split the proof of (\ref{eq:norm}) into two steps. First we suppose that $n_2=0$. Then by Proposition~\ref{pro:n20} $f$ is an isoparametric eigencubic of type $(\ell+1,0)$, $\ell=1,2,4,8$. In particular, $f$ is harmonic and  satisfies the eiconal equation $|\nabla f|^2=cx^4$ for some $c>0$. Hence
$$
L(f)=|\nabla f|^2\Delta f-\frac{1}{2}\scal{\nabla |\nabla f|^2}{\nabla f}=-6cfx^2,
$$
hence $\lambda(f)=-6c$. On the other hand, by using the harmonicity of $f$ again, we get
\begin{equation*}\label{Heq}
\trace \mathrm{Hess}^2(f)=\sum_{i,j=1}^n f_{x_ix_j}^2\equiv \sum_{i=1}^{n}f_{x_i}\Delta f_{x_i}+\sum_{i,j=1}^n f_{x_ix_j}^2=\frac{1}{2}\Delta |\nabla f|^2=2c(n+2)x^2.
\end{equation*}
Since $n=3\ell+2\equiv 3n_1+n_2-1$, we latter identity is equivalent to (\ref{eq:norm}).

Now suppose that $n_2=3\ell +2$, $\ell\in \{1,2,4,8\}$, and $f$ is written in the normal form (\ref{given}) with $\Delta \psi_{030}=0$.
Let $\R{n}=V_0\oplus V_1\oplus V_2\oplus V_3$, $V_0=\mathrm{span}(e_n)$, be the associated with (\ref{given}) orthogonal decomposition. Then
$$
\trace \mathrm{Hess}^2 f=\sum_{i,j=0}^3\trace f_{V_iV_j}f_{V_jV_i}\equiv \sum_{i,j=0}^3 T_{ij}, \qquad T_{ij}=T_{ji}.
$$
Here $f_{V_iV_j}$ stands for the submatrix of the Hessian of $f$ with entries $f_{u,v}$, where $u$ and $v$ run orthogonal coordinates in $V_i$ and $V_j$ respectively. We have $f_{x_nx_n}=6x_n$, $f_{x_n\xi}=-6\xi^\TOP$, $f_{x_n\eta}=-3\eta^\TOP$, $f_{x_n\zeta}=3\zeta^\TOP$, $f_{\xi\xi}=-6x_n\mathbf{1}_{V_1}$, $f_{\xi_i\eta_j}=3e_i^\TOP N_j \zeta$, where $\{e_i\}$ is an orthonormal basis in $V_1$, and the matrices $N_i$ are defined by (\ref{PQmatrix1}). This yields
\begin{equation*}
\begin{split}
T_{00}&=36x_n^2,\qquad T_{01}=36\xi^2,\qquad T_{03}=9\eta^2,\\
T_{04}&=9\zeta^2,\qquad T_{11}=36n_1x_n^2.
\end{split}
\end{equation*}
We also have
$$
T_{12}= 9\sum_{i=1}^{n_1}\sum_{j=1}^{n_2}\zeta^\TOP N_j^\TOP e_i\cdot e_i^\TOP N_j \zeta=9\zeta^\TOP (\sum_{j=1}^{n_2}N_j^\TOP N_j) \zeta.
$$
On the other hand, since $\trace N_\eta N_\eta^\TOP=n_1\eta^2$ by virtue of (\ref{mm220}) and $\sum_{i=1}^{n_1}Q_i N_\eta^\TOP e_i=\frac{1}{2}\Delta_\xi (Q_\xi N_\eta^\TOP \xi)= 0$ by virtue  of (\ref{mm211}), we find
\begin{equation*}
\begin{split}
T_{13}&=9\sum_{i=1}^{n_1}\sum_{k=1}^{n_3} (e_i^\TOP N_\eta\epsilon_k+\sqrt{3}\zeta^\TOP Q_i \epsilon_k)^2=9\trace N_\eta N_\eta^\TOP+18\sqrt{3}\zeta^\TOP \sum_{i=1}^{n_1}Q_i N_\eta^\TOP e_i+27\zeta^\TOP (\sum_{i=1}^{n_1}Q_i^2)\zeta\\
&=9n_1\eta^2+27\zeta^\TOP (\sum_{i=1}^{n_1}Q_i^2)\zeta,
\end{split}
\end{equation*}
where $\{\epsilon_k\}$ is an orthonromal basis in $V_3$.

By the above, $\psi_{030}$ is harmonic, hence (\ref{zz040}) yields
$$
\trace ((\psi_{030})_{\eta\eta})^2=\frac{1}{2}\Delta_\eta(\psi_{030})_\eta^2=\frac{9}{4}\Delta_{\eta}\eta^4=9(n_2+2)\eta^2,
$$
hence
\begin{equation*}
\begin{split}
T_{22}&=\trace ((\psi_{030})_{\eta\eta}-3x_n\mathbf{1}_{V_2})^2=9(n_2+2)\eta^2+9n_2x_n^2.
\end{split}
\end{equation*}

Next, $f_{\eta_j\zeta_k}=3(\xi N_j+\sqrt{2}\zeta^\TOP R_j) \epsilon_k$ yields
\begin{equation*}
\begin{split}
T_{23}&=9\sum_{j=1}^{n_2} |N_j^\TOP\xi+\sqrt{2}R_j \zeta|^2=9\sum_{j=1}^{n_2}(\xi^\TOP N_iN_i^\TOP \xi
+2\sqrt{2} \xi^\TOP N_iR_j \zeta+2\zeta^\TOP R_i^2\zeta).
\end{split}
\end{equation*}
From (\ref{mm220}) $N_iN_i^\TOP=\mathbf{1}_{V_1}$. On the other hand, $\Delta_\eta H_i=\frac{\sqrt{2}}{3}\partial_{\eta_i}\Delta_\eta \psi_{030}=0$, hence (\ref{mm121}) yields $0=\Delta_\eta(N_H +N_\eta R_\eta) = \sum_{j=1}^{n_2}N_jR_j$, thus, $$T_{23}=9n_2\xi^2+18\sum_{j=1}^{n_2}\zeta^\TOP R_j^2\zeta.$$

Finally, $f_{\zeta_i\zeta_j}=3x_0\delta_{ij}+3\sqrt{3}Q_\xi+3\sqrt{2}R_\eta$ yields
$$
T_{33}=9n_3x_0^2+18\sqrt{3}\trace Q_\xi+18\sqrt{2}\trace R_\eta+27\trace Q_\xi^2+18\sqrt{6}\trace Q_\xi R_\eta +18\trace R_\eta^2.
$$
From (\ref{nummm1})
$$
n_3 \xi^2=\trace (P_\xi^\TOP P_\xi+Q_\xi ^2)\trace P_\xi P_\xi^\TOP+\trace Q_\xi ^2=n_2\xi^2+\trace Q_\xi ^2,
$$
hence $\trace Q_\xi ^2=(n_3-n_2)\xi^2=2(n_1-1)\xi^2$. Similarly, the harmonicity of $\psi_{030}$ yields from (\ref{Psi1}) $\trace R_i=\trace Q_j=0$, hence taking the trace in (\ref{mm022})
$$
0=2\trace N_\eta^\TOP N_\eta-2\trace R_\eta^2+\eta^2\trace \mathbf{1}_{V_3}=(2n_1+n_3)\eta^2-2\trace R_\eta^2,
$$
which yields
$$
T_{33}=9n_3x_0^2+54(n_1-1)\xi^2+18\sqrt{6}\trace Q_\xi R_\eta +9(4n_1+n_2-2)\eta^2.
$$

Summing up the found relations, we obtain
\begin{equation*}
\begin{split}
\sum_{i,j=0}^3 T_{ij}=&18(3n_1+n_2+1)(x_n^2+\xi^2+\eta^2)+18\zeta^2+18\zeta^\TOP (\sum_{j=1}^{n_2}N_j^\TOP N_j) \zeta\\
&+54\zeta^\TOP (\sum_{i=1}^{n_1}Q_i^2)\zeta +36\sum_{j=1}^{n_2}\zeta^\TOP R_i^2\zeta+18\sqrt{6}\trace Q_\xi R_\eta.
\end{split}
\end{equation*}
Applying the $\xi$-Laplacian to (\ref{mm202}) and the $\zeta$-Laplacian to (\ref{mm004}), and taking into account that $\trace Q_i=\trace R_j=0$, we find respectively
$$
\zeta^\TOP(\sum_{j=1}^{n_2}N_j^\TOP N_j+\sum_{i=1}^{n_1}Q_i^2)\zeta=n_1\zeta^2,
\qquad
\zeta^\TOP(\sum_{j=1}^{n_2}R_j^2+\sum_{i=1}^{n_1}Q_i^2)\zeta=\frac{n_3+2}{2}\zeta^2,
$$
which yields
$$
18\zeta^\TOP (\sum_{j=1}^{n_2}N_j^\TOP N_j) \zeta+54\zeta^\TOP (\sum_{i=1}^{n_1}Q_i^2)\zeta
+36\sum_{j=1}^{n_2}\zeta^\TOP R_i^2\zeta=18(n_1+n_3+2)\zeta^2.
$$
Thus,
\begin{equation*}
\begin{split}
\sum_{i,j=0}^3 T_{ij}&=18(3n_1+n_2+1)(x_n^2+\xi^2+\eta^2+\zeta^2)+18\sqrt{6}\trace Q_\xi R_\eta.
\end{split}
\end{equation*}
In order to prove (\ref{eq:norm}) it remains to show only that $\trace Q_\xi R_\eta=0$, or equivalently that $\trace Q_iR_j=0$ for all $i,j$. To this end, let us fix an index $i$, $1\le i\le n_1$. Then (\ref{nummm1}) yields $Q_\xi^3=\xi^2Q_\xi$, hence $Q_i^3=Q_i$. This means that $Q_i$ has three eigenvalues: $\pm1 $ and $0$. Let $W_{\pm }$ and $W_0$ denote the corresponding eigenspaces and let $\zeta=(w_+, w_-,w_0)$ be the   vector decomposition corresponding to the decomposition $\zeta\in V_3=W_+\oplus W_-\oplus W_0$. Then (\ref{mm004}) yields
$$
\sum_{k=1,k\ne i}^{n_1}(\zeta^\TOP Q_k\zeta)^2+\sum_{j=1}^{n_2}(\zeta^\TOP R_j \zeta)^2=\zeta^4-(w_+^2-w_-^2)^2=4w_+^2w_-^2+2w_0^2(w_+^2+w_-^2)+w_0^4,
$$
implying that $R_j$ has the following block structure:
$$
R_j=\left(
        \begin{array}{ccc}
          \mathbf{0}_{W_+} & * &* \\
          * & \mathbf{0}_{W_-} & * \\
          * & * & M_j
        \end{array}
      \right),
$$
where $M_j$ is a symmetric endomorphism of $ W_0$ with $\trace M_j=\trace R_j=0$. We have
$$
Q_iR_j=\left(
        \begin{array}{ccc}
          \mathbf{0}_{W_+} & * &* \\
          * & \mathbf{0}_{W_-} & * \\
          * & * & -\frac{1}{2}M_j
        \end{array}
      \right),
$$
which yields that $\trace Q_iR_j=-\frac{1}{2}\trace M_j=0 $. This finishes the proof of the theorem.
\end{proof}

\subsection{Completion of the proof of Theorem~\ref{thC} and Theorem~\ref{thD}}\label{sec:finish}

Now we are able to finish the proof of Theorem~\ref{thD}. The following two corollaries establish the implications (c)$\Rightarrow$(a) and (b)$\Rightarrow$(a) in Theorem~\ref{thD}, respectively. Furthermore, Corollary~\ref{cor:reversal} also establishes the `only-if' part in Theorem~\ref{thC}.

\begin{corollary}
\label{cor:exc_cliff}
If a radial eigencubic $f$ possess the quadratic trace identity and $n_2\in \{0,5,8,14,26\}$ then $f$ is exceptional.

\end{corollary}

\begin{proof}
We argue by contradiction and assume that there is an eigencubic $f$ of Clifford type satisfying the quadratic trace identity. Then $\sigma_2(f)$ has a single eigenvalue, and it follows  from Proposition~\ref{pro:tau} that for the associated to $f$ Clifford system $\mathcal{A}$ the equality $m=q+1$ holds. But in that case (\ref{Cliffordtype}) yields $n_2=m+1-q=2$, a contradiction.

\end{proof}

\begin{corollary}\label{cor:reversal}
Let $f$ be a radial eigencubic in the normal form $(\ref{given})$. If $\psi_{030}$ is either identically zero or irreducible then $f$ is an exceptional eigencubic. Equivalently, if $f$ is a radial eigencubic of Clifford type then $\psi_{030}\not\equiv 0$ and reducible.
\end{corollary}

\begin{proof}
If $\psi_{030}\equiv 0$ then in view of (\ref{zz040}) we have $n_2=0$, hence by Proposition~\ref{pro:n20} $f$ is exceptional. If $\psi_{030}\not\equiv 0$ and irreducible then by the Eiconal Cubic Theorem, $n_2=3\ell+2$, $\ell=1,2,4,8$, and $\Delta\psi_{030}=0$, in particular Proposition~\ref{pro:norm} yields that $f$ possess the quadratic trace identity. If $n\ne 3m$, $m\in \{1,2,4,8\}$ then by Corollary~\ref{cor:exc_cliff}, $f$ is exceptional. Now suppose that $n=3k$, where $k\in \{1,2,4,8\}$, and assume that $f$ is an eigencubic of Clifford type. Then we have from (\ref{a111}) that $n_1=k-2\ell-1$ and $(\ref{eq:norm})$ yields
$$
\trace \mathrm{Hess}^2(f)=-\frac{1}{3}(3n_1+n_2+1)\,\lambda x^2=(\ell-k)\,\lambda x^2,
$$
hence
$$
\sigma_2(f)=-\frac{1}{\lambda}\trace \mathrm{Hess}^2(f)=(k-\ell)x^2.
$$
By our assumption, $f$ an eigencubic is of Clifford type, say, given by (\ref{Clifffff}). Then,  $n=2m+q+1$, where by Proposition~\ref{pro:tau}, $q=n_1=k-2\ell-1$, hence $n\equiv 3k=2(m-\ell)+k$. We find $k=m-\ell$. On the other hand, arguing as in Corollary~\ref{cor:exc_cliff}, we see that $m=q+1=k-2\ell$, hence $k=m+2\ell$. Since $\ell\ne 0$ we get a contradiction which proves that $f$ is an exceptional radial eigencubic.
\end{proof}

\section{Examples of exceptional eigencubics}\label{sec:someexamples}

In section~\ref{sec:Cartan} we already discussed the Cartan isoparametric polynomials which are also exceptional radial  eigencubics of types $(\ell+1,0)$, $\ell=1,2,4,8$. Below we exhibit more examples of exceptional eigencubics. These examples cover all realizable types of exceptional eigencubics presented in Table~\ref{tabsBas} above.

\subsection{Hsiang's trick}

\label{sec:HSIANG}
For our further purposes, we describe the Hsiang construction  with minor modifications. Let  $\mathfrak{G}'(k,\R{})$ be the vector space of quadratic forms of $k\ge 2$ real variables with trace zero. We identify $\mathfrak{G}'(k,\R{})$ with the vector space of  $k\times k$ real symmetric matrices of trace zero equipped with the  scalar product $\scal{X}{Y}= \frac{1}{2}\trace XY$, $X,Y\in \mathfrak{G}'(k,\R{})$. Then the orthogonal group $\SO(k)$ acts naturally on $\mathfrak{G}'(k,\R{})$ as  substitutions, and $\mathfrak{G}'(k,\R{})$ is invariant with respect to the action. Consider an isometry $i:\R{N}\to \mathfrak{G}'(k,\R{})$, where $N=\frac{k^2+k-2}{2}$, and let $\Gamma=i^*(\SO(k))$ be the corresponding pullback of  $\SO(k)$ into $\OO(N)$.
It is well-known that the coefficients of the characteristic polynomial
$$
\det (X-\lambda\mathbf{1})=\lambda^k+b_2(X)\lambda-b_3(X)+\ldots+(-1)^kb_k(X), \quad X\in \mathfrak{G}'(k,\R{}),
$$
form a complete set of basic invariants with respect to the $\SO(k)$ action \cite{Wallach}.
Then the polynomial forms $\beta_k(x)=i^*b_k\in \mathbb{R}[x_1,\ldots,x_N]$ are invariant under the action of the group $\Gamma$, and form a complete set of invariants:
\begin{equation}\label{invar}
\R{}[x_1,\ldots,x_N]^\Gamma=\R{}[\beta_2,\ldots,\beta_k]
\end{equation}
(as usually, $\R{}[x_1,\ldots,x_N]^\Gamma$ denotes the subring of $\Gamma$-invariant polynomials).
Hsiang proves that $L$ is an invariant operator in the sense that for any element $g\in \OO(N)$, the operator $L$ commutes with the linear substitution $g$: $g^* L=L g^*$. Since $\Gamma\subset \OO(N)$ and $L$ is invariant, we have $L(\beta_k)\in  \R{}[x_1,\ldots,x_N]^\Gamma$, hence  by (\ref{invar})
\begin{equation}\label{Bk}
L(\beta_k)\in \R{}[\beta_2,\ldots,\beta_k].
\end{equation}

We have $b_2(X)=\frac{1}{2}((\trace X)^2-\trace X^2)=-\frac{1}{2}\trace X^2$, therefore
$$
\beta_2=i^*(b_2)=-x^2.
$$
On the other hand, $\deg L(\beta_3)=5$, thus by (\ref{Bk}) there are $c_1,c_2\in \R{}$ such that
\begin{equation}\label{Hsiang1}
L(\beta_3)=c_1\beta_2\beta_3+c_2\beta_5.
\end{equation}

Now suppose $3\le k\le 4$. Then $\beta_5\equiv 0$, hence (\ref{Hsiang1}) reads as follows:
\begin{equation}\label{B3}
L(\beta_3)=c_1x^2\beta_3,
\end{equation}
which is equivalent to that $\beta_3$ is a \textit{radial} eigenfunction in $\R{N}$.  This yields two (irreducible) eigencubics: in $\R{5}$ and in $\R{9}$, for $k=3$ and $k=4$ respectively. The first example is easily identified with the Cartan isoparametric cubic $\theta_{1}$ in $\R{5}$, see (\ref{CartanFormula0}) above. Indeed, applying the same argument to the Laplacian and the length of the gradient which obviously are invariant operators, we find in view of $\deg\Delta \beta_3=1$  and $\deg |\nabla \beta_3|^2=4$ that
\begin{equation}\label{Hsiang2}
\Delta \beta_3=0,\qquad |\nabla \beta_3|^2=c_3\beta_2^2+c_4\beta_4,
\end{equation}
and in view of $k=3$ we have $\beta_4\equiv 0$, so that (\ref{Hsiang2}) becomes equivalent to the M\"{u}nzner-Cartan differential equations (\ref{Mun2}) for $\theta_1$. This can also be done  explicitly if one consider a map $i:\R{5}\to \mathfrak{G}'(3,\R{})$ given by
\begin{equation}\label{MAT}
X=i(x)=\frac{1}{\sqrt{2}}\left(\begin{array}{ccc}
x_{4}-\frac{1}{\sqrt{3}}x_{5} & x_{3} & x_2 \\
x_3 & -x_4-\frac{1}{\sqrt{3}}x_{5} &x_1\\
x_2 & x_1 &\frac{2}{\sqrt{3}}x_5 \\
\end{array}\right),\qquad x\in \R{5}.
\end{equation}
Then $\scal{X}{X}= \frac{1}{2}\trace X^2=x^2$, hence $i$ is indeed an isometry. Hence  $\beta_3\equiv \det X$ provides  an explicit determinantal representation for $\theta_1$ (up to a constant factor). An important feature of the obtained determinantal representation (\ref{MAT}) is that (in view of $\trace X=0$)
\begin{equation}\label{traccc}
\beta_3(x)=-\frac{1}{3}\trace X^3.
\end{equation}

In case $k=4$, the quartic form  $\beta_4$ is no more trivial, hence (\ref{B3}) yields a {non-homogeneous} radial eigencubic  $\beta_3$  in $\R{9}$ discovered by Hsiang in \cite{Hsiang67}. In this case,
$$
\beta_3=-\frac{1}{6}(\trace X)^3+\frac{1}{2}\trace X\trace X^2 -\frac{1}{3}\trace X^3=-\frac{1}{3}\trace X^3,
$$
hence $\beta_3$ possess also a trace identity like (\ref{traccc}). Up to a congruence, $\beta_3$ coincides with the exceptional radial eigencubic of type $(0,5)$ constructed in Example~\ref{ex:2} below.

The above constructions can be repeated literally for the corresponding Hermitian analogues  $\mathfrak{G}'(3,\Com{})\cong \R{8}$ and $\mathfrak{G}'(4,\Com{})\cong \R{15}$ with the special unitary group $SU(3)$ acting on them.  This yields respectively the Cartan isoparametric cubic $\theta_2$ of type $(3,0)$ in $\R{8}$ and the Hsiang (non-homogeneous) exceptional eigencubic of type $(0,8)$ in $\R{15}$. For $k=3$ it is still possible to obtain the quaternionic (non-commutative) and octonionic (non-associative) counterparts of the above constructions by using the maximal orbits of ${Sp}(3)$ and ${F}_4$ instead of $\SO(3)$ on the corresponding Hermitian matrices of trace zero. This yields the Cartan isoparametric eigencubics $\theta_4$  and $\theta_8$, respectively. For $k=4$, however, there is only the quaternionic counterpart in $\R{27}$ of type $(0,14)$ (some care is  needed to appropriately interpret the trace of the corresponding matrix).

\subsection{The case $n_1=0$}The following proposition shows that all types in Table~\ref{tabs} with $n_1=0$ except for $(0,26)$ are realizable.

\begin{proposition}\label{pro:n1=0}
For $n_1=0$, there exist only exceptional eigencubic of type $(0,3\ell+2)$, where $\ell=1,2,4$.
\end{proposition}

\begin{proof}

Let us suppose that $f\in E(0,n_2)$ is an exceptional  eigencubic given in the normal form (\ref{given}). Then $n_2\ne0$, hence the condition $n_2=3\ell +2$, $\ell=1,2,4,8$, yields $n_3=n_2-2=3\ell $,  and the condition $\dim V_1=n_1=0$  yields $\psi=\psi_{012}+\psi_{030}$. Let us consider
$$
 \chi(\eta,\zeta)=\psi_{030}-\psi_{012}.
$$
Then from (\ref{Psi1})  and $\Delta_\eta\psi_{030}=0$ we see that $\chi$ is harmonic. Furthermore, by virtue of (\ref{zz004}), (\ref{zz040}) and (\ref{zz022})
\begin{equation}
\begin{split}
|\nabla \chi|^2&=(\psi_{030})^2_\eta-2(\psi_{030})^\TOP_\eta (\psi_{012})_\eta+(\psi_{012})^2_\eta +(\psi_{012})^2_\zeta=\frac{9}{2}(\eta^2+\zeta^2)^2,
\end{split}
\end{equation}
hence $\sqrt{2}\chi$ satisfies that Cartan-M\"unzner equations (\ref{Mun2}). Since $\ell\ge1$, the Cartan theorem implies that the dimension $\dim V_2+\dim V_3=2n_2-2=6\ell+2$ can take only values $5,8,14,26$. This implies $\ell\in \{1,2,4\}$.

In the converse direction, let us assume that $\ell\in \{1,2,4\}$ and show how to construct an exceptional eigencubic of type $(0,3\ell+2)$. To this end, we consider the division algebras $\mathbb{F}_{2\ell}$ and $\mathbb{F}_{\ell}$. Regarding $\mathbb{F}_{2\ell}$ as the Cayley-Dickson doubling of $\mathbb{F}_{\ell}$ \cite{Baez}, we write
\begin{equation}\label{FeF}
\mathbb{F}_{2\ell}=\mathbb{F}_{\ell}\oplus \mathbb{F}_{\ell},
\end{equation}
where the multiplication and conjugation on $\mathbb{F}_{2\ell}$ is given by
\begin{equation}\label{multipli}
(a\oplus b)(c\oplus d) =(ac-d\bar b)\oplus (\bar a d+cb), \qquad \overline{(a\oplus b)}=(\bar a, -b).
\end{equation}
Let $\gamma_{i}: \mathbb{F}_{2\ell}\to \mathbb{F}_\ell$ denote the canonical projection on the $i$th component. By identifying $\mathbb{F}_\ell$ and $\R{\ell}$, we consider the induced projections
\begin{equation}\label{gamma1}
\widetilde \gamma_i:\,\R{6\ell+2}\cong\R{2}\oplus \mathbb{F}_{2\ell}\oplus \mathbb{F}_{2\ell}\oplus \mathbb{F}_{2\ell}\,\,\rightarrow\,\,
\R{3\ell+2}\cong\R{2}\oplus \mathbb{F}_{\ell}\oplus \mathbb{F}_{\ell}\oplus \mathbb{F}_{\ell}.
\end{equation}
Then it follows readily from the definition of the Cartan polynomials (\ref{CartanFormula0}) that $\theta_\ell=\theta_{2\ell}\circ \widetilde{\gamma}_1$. Let $x=(x_1,x_2,z_1,z_2,z_3)\in \R{2}\oplus \mathbb{F}_{2\ell}\oplus \mathbb{F}_{2\ell}\oplus \mathbb{F}_{2\ell}=\R{6\ell+2}$ and let also $z_i=z_i'\oplus z_i''$ according to (\ref{FeF}). Then we denote
\begin{equation}\label{etazeta}
\eta:=(x_1,x_2,z_1',z_2',z_3')\in \R{3\ell+2},\qquad \zeta:=(z_1'',z_2'',z_3'')\in \R{3\ell}.
\end{equation}
In this notation,
\begin{equation*}
\begin{split}
\theta_{2\ell}(\eta,\zeta)& =x_{1}^3+\frac{3x_1}{2}(2|z_3|^2-|z_1|^2-|z_2|^2+2x_{2}^2)
+\frac{3\sqrt{3}}{2}[x_{2}(|z_1|^2-|z_2|^2)+\re z_1z_2z_3],\\
\end{split}
\end{equation*}
and
\begin{equation}\label{etha}
\theta_{\ell}(\eta)=\theta_{2\ell}(\eta,0),
\end{equation}
where the real part is defined by (\ref{realpart}).
We claim that
\begin{equation}\label{Feta}
F(\eta,\zeta,x_{3\ell+3}):=x_{3\ell+3}^3+\frac{3}{2}(\zeta^2-\eta^2)x_{3\ell+3}+\frac{1}{\sqrt{2}}(\theta_{2\ell}(\eta,\zeta)-2\theta_{2\ell}(\eta,0)).
\end{equation}
is the required exceptional eigencubic of type $(0,3\ell+2)$ in $\R{3\ell+3}$. First we observe that the above expression is written in the norm form, where $\phi=\frac{3}{2}(\zeta^2-\eta^2)$ and $\psi=\frac{\theta_{2\ell}(\eta,\zeta)-2\theta_{2\ell}(\eta,0)}{\sqrt{2}}$. Setting $\eta\in V_2=\R{3\ell+2}$ and $\zeta\in V_3=\R{3\ell}$, we determine the decomposition of $\psi$ into the corresponding homogeneous parts $\eta^i\otimes\zeta^{3-i}$. To this end, we denote $\psi=\psi_{030}+\psi_{012},$ where
\begin{equation}\label{eqpsi}
\psi_{030}=-\frac{1}{\sqrt{2}}\theta_{2\ell}(\eta,0),\,\,\,\,\, \psi_{012}:=\frac{1}{\sqrt{2}}(\theta_{2\ell}(\eta,\zeta)-\theta_{2\ell}(\eta,0)),
\end{equation}
and notice that $\psi_{030}\in \eta\otimes\eta\otimes \eta$ by the definition and also $\psi_{012}\in \eta\otimes \zeta\otimes \zeta.$ Indeed, to verify the latter identity, we write
$$
\psi_{012}=\frac{3x_1}{2}(2|z_3''|^2-|z_1''|^2-|z_2''|^2)
+\frac{3\sqrt{3}}{2}[x_{2}(|z_1''|^2-|z_2''|^2)+D]
$$
where $D=\re z_1z_2z_3-\re z_1'z_2'z_3'$.
Observe that the multiplication of $z_i'\in \mathbb{F}_{\ell}$ is  associative for $\ell\le 4$, and that in view of (\ref{multipli})
\begin{equation}\label{imagin}
\re (0\oplus b)=0, \qquad b\in \mathbb{F}_\ell.
\end{equation}
Then (\ref{multipli}) yields
\begin{equation*}
\begin{split}
\re (z_1'\oplus z_1'')\cdot(z_2'\oplus z_2'')\cdot(z_3'\oplus z_3'')\re (z_1'z_2'z_3'-z_2''\bar z_1''z_3'-z_3''\bar z_1''\bar z_2'-z_3''\bar z_2''z_1'),
\end{split}
\end{equation*}
thus
$$
D=-\re (z_2''\bar z_1''z_3'+z_3''\bar z_1''\bar z_2'+z_3''\bar z_2''z_1')\in \eta\otimes \zeta^{2},
$$
which, in vie  of (\ref{etazeta}), yields $\psi_{012}\in \eta\otimes \zeta^{2}$.

In order to show that $F$ is an eigencubic, we need by Proposition~\ref{pro:main} to verify that $\phi$ and $\psi$ satisfy equations (\ref{e5new}) and (\ref{e6new}). The former equation is equivalent to the system  in Proposition~\ref{pr:hidden}, where equations (\ref{zz202})--(\ref{zz211}), (\ref{zz121}) are trivial because $n_1=0$, equation (\ref{zz040}) is satisfied by our choice of $\psi_{030}$ in (\ref{eqpsi}):
\begin{equation}\label{ZZ040}
(\psi_{030})_\eta^2=\frac{9}{2}|\nabla\theta_{\ell}(\eta)|=\frac{9}{2}\eta^4,
\end{equation}
and equations (\ref{zz004}) and (\ref{zz022}) have  the following form:
\begin{equation}\label{ZZ004}
6(\psi_{012})_\eta^2=27\,\zeta^4,
\end{equation}
and
\begin{equation}\label{ZZ022}
2(\psi_{030})_\eta^\TOP (\psi_{012})_\eta-(\psi_{012})_\zeta^2=-9\,\zeta^2\eta^2,
\end{equation}
respectively. But these equations follows immediately by collecting the terms by homogeneity in the identity
$$
\frac{9}{2}(\eta^2+\zeta^2)^2=(\psi_{012})^2_\eta-2(\psi_{012})^\TOP_\eta (\psi_{030})_\eta+(\psi_{030})^2_\eta +(\psi_{012})^2_\zeta,
$$
which is obtained from $|\nabla \theta_{2\ell}(\eta,\zeta)|^2=9(\eta^2+\zeta^2)^2$, where $\theta_{2\ell}(\eta,\zeta)=\sqrt{2}(\psi_{012}-\psi_{030})$  by virtue of (\ref{eqpsi}).

Thus, it remains only to verify (\ref{e6new}). To this end, we note that by (\ref{ZZ022}), (\ref{ZZ004}) and (\ref{ZZ040})
\begin{equation*}
\begin{split}
|\nabla \psi|^2&=(\psi_{012})_{\zeta}^2+(\psi_{012})_{\eta}^2+2(\psi_{012})^\TOP_\eta (\psi_{030})_\eta+(\psi_{030})_{\eta}^2 \\ 
&=\frac{9}{2}(\eta^2+\zeta^2)^2+4(\psi_{012})^\TOP_\eta (\psi_{030})_\eta,
\end{split}
\end{equation*}
hence
\begin{equation*}
\begin{split}
\nabla\psi\cdot \nabla(|\nabla \psi|^2)&=54\psi(\eta^2+\zeta^2)+4\psi_\zeta^\TOP(\psi_{012})_{\zeta\eta} (\psi_{030})_\eta+4\psi_\eta^\TOP(\psi_{030})_{\eta\eta} (\psi_{012})_\eta\\
&=54\psi(\eta^2+\zeta^2)+4(\psi_{012})_\zeta^\TOP(\psi_{012})_{\zeta\eta} (\psi_{030})_\eta+4(\psi_{012}+\psi_{030})_\eta^\TOP(\psi_{030})_{\eta\eta} (\psi_{012})_\eta
\end{split}
\end{equation*}
We have by (\ref{ZZ022}) and (\ref{ZZ040})
\begin{equation*}
\begin{split}
4(\psi_{012})_\zeta^\TOP(\psi_{012})_{\zeta\eta} \,(\psi_{030})_\eta&2((\psi_{012})_\zeta^2)_\eta^\TOP(\psi_{030})_\eta\\
&=2\bigl(2(\psi_{030})_\eta^\TOP (\psi_{012})_\eta+9\zeta^2\eta^2\bigr)_\eta^\TOP\,(\psi_{030})_\eta\\
&=4(\psi_{012})_\eta^\TOP (\psi_{030})_{\eta\eta}(\psi_{030})_\eta+108\zeta^2\psi_{030}\\
&=2(\psi_{012})_\eta^\TOP ((\psi_{030})_{\eta}^2)_\eta+108\zeta^2\psi_{030}\\
&=36\eta^2 \psi_{012}+108\zeta^2\psi_{030}\\
\end{split}
\end{equation*}
and similarly
\begin{equation*}
\begin{split}
4(\psi_{030})_\eta^\TOP(\psi_{030})_{\eta\eta} (\psi_{012})_\eta=36\eta^2 \psi_{012}.
\end{split}
\end{equation*}
We also have by (\ref{ZZ022}) and (\ref{ZZ004})
\begin{equation*}
\begin{split}
4(\psi_{012})_\eta^\TOP(\psi_{030})_{\eta\eta} (\psi_{012})_\eta&=4\bigl((\psi_{012})_\eta^\TOP(\psi_{030})_\eta\bigr)_\eta^\TOP\,(\psi_{012})_\eta\\
&=2\bigl(-9\zeta^2\eta^2+(\psi_{012})_\zeta^2\bigr)_\eta^\TOP\,(\psi_{012})_\eta\\
&=-36\zeta^2\psi_{012}+4(\psi_{012})_\zeta^\TOP (\psi_{012})_{\zeta\eta}(\psi_{012})_\eta\\
&=-36\zeta^2\psi_{012}+2(\psi_{012})_\zeta^\TOP \bigl((\psi_{012})_{\eta}^2\bigr)_\zeta\\
&=36\zeta^2\psi_{012}\\
\end{split}
\end{equation*}

Combining the found identities we obtain
\begin{equation*}
\begin{split}
\psi_{\bar x}^{\TOP} \psi_{{\bar x}{\bar x}}\psi_{\bar x}\equiv \frac{1}{2}\nabla\psi\cdot \nabla(|\nabla \psi|^2)&=27\psi+36\eta^2 \psi_{012}+54\zeta^2\psi_{030}+18\zeta^2\psi_{012}\\
&=27\psi_{030}(\eta^2+3\zeta^2)+9\psi_{030}(7\eta^2+5\zeta^2)\\
\end{split}
\end{equation*}
Furthermore,
$$
2\phi\, \phi_{\bar x}^{\TOP} \psi_{{\bar x}}=9(\zeta^2-\eta^2)(\zeta^\TOP\psi_\zeta-\eta^\TOP\psi_\eta)9(\zeta^2-\eta^2)(\psi_{012}-3\psi_{030}),
$$
which finally yields
$$
2\phi\, \phi_{\bar x}^{\TOP} \psi_{{\bar x}}+\psi_{\bar x}^{\TOP} \psi_{{\bar x}{\bar x}}\psi_{\bar x}=54(\eta^2+\zeta^2)\psi
$$
and proves (\ref{e6new}). Thus, $F$ is indeed a radial eigencubic.  From (\ref{Feta}) we see that $(n_1,n_2)=(0,3\ell+2)$. On the other hand, by our construction, $\Delta\psi_{030}=\frac{1}{\sqrt{2}}\Delta\theta_\ell(\eta)=0$, hence Corollary~\ref{cor:reversal} yields that $F$ is exceptional.
The proposition is proved completely.
\end{proof}

\begin{example}\label{ex:2}
Let us consider the following cubic polynomial in $\R{9}$ given by
\begin{equation*}\label{eigf2}
f(x):=\re \prod_{s=1}^3(x_{3s-2}\mathrm{i}+x_{3s-1}\mathrm{j}+x_{3s}\mathrm{k})\equiv \det \left(
        \begin{array}{ccc}
          x_1&x_2&x_3 \\
          x_4&x_5&x_6 \\
          x_7&x_8&x_9 \\
        \end{array}
      \right),
\end{equation*}
where $\mathrm{i},\mathrm{j},\mathrm{k}$ is the standard basis elements in quaternion algebra $\mathbb{H}=\mathbb{F}_4$. Then it is straightforward to see that $f$  satisfies (\ref{mainlambda}) with $\lambda=-2$. Observe also that $f$ satisfies the quadratic trace identity
\begin{equation}\label{label}
\sigma_2(f)=2\sum_{i=1}^9 x_i^2,
\end{equation}
hence, by Corollary~\ref{cor:exc_cliff}, $f$ is an exceptional eigencubic. From (\ref{eq:norm}) and (\ref{label}), $3n_1+n_2+1=6$ and in virtue of (\ref{a111}) $9=n=3n_1+2n_2-1$, hence $n_1=0$ and $n_2=5$. Thus $f$ is an exceptional eigencubic of type $(0,5)$. It was already mentioned in section~\ref{sec:HSIANG} that $f$ is congruent to the Hsiang example in $\R{9}$ mentioned in \cite{Hsiang67}.
\end{example}

\subsection{The case $n_1=1$.}

\begin{proposition}\label{pro:n2=1}
If $n_1=1$ then there are exactly four (congruence classes of) exceptional eigencubics of types $E(1,3\ell +2)$, $\ell=1,2,4,8$. Any such a cubic is congruent to
$$
f(t)=\bigl[t_{3\ell+3}^3-\frac{3}{2}t_{3\ell+3}(t_1^2+\ldots t_{3\ell+2}^2)+\frac{1}{\sqrt{2}}\theta_\ell(t_1,\ldots, t_{3\ell+2})\bigr]_{\mathbb{C}},\quad t\in \R{3\ell +3},
$$
where
$$
[g]_{\mathbb{C}}(x,y):=\frac{1}{2}(g(x+\I y)+g(x-\I y)), \qquad x,y\in \R{3\ell+3}
$$ is the complex doubling of a polynomial $g$.

\end{proposition}

\begin{proof}
We assume that $f$ is given in the normal form (\ref{given}). By the assumption  $\dim V_1=1$, hence  $n_3=n_2$. Since $f$ is exceptional, we also have $n_2=3\ell +2$, where $\ell=1,2,4,8$, see Table~\ref{tabs}. Furthermore, by identifying the coefficient of $\xi_1^2$ in the first identity in  (\ref{nummm1}) one finds
$$
P_1^\TOP P_1+Q_1^2 = \mathbf{1}_{V_3}, \qquad
P_1 P_1^\TOP = \mathbf{1}_{V_2},
$$
which yields $\trace Q_1^2=n_3-n_2=0$, hence $Q_1=0$. Thus $\psi_{102}\equiv 0$.

In the notation of section~\ref{sec:def}, $\psi_{111}=3\xi_1 N_\eta \zeta$, where $N_i$, $1\le i\le 3\ell +2$, are matrices of size $1\times (3\ell+2)$. By (\ref{mm202}) and (\ref{mm220}) we see that $\{N_1^\TOP\,\ldots, N_{3\ell +2}^\TOP\}$ forms an orthonormal basis in $V_3=\R{3\ell +2}$. Using the freedom to choose the orthonormal basis elements, we can ensure that $N_i^\TOP=e_i$, where $\{e_i\}$ is the standard orthonormal basis in $\R{3\ell +2}$. This  yields
$$
N_\eta = N_1^\TOP\eta_1+\ldots+ N_{n_2}^\TOP\eta_{n_2}=\eta^\TOP,\qquad \eta\in V_2,
$$
hence $\psi_{111}=3\xi_1 \eta^\TOP\zeta$ and we have for the normal form
\begin{equation}\label{norrr}
f=x_n^3-3\xi_1x_n^2+\frac{3}{2}(\zeta^2-\eta^2)x_n+3\xi_1\eta^\TOP \zeta+\psi_{030}+\psi_{012}
\end{equation}
Moreover, in this notation (\ref{mm121}) reads as
\begin{equation}\label{terra}
H^\TOP=-\eta^\TOP R_\eta,
\end{equation}
where $H(\eta)=\frac{\sqrt{2}}{3}\nabla_\eta \psi_{030}(\eta)$ and $\psi_{012}=\frac{3\sqrt{2}}{2}\zeta^\TOP R_\eta \zeta$. This yields
$$
\frac{\sqrt{2}}{3}(\psi_{030})_{\eta_i\eta_j}= (H_i)_{\eta_j}=(H_i)_{\eta_j},
$$
and
$$
\eta^\TOP R_j e_i +e_j^\TOP R_\eta e_i=\eta^\TOP R_i e_j +e_i^\TOP R_\eta e_j.
$$
By virtue of symmetricity of matrix $R_\eta$,
\begin{equation}\label{Poinc}
\frac{1}{2}(\partial_{\eta_j}r_i-\partial_{\eta_i}r_j)=\eta^\TOP(R_ie_j-R_je_i)=0,
\end{equation}
where $r_i=\eta^\TOP R_i\eta$. The latter relation yields that there exists a homogeneous cubic  polynomial $r=r(\eta)$ such that $r_i=\partial_{\eta_i}r$. By the homogeneity of $r$ and (\ref{terra}) we have
$$
r(\eta)=\frac{1}{3}\sum_{i=1}^{n_2}\eta_i \partial_{\eta_i}r\equiv \frac{1}{3}\sum_{i=1}^{n_2}\eta_i r_i=\frac{1}{3}\eta^\TOP R_\eta\eta=-\frac{1}{3}H^\TOP\eta=-\frac{\sqrt{2}}{9}(\psi_{030})_\eta^\TOP \eta=-\frac{\sqrt{2}}{3}\psi_{030}.
$$
This yields, in particular,
\begin{equation}\label{RetA}
\psi_{030}=-\frac{1}{\sqrt{2}}\eta^\TOP R_\eta\eta.
\end{equation}
By choosing an orthonormal basis in $V_2=\R{3\ell+2}$ we may ensure (in view of the Eiconal Cubic Theorem and (\ref{zz004})) that  $\psi_{030}(\eta)=\frac{1}{\sqrt{2}}\theta_\ell(\eta)$, where $\theta_\ell(\eta)$ is the Cartan polynomial (\ref{CartanFormula0}). Then (\ref{RetA}) becomes $\eta^\TOP R_\eta\eta=-\theta_\ell(\eta).$ Setting $\eta$ by $\eta\pm \I \zeta$ in the latter identity, where $\I^2=-1$, and summing we get
$$
\frac{1}{2}(\theta_\ell(\eta+\I \zeta)+\theta_\ell(\eta-\I \zeta))=-\eta^\TOP R_\eta\eta+2\zeta^\TOP R_\zeta\eta+\zeta^\TOP R_\eta\zeta.
$$
From (\ref{Poinc}), $\zeta^\TOP R_ie_j=\zeta^\TOP R_je_i$, hence $\zeta^\TOP R_\eta \zeta=\zeta^\TOP R_\zeta \eta$, and the above relation yields
\begin{equation}\label{inview}
\frac{1}{2}(\theta_\ell(\eta+\I \zeta)+\theta_\ell(\eta-\I \zeta))=\sqrt{2}\psi_{030}+3\zeta^\TOP R_\eta \zeta\equiv \sqrt{2}(\psi_{030}+\psi_{012}).
\end{equation}

Now let us define a new cubic polynomial
$$
g(t):=t_{3\ell+3}^3-\frac{3}{2}t_{3\ell+3}(t_1^2+\ldots t_{3\ell+2}^2)+\frac{1}{\sqrt{2}}\theta_\ell(t_1,\ldots, t_{3\ell+2}).
$$
Then, setting $t_{\pm}=(\eta_1\pm\zeta_1\I,\ldots, \eta_{3\ell +2}\pm\zeta_{3\ell+2}\I,x_n\pm\xi\I)$, one easily finds that by virtue of (\ref{inview}) and (\ref{norrr}) that
\begin{equation*}
\begin{split}
[g]_{\mathbb{C}}&\equiv\frac{1}{2}(g(t_-)+g(t_+))\\
&=x_n^3-3\xi_1(x_n^2-\eta^\TOP \zeta)+\frac{3}{2}(\zeta^2-\eta^2)x_n+\frac{1}{2\sqrt{2}}(\theta_\ell(\eta+\I \zeta)+\theta_\ell(\eta-\I \zeta))\\
&=f,
\end{split}
\end{equation*}
which yields the required representation and  proves the proposition.
\end{proof}

\subsection{An exceptional eigencubic of type $(4,5)$ in $\R{21}$}\label{sec:octonions}

Now we exhibit an example of an exceptional eigencubic of type $E(4,5)$ in $\R{21}$.  Let $o_0=1,o_1,\ldots, o_7$ be a basis for the octonion algebra, $\mathbb{F}_8=\mathbb{O}$. The multiplication table is given by $o_i\cdot o_{i+1}=o_{i+3}$, where the indices are permuted cyclically and translated modulo 7, see for instance \cite{Baez}.  For any vector $u=(u_1,\ldots,u_7)\in \R{7}$ we denote by $o_u$ the imaginary octonion $o_u=u_1o_1+\ldots+u_7o_7$. Let us consider the cubic form
\begin{equation}\label{fO}
f\equiv f_{\mathbb{F}_8}=\re  o_u(o_v o_w), \qquad x=(u,v,w)\in \R{21},
\end{equation}
where $u=(x_1,\ldots,x_7)$, $v=(x_{8},\ldots,x_{14})$, $w=(x_{15},\ldots,x_{21})$. We shall need the following known properties of the real part (see, for instance, Corollary~15.12 in \cite{Adams}): for any three octonions $\alpha,\beta,\gamma\in \mathbb{F}_8$
\begin{equation}\label{commuta}
\re \alpha\beta=\re \beta\alpha,
\end{equation}
\begin{equation}\label{assoc}
\re (\alpha\beta)\gamma=\re \alpha(\beta\gamma),
\end{equation}
\begin{equation}\label{cyclic}
\re (\alpha\beta)\gamma=\re (\beta\gamma)\alpha.
\end{equation}
Since $\re \alpha=\re \bar \alpha$, where $\bar \alpha$ denotes the conjugate octonion, we have for \textit{imaginary} octonions
$$
\re o_u(o_v o_w)=\re \overline{o_u(o_v o_w)})=\re \overline{(o_v o_w})\bar o_u=-\re  (o_w o_v)o_u=-\re (o_w o_v)o_u.
$$
This yields by virtue of (\ref{cyclic}) and (\ref{assoc})
\begin{equation}\label{antis}
\re o_u(o_v o_w)=-\re  o_u(o_w o_v),
\end{equation}
i.e. the cubic form $f$  is an alternating function in the three variables.

\begin{proposition}
The cubic form $f=\re  o_u(o_v o_w)$ is an exceptional eigencubic in $\R{21}$ of class $E(4,5)$.
\end{proposition}

\begin{proof}
Since $f$ is a trilinear form in $u,v,w$, it is harmonic. We have
$$
-L(f)=2(\sum_{i,j=1}^7f_{u_i}f_{v_j}f_{u_iv_j}+\sum_{i,k=1}^7f_{u_i}f_{w_k}f_{u_iw_k}+\sum_{j,k=1}^7f_{v_j}f_{w_k}f_{v_jw_k}).
$$
By symmetry, it suffices to find the first sum. We have
$$
\sum_{i,j=1}^7f_{u_i}f_{v_j}f_{u_iv_j}=\sum_{i,j=1}^7\re  o_i(o_vo_w)\,\re  o_i (o_j o_w)\,\re  o_u (o_j o_w).
$$
Notice that for $\alpha,\beta\in \mathbb{F}_8$
\begin{equation}\label{conjuga}
\sum_{i,j=1}^7\re  o_i\alpha\,\re  o_i \beta=\sum_{i,j=1}^7\alpha_i\beta_i\equiv \re \alpha\bar \beta-\re \alpha\,\re \beta,
\end{equation}
Hence
\begin{equation*}
\begin{split}
\sum_{i,j=1}^7f_{u_i}f_{v_j}f_{u_iv_j}&\sum_{j=1}^7[\re  ((o_vo_w) \overline{(o_j o_w)})-\re  o_vo_w\,\re o_j o_w]\re  o_u (o_j o_w)\equiv \Sigma_1-\Sigma_2
\end{split}
\end{equation*}

We have
\begin{equation*}
\begin{split}
\Sigma_1=\sum_{j=1}^7\re  ((o_vo_w) (o_w o_j))\,\re  o_u (o_j o_w)\sum_{j=1}^7\re  o_v(o_w(o_w o_j))\,\re  o_u (o_j o_w)
\end{split}
\end{equation*}
Since the octonions are alternative, that is, products involving no more than 2 independent octonions do associate, we find
$$
\re o_v (o_w(o_w o_j))= \re o_v ((o_wo_w) o_j)=-\re o_v ((o_w\overline{o}_w) o_j)=-w^2\re o_v o_j,
$$
where $o_w\overline{o}_w=|w|^2\equiv w^2$ is the norm of the vector $w$. By (\ref{commuta}), $\re o_v o_j=\re o_v o_j$ and by (\ref{assoc}), $\re  o_u (o_j o_w)=\re  o_j(o_w o_u)$, so applying (\ref{conjuga}) and (\ref{antis}) we obtain
\begin{equation*}
\begin{split}
\Sigma_1&=-w^2\sum_{j=1}^7\re  o_jo_v\,\re  o_j (o_w o_u)=-w^2\re  o_v\overline{(o_w o_u})\\
&=-w^2\re  o_v(o_u o_w)=w^2\re  o_v(o_w o_u)=w^2\re  o_u(o_v o_w)\\
&= w^2 f.
\end{split}
\end{equation*}

Similarly we obtain
\begin{equation*}
\begin{split}
\Sigma_2&=\sum_{j=1}^7\re  o_vo_w \, \re o_w o_j\,\re  o_u (o_j o_w)\sum_{j=1}^7\re  o_vo_w \, \re o_j o_w\,\re  o_j (o_w o_u)\\
&=\re  o_vo_w \, \re o_w\overline{(o_w o_u})=0\\
\end{split}
\end{equation*}
because $\re o_w\overline{(o_w o_u})=\re o_w o_u o_w=-w^2\re o_u=0$.

Combining the found identities together, we get
$$
L(f)=-2(u^2+v^2+w^2)f\equiv -2 x^2 f,
$$
hence $f$ is a radial eigencubic with $\lambda(f)=-2$.
Let us verify that $f$ is indeed an exceptional eigencubic. We have
$$
\trace \mathrm{Hess}^2 f=2(\sum_{i,j=1}^7f^2_{u_iv_j}+\sum_{i,k=1}^7f_{u_iw_k}^2+\sum_{j,k=1}^7f_{v_jw_k}^2),
$$
where
\begin{equation*}
\begin{split}
\sum_{i,j=1}^7f_{u_iv_j}^2&=\sum_{i,j=1}^7\re  o_i (o_j o_w)\,\re  o_i (o_j o_w)=\sum_{j=1}^7\bigl(\re  (o_j o_w)(\overline{o_j o_w})-\re o_j o_w\, \re o_jo_w\bigr)\\
&=\sum_{j=1}^7(\re  (o_j o_w)(o_w o_j)-w_j^2)=6w^2.
\end{split}
\end{equation*}
Thus,
\begin{equation}\label{squarenorm}
\trace (\mathrm{Hess} f)^2=12x^2.
\end{equation}

By Corollary~\ref{cor:exc_cliff}, $f$ is an exceptional eigencubic. In view of Table~\ref{tabs} and Proposition~\ref{pr:n1=2}, we have the only possibility: $n_1=4$. In order to see this also directly, we observe that by virtue of $\lambda(f)=-2$ and (\ref{squarenorm}), $\sigma_2(f)=6x^2$, so that (\ref{eq:norm}) yields $3n_1+n_2+1=18$. Since $21=n=3n_1+2n_2-1$, we find $n_1=4$ and $n_2=5$, as required.

\end{proof}

\subsection{Some remarks}
The four Cartan isoparametric eigencubics $\theta_\ell$ are well-known and  appear in various mathematical contexts. We  mention only a recent interest in  $\theta_\ell$ in special Riemannian geometries satisfying the nearly integrabilty condition \cite{Nur}, \cite{GN}, explicit solutions to the Ginzburg-Landau system \cite{Farina}, \cite{GeXie}, the harmonic analysis of cubic isoparametric minimal hypersurfaces \cite{BS1}, \cite{BS2}.

We would  like to mention that the  octonionic trilinear form discussed  in section~\ref{sec:octonions}  and its quaternionic analogue occur also in calibrated geometries (as an associative calibration on $\mathrm{Im}\, \mathbb{O}$) \cite{HL} and in constructing of singular solutions of Hessian fully nonlinear elliptic equations  \cite{NV1}, \cite{NV2}.

\section{Radial eigencubics and isoparametric quartics}\label{sec:isop}

\subsection{The degenerate form} By Proposition~\ref{pro:n20}, any radial eigencubic having property $n_2=0$ is exactly on of the  four Cartan's isoparametric eigencubics $\theta_\ell$. In this section we study non-isoparametric eigencubics, i.e. those with $n_2\ne 0$. An advantage of working with the degenerate form is that it is more symmetric and can easily be converted into a purely matrix representation. Another important aspect of the degenerate form is that it establishes a correspondence between eigencubics and isoparametric eigencubics which will be studied in the next section.

\begin{proposition}\label{th:dege}
Any non-isoparametric radial eigencubic $f$ in $\R{n}$, normalized by $\lambda(f)=-8$, in some orthogonal coordinates has the form
\begin{equation}\label{degeeabc}
f=(u^2-v^2)x_n+a(u,w)+b(y,w)+c(u,y,w),
\end{equation}
where $u=(x_1,\ldots, x_m)$, $v=(x_{m+1},\ldots, x_{2m})$, $w=(x_{2m+1},\ldots, x_{n-1})$, and the cubic forms
$a\in u\otimes w^2$, $b\in v\otimes w^2$, $c\in u\otimes v\otimes w$ satisfy the system
\begin{eqnarray}
&&\Delta_w a=\Delta_w b=0,\label{harmabw}\\
&&a_u^2=b_v^2,\label{r17}\\
&&a_w^\TOP   c_w=b_w^\TOP   c_w=0,\label{r15}\\
&&c_v^2+2a_w^2=4u^2w^2,\label{r16} \\
&&c_u^2+2b_w^2=4v^2w^2,\label{r16r}\\
&&c_u^{\TOP}a_{uw}b_w=c_v^{\TOP}b_{vw}a_w=0, \label{s023}\\
&&a_u^{\TOP}c_{uv}b_v=0.  \label{s005}\\
&&2c_w^{\TOP}c_{wu}a_u+2a_w^{\TOP}b_{ww}b_w+b_w^{\TOP}a_{ww}b_w=12av^2,  \label{s122}\\
&&2c_w^{\TOP}c_{wv}b_v+2b_w^{\TOP}a_{ww}a_w+a_w^{\TOP}b_{ww}a_w=12bu^2,  \label{s122b}
\end{eqnarray}

Conversely, any cubic polynomial (\ref{degeeabc}) satisfying the above system is a non-isoparametric eigencubic.
\end{proposition}

\begin{proof}
First notice that in some orthogonal coordinates $f$ is linear with respect to some coordinate function. Indeed, let (\ref{given}) be the normal form of $f$. Then by our assumption $V_2\ne \emptyset$ and by (\ref{zz040}) $(\psi_{030})_\eta^2=\frac{9}{2}\,\eta^4$, hence
by writing $\psi_{030}$ in the normal form
$$
\psi_{030}=\frac{\sqrt{2}}{2}\eta_{1}^3+\eta_{1}\phi(\bar \eta)+\psi(\bar\eta), \quad \bar\eta=(\eta_2,\ldots,\eta_{n_2}),
$$
we get from (\ref{given}) that
\begin{equation}\label{s59}
\begin{split}
f&=x_n^3+\phi x_n+\psi_{030}+\psi_{111}+\psi_{102}+\psi_{012} (x_n+\eta_{1}\sqrt{2})(x_n-\frac{\eta_{1}}{\sqrt{2}})^2+F(x),
\end{split}
\end{equation}
where $F$ is a linear form in the variables $\eta_{1}$ and $x_{n}$. Then applying rotation $u=\frac{x_n+\eta_{1}\sqrt{2}}{\sqrt{3}}$, $v=\frac{x_n\sqrt{2}-\eta_{1}}{\sqrt{3}}$ in the $(\eta_1,x_n)$-plane
we conclude that $f$ becomes linear in $u$ in the new coordinates.

Thus, we may assume without loss of generality that $f$
\begin{equation}\label{fuu}
f=x_n\Phi(\barr )+\Psi(\barr ), \qquad \barr =(x_1,\ldots,x_{n-1}).
\end{equation}
Moreover, we shall assume that $f$ is normalized by $\lambda(f)=-8$. Then using (\ref{fuu}) we obtain from (\ref{mainlambda}) by identifying the coefficients of $x_n^k$
\begin{eqnarray}
&&\label{r1} \Phi_{\barr }^\TOP  \Phi_{\barr \barr }\Phi_{\barr }=8\Phi,\\
&&2\Phi_{\barr }^\TOP  \Phi_{\barr \barr }\Psi_{\barr } + \Phi_{\barr }^\TOP  \Psi_{\barr \barr }\Phi_{\barr }=8\Psi,
\label{r2}\\
&&2\Phi\, \Phi_{\barr }^2 + 2\Phi_{\barr }^\TOP  \Psi_{\barr \barr }\Psi_{\barr }+\Psi_{\barr }^\TOP  \Phi_{\barr \barr }\Psi_u=8 \barr ^2\Phi,
\label{r3}\\
&&2\Phi\,\Phi_{\barr }^\TOP  \Psi_{\barr } + \Psi_{\barr }^\TOP  \Psi_{\barr \barr }\Psi_{\barr }=8 \barr ^2\Psi.
\label{r4}
\end{eqnarray}
Since $f$ is harmonic, we also have
\begin{equation}\label{eq:harm}
\Delta \Phi( \barr )=\Delta \Psi(\barr )=0.
\end{equation}
By virtue of (\ref{r1}) we see that that the eigenvalues of $\Phi$ are $\pm 1$ and $0$, and by (\ref{eq:harm}) the multiplicities of $\pm1$ are equal. Denote by $U\oplus V\oplus W$ the associated with $\Phi$ eigen decomposition of $\R{n}\ominus \mathrm{span} (e_n)$, and let
$\barr=u\oplus v\oplus w$ denote the corresponding decomposition of a typical vector $\bar u\in V$. Then
$\Phi= u^2-v^2,$
and
\begin{equation*}
\Psi=\sum_{q}
\Psi_{q },\qquad \Psi_{q}\in u^{q_1 }\otimes v^{q_2 }\otimes w^{q_3 },
\end{equation*}
where $q=(q_1,q_2,q_3)$ and $q_i\ge0$, $q_1+q_2+q_3=3$. Applying the explicit form of $\Phi$, we obtain from  (\ref{r2})
$$
\sum_{q}
((q_1 -q_2 )^2+q_1 +q_2 -2)\,\Psi_{q }=0
$$
and since the non-zero components $\Psi_q$ are linearly independent their coefficients must be zero. This yields
\begin{equation}\label{r11}
\Psi=\Psi_{102}+\Psi_{012}+\Psi_{111}\equiv a+b+c.
\end{equation}
Since $\Delta_{\barr }c=0$ and $a$ and $b$ are linear in $u$ and $v$ respectively, (\ref{harmabw}) follows from (\ref{eq:harm}).
Furthermore,  (\ref{r3}) is equivalent to $(\ref{r17})$-$(\ref{r16r})$. Indeed, we have from (\ref{r11}), $\Psi_{\barr }^2=(a_u+c_u)^2+(b_v+c_v)^2+(a_w+b_w+c_w)^2.$
By Euler's homogeneity theorem, one readily finds that
\begin{equation}\label{rs}
\begin{split}
\Phi_{\barr }^\TOP  \Psi_{{\barr }{\barr }}\Psi_{\barr }&
=2(b_v^\TOP   c_v-a_u^\TOP   c_u+c_v^2-c_u^2+(a_w-b_w)^\TOP   c_w+a_w^2-b_w^2).
\end{split}
\end{equation}
By using the explicit form of $\Phi$ we rewrite  (\ref{r3}) as
\begin{equation*}
\begin{split}
\Phi_{\barr }^\TOP  \Psi_{{\barr }{\barr }}\Psi_{\barr }+\Psi_u^2-\Psi_v^2=4(u^2-v^2)w^2,
\end{split}
\end{equation*}
so applying (\ref{rs})  we obtain
\begin{equation*}
\begin{split}
c_v^2+2a_w^2-c_u^2-2b_w^2+2a_w^\TOP   c_w-2b_w^\TOP   c_w+a_u^2-b_v^2=4(u^2-v^2)w^2.
\end{split}
\end{equation*}
Identifying the terms by homogeneity  shows that (\ref{r3}) is equivalent to (\ref{r17})--(\ref{r16r}).

We proceed similarly with (\ref{r4}). Since the Hessian matrices  $\Psi_{uu}$ and $\Psi_{vv}$ are identically zero and  $2\Phi\,\Phi_{\barr }^\TOP  \Psi_{\barr }=4(u^2-v^2)(a-b)$, (\ref{r4}) becomes
\begin{equation*}
\begin{split}
2\Psi_u^\TOP \Psi_{uv}\Psi_v &+2\Psi_u^\TOP \Psi_{uw}\Psi_w+2\Psi_v^\TOP \Psi_{vw}\Psi_w+\Psi_w^\TOP \Psi_{ww}\Psi_w\\
&=4a(u^2+3v^2+2w^2)+4b(3u^2+v^2+2w^2)+8c(u^2+v^2+w^2).
\end{split}
\end{equation*}
By using (\ref{r11}) and identifying the  terms by homogeneity, we obtain
\begin{eqnarray}
&&   a_u^{\TOP}c_{uv}b_v=0,  \label{s005h}\\
&&{2c_u^{\TOP}c_{uv}c_v+(a_u^2+b_v^2)^{\TOP}_wc_w+2(a_u^{\TOP}c_{u})^{\TOP}_wa_w +2(b_v^{\TOP}c_{v})^{\TOP}_wb_w}=8cw^2,  \label{s113h}\\
&&{  (a_w^{\TOP}b_{w})^{\TOP}_wc_w}=0,  \label{s221h}\\
&&{c_wa_{ww}c_w}=0,  \label{s320h}\\
&&   (a_u^{\TOP}c_{u})_wb_w=0,  \label{s023h}\\
&&{  (c_v^2+a_w^2)^{\TOP}_wc_w}=8cu^2,  \label{s311h}\\
&&{ (2c_v^2+a_w^2)^{\TOP}_wa_w}=8au^2,  \label{s302h}\\
&&{  2a_u^{\TOP}c_{uv}c_v+(a_u^2+b_v^2)^{\TOP}_wa_w}=8aw^2,  \label{s104h}\\
&&   2(a_w^{\TOP}c_{w})^{\TOP}_uc_u+2c_w^{\TOP}c_{wu}a_u+2a_w^{\TOP}b_{ww}b_w+b_w^{\TOP}a_{ww}b_w=12av^2,  \label{s122h}
\end{eqnarray}
and additionally six   equations obtained from (\ref{s320h})--(\ref{s122h}) by  permutation $(u,a) \leftrightarrow (v,b)$.
Then identities (\ref{s320h}), (\ref{s311h}), (\ref{s104h}), (\ref{s113h}), (\ref{s221h})  as well as their permutations  are corollaries of (\ref{r15})--(\ref{r17}).
Indeed, (\ref{s221h}) and (\ref{s320h}) easily follows from (\ref{r15}). Using the first identity in (\ref{r16}),
$$
(c_v^2+a_w^2)^\TOP_wc_w=(4u^2w^2-a_w^2)^\TOP_wc_w=8u^2c-2a_wc_{ww}c_w=8cu^2,
$$
so that (\ref{r15}) and (\ref{r16}) implies  (\ref{s311h}). Similarly,
$$
2a_u^{\TOP}c_{uv}c_v=a_u^{\TOP}(c_v^2)_u=a_u^{\TOP}(4u^2w^2-2a_w^2)_u=8aw^2-4a_ua_{uw}a_w,
$$
and by virtue of (\ref{r17}) we have
$
(a_u^2+b_v^2)^{\TOP}_wa_w=2(a_u^2)^{\TOP}_wa_w=4a_u^{\TOP}a_{uw}a_w.
$
Thus, (\ref{s104h}) is a corollary of (\ref{r16}) and (\ref{r17}). Finally, by virtue of (\ref{r15}) and (\ref{r16})
\begin{equation*}
\begin{split}
c_u^{\TOP}c_{uv}c_v+(a_u^2)_wc_w+2(a_u^{\TOP}c_{u})^{\TOP}_wa_w\frac{1}{2}c_u^{\TOP}(c_{v}^2+ 2a_{w}^2)_u +2a_u^\TOP( a_{w}c_w)_u,
\end{split}
\end{equation*}
and after permutation $(u,a) \leftrightarrow (v,b)$
$$
c_u^{\TOP}c_{uv}c_v+(b_v^2)_wc_w+2(b_v^{\TOP}c_{v})^{\TOP}_wb_w=4cw^2.
$$
Summing these identities  yields (\ref{s113h}). Next, (\ref{s005h}) is equivalent to (\ref{s005}) and (\ref{s122}) is equivalent to (\ref{s122h}) modulo (\ref{r15}). Finally, (\ref{s023}) is equivalent to (\ref{s023h}) modulo (\ref{r15}). For instance,
$$
0=(a_u^{\TOP}c_{u})_wb_w=a_u^{\TOP}c_{uw}b_w+c_u^{\TOP}a_{uw}b_w=a_u^{\TOP}(c_{u}^\TOP b_w)_u+c_u^{\TOP}a_{uw}b_w=c_u^{\TOP}a_{uw}b_w.
$$
The proposition is proved completely.
\end{proof}

\subsection{A matrix representation of the degenerate form}

Now we are going to explain how to extract the type of a non-isoparametric eigencubic from its degenerate  representation (\ref{degeeabc}).
To this end we convert (\ref{degeeabc}) into the matrix form
\begin{equation}\label{degee}
\begin{split}
f&=(u^2-v^2)x_n+\frac{1}{\sqrt{2}}w^\TOP   A_uw+ \frac{1}{\sqrt{2}}w^\TOP   B_vw+2v^\TOP   M_u w.
\end{split}
\end{equation}
where
\begin{equation}\label{abceq}
a=\frac{1}{\sqrt{2}}w^\TOP   A_uw,\qquad b= \frac{1}{\sqrt{2}}w^\TOP   B_vw, \qquad c= 2v^\TOP   M_u w= 2u^\TOP   N_v w.
\end{equation}
Here $A_u=\sum_{i=1}^m u_iA_i$ etc., and  $A_i, B_i$ are symmetric matrices of size $r\times r$, $r=n-1-2m$, and matrices $M_i, N_i$ are of size $m\times r$, $1\le i\le m$.
Then (\ref{harmabw})--(\ref{r16r}) readily yield
\begin{eqnarray}
&&\sum_{i=1}^m (w^\TOP A_i w)^2=\sum_{i=1}^m (w^\TOP B_i w)^2,\\
&&\trace A_u=\trace B_v=0,\label{tracee}\\
\label{Axcube}
&&A_u^3=u^2A_u, \qquad \qquad \qquad  B_v^3=v^2B_v,\\
\label{ABtrace}
&&\trace A_u^2B_v=\trace A_uB_v^2=0,\\
\label{AM}
&&A_u^2+M_u^\TOP M_u=\mathbf{1}_{m},\qquad \,\,\,B_v^2+N_v^\TOP N_v=\mathbf{1}_{m},\\
\label{AM1}
&&A_uM_u^\TOP=B_vN_v^\TOP=0.
\end{eqnarray}

\begin{definition*}
The representations (\ref{degeeabc}) and (\ref{degee}) are called the \textit{degenerate forms} of the corresponding radial (non-isoparametric) eigencubic $f$.
\end{definition*}

\begin{proposition}\label{pr:AB}
In the above notation, there exist nonnegative integers $p$ and $q$  such that
$2p+q=r$  and for any  $u,v\in \R{m}$, $u,v\ne 0$, for any $i$, $1\le i\le m$,  the spectrum of any of the matrices $\frac{1}{|u|}A_u$, $\frac{1}{|v|}B_v$, $A_i$ and $B_i$ is $\{\pm1,0\}$, where each eigenvalue $1$ and $-1$ has multiplicity $p$, and the eigenvalue $0$ has multiplicity $q$.

\end{proposition}

\begin{proof}
We may assume that $m\ge 1$, otherwise the statement of the proposition is trivial. The assumption $A_u^3=u^2 A_u$ implies that for any $u\ne 0$, the eigenvalues of $A_u$ are $\pm|u|$ and $0$. Denote by $p^{\pm}(u)$ and $q(u)$ the multiplicities of $\pm|u|$ and $0$ respectively. By virtue of (\ref{tracee}),$\trace A_u=p^+(u) - p^-(u)=0.$ On the other hand,$\trace A^2_u=p^+(u) + p^-(u).$
The continuous functions $p^+(u)$ and $p^-(u)$ are integer-valued, thus they are identically constant. Denote by $p$ the common value of $p^\pm(u)$. Then $q(u)=r-2p\equiv q$ also is a constant. Since each matrix $A_i$ is the specialization of $A_{u}$ when $u$ is the coordinate vector $e_i$, we conclude that the spectrum of $A_i$ is the same as that of $A_u$.

Similarly, one defines $p'$ and $q'$ for $\frac{1}{|v|}B_v$, $v\ne 0$. To see that  $p=p'$ and $q=q'$,  we apply an iterated  Laplacian to  (\ref{r16}). We find by virtue of  harmonicity of $c$ and (\ref{r16}), $\Delta_v c^2=8w^2u^2-4a_w^2$, hence in view of $r=2p+q$
\begin{equation}\label{ser1}
\begin{split}
\Delta_w \Delta_v c^2&=16(2p+q)u^2-8\trace a_{ww}^2
=16qu^2.
\end{split}
\end{equation}
Applying the $u$-Laplacian to (\ref{ser1}) we get
\begin{equation*}\label{cc}
\Delta_u\Delta_w\Delta_v c^2=32mq.
\end{equation*}
Arguing similarly one finds $\Delta_v\Delta_w\Delta_u c^2=32mq'$ which yields$q=q'$. Since $2p+q=2p'+q'$ we also have $p=p'$. The proposition is proved.
\end{proof}

\begin{definition*}
Given a degenerate form of a radial (non-isoparametric) eigencubic, we associate by virtue of Proposition~\ref{pr:AB} a triple of integers $(m,p,q)$ called the \textit{signature} of the degenerate form.
\end{definition*}

Observe that the type of the degenerate form is determined by virtue of the following formulas which are the corollary of the above definitions:
\begin{equation}\label{traceABsq}
\trace A_u^2=2pu^2, \qquad \trace B_v^2=2pv^2.
\end{equation}
\begin{equation}\label{ser10}
\begin{split}
\trace M_u^\TOP M_u\equiv \frac{1}{4}\trace c_{wv}c_{vw}=qu^2, \qquad
\trace N_v^\TOP N_v\equiv \frac{1}{4}\trace c_{wu}c_{uw}=qv^2.
\end{split}
\end{equation}

The connection between the signature of the degenerate form and the type of an arbitrary non-isoparametric  radial eigencubic is given in the following
\begin{proposition}\label{pr:degtonorm}
Let $f$ be a non-isoparametric eigencubic of type $(n_1,n_2)$ which has the degenerate form of signature $(m,p,q)$. Then
\begin{equation}\label{degtonorm}
n_1=q,\quad n_2=m+p+1-q,\quad n_3=m+p-1+q.
\end{equation}
\end{proposition}

\begin{proof}
We only show that $n_1=q$ because the remaining two identities in (\ref{degtonorm}) readily follow from (\ref{a111}) and $n=2m+2p+q+1$. To this end, we assume that $f$ is given by (\ref{degeeabc}) so that
$$
\mathrm{Hess}(f)\left(
        \begin{array}{ll}
          f_{x_nx_n} & f_{x_n\bar x} \\
          f_{\bar x x_n} & f_{\barr\barr}\\
        \end{array}
      \right)
\equiv \left(
        \begin{array}{ll}
          0\,\, & (\Phi_{\bar x})^\TOP \\
          \Phi_{\bar x}& x_n\Phi_{\bar x\barr }+\Psi_{\barr\barr} \\
        \end{array}
      \right),
$$
so that
$$
\trace \mathrm{Hess}^3(f)=x_n^3\trace \Phi^3_{\bar x \bar x}+3x_n^2\trace \Phi_{\barr\barr}^2\Psi_{\barr\barr}
+3x_n(\trace \Phi_{\barr\barr}\Psi_{\barr\barr}^2+\trace \Phi_{\barr}^\TOP\Phi_{\barr\barr}\Phi_{\barr})+\ldots
$$
where the dots stands for the terms which contain no $x_n$. Comparing this with (\ref{trace_iden}) and (\ref{degeeabc}) yields by virtue of $\lambda(f)=-8$ the tautological identities $\trace \Phi^3_{\bar x \bar x}=\trace \Phi_{\barr\barr}^2\Psi_{\barr\barr}=0$ and additionally
\begin{equation}\label{n1n1q}
\trace \Phi_{\barr\barr}\Psi_{\barr\barr}^2+\trace \Phi_{\barr}^\TOP\Phi_{\barr\barr}\Phi_{\barr}=8(n_1-1)(v^2-u^2).
\end{equation}
In view of $\Phi=u^2-v^2$, we find $\trace \Phi_{\barr}^\TOP\Phi_{\barr\barr}\Phi_{\barr}=8(v^2-u^2)$.
Furthermore,
\begin{equation*}
\begin{split}
\trace a_{wu}a_{uw}&=\frac{1}{2}\Delta_w a_u^2=\text{by (\ref{r17})}\frac{1}{2}\Delta_w b_v^2=\trace b_{wv}b_{vw},\\
\trace a_{wu}c_{uw}&=\trace b_{wv}c_{vw}=\text{by (\ref{r15})}=0,\\
\trace c_{wu}c_{uw}&=\text{by (\ref{ser10})}=4qv^2.
\end{split}
\end{equation*}
and taking into account that $\Psi_{uu}=\Psi_{vv}=0$,
\begin{equation*}
\begin{split}
\trace \Psi_{{\bar x}{\bar x}}\Phi_{{\bar x}{\bar x}}\Psi_{{\bar x}{\bar x}}&2\trace (\Psi_{wu}\Psi_{uw}-\Psi_{wv}\Psi_{vw})\\
&=2\trace (a_{wu}+c_{wu})(a_{uw}+c_{uw})-2\trace (b_{wv}+c_{wv})(b_{vw}+c_{vw})\\
&=2(\trace c_{wu}c_{uw}-\trace c_{wv}c_{vw})=-8q(u^2-v^2).
\end{split}
\end{equation*}
Substituting the found relations into (\ref{n1n1q}) yields $n_1=q$.
\end{proof}

In the converse direction, the the signature of the degenerate form of a radial eigencubic is not well determined in general by its type. However, it is well determined, in the most important for us case of exceptional eigencubics.

\begin{corollary}\label{cor:normdeg}
Let $f$ be an exceptional eigencubic of type $(n_1,n_2)$, $n_2\ne 0$. Then $n_2=3\ell+2$, $\ell\in \{1,2,3,4\}$ and the signature of the degenerate form of $f$ is determined  by the formula
\begin{equation}\label{normtodeg}
(m,p,q)=(\ell+n_1+1, 2\ell, n_1).
\end{equation}
\end{corollary}

\begin{proof}
The first statement follows immediately from Proposition~\ref{pro:finiten2} and the inequality $n_2\ne0$ for non-isoparametric eigencubics. Next, we assume that $f$ is given by (\ref{degee}). Then
$$
\trace \mathrm{Hess}^2(f)=8mx_n^2+\ldots
$$
where the dots stands for the terms with degree of $x_n$ lower than 2. By Proposition~\ref{th:dege}, $\lambda(f)=-8$, hence (\ref{eq:norm}) together with $n_2=3\ell+2$  yield
$m=n_1+\ell+1$.  Applying (\ref{degtonorm}) we get (\ref{normtodeg}).
\end{proof}

\subsection{Proof of Theorem~\ref{thE}}
First note a simple corollary of (\ref{Muntzer4}):
\begin{equation}\label{Muntzer100}
h\in \mathrm{Isop}(m_1,m_2) \quad \Leftrightarrow \quad -h\in \mathrm{Isop}(m_2,m_1).
\end{equation}
Let us define two auxiliary polynomials
$$
\widehat h_0(t)=\frac{c_w^2}{4}, \qquad \widehat h_1(t)=\frac{1}{4}(u^2-v^2)^2+\frac{c_w^2}{4},
$$
which are quartic polynomials in the variable $t=(u,v)\in \R{2m}$. From (\ref{r16r}) we obtain
\begin{equation*}
\begin{split}
\frac{1}{4}|\nabla_{u}|c_w|^2|^2=&\sum_{i,j,k}c_{w_i} c_{w_iu_k}c_{u_kw_j}c_{w_j}
=\frac{1}{2}\sum_{i,j}(c_{u}^2)_{w_iw_j}c_{w_i}c_{w_j}\\
&=\sum_{i,j}(4u^2\delta_{ij}c_{w_i}c_{w_j}-2(b_w^2)_{w_iw_j}c_{w_i}c_{w_j})\\
&=4u^2c_w^2- 4\sum_{i,j}b_{w_kw_i}b_{w_kw_j}c_{w_i}c_{w_j}\\
&=4u^2c_w^2- 4\sum_{i,j}(b_w^\TOP c_w)_{w_k}(b_w^\TOP c_w)_{w_k}=4u^2c_w^2\\
\end{split}
\end{equation*}
where the last equality is by (\ref{r15}).
Thus,
\begin{equation}\label{h00}
|\nabla_t \widehat h_0|^2=\frac{1}{64}\bigl(|\nabla_{u}|c_w|^2|^2+|\nabla_{u}|c_w|^2|^2\bigr)=\frac{1}{16}(u^2+v^2)c_w^2=4u^2\widehat h_0.
\end{equation}
Since $\Delta_w c^2=2c_w^2$, we obtain from (\ref{ser1})
$$
\Delta_t \widehat h_0=\Delta_u\widehat h_0+\Delta_v\widehat h_0=\frac{1}{8}(\Delta_w \Delta_u c^2+\Delta_w \Delta_v c^2)=2qt^2.
$$
This implies that $h_0$ satisfies the Cartan-M\"unzner equations:
\begin{equation*}\label{isop-gamma}
|\nabla_u h_0|^2=16t^6, \qquad
\Delta_u h_0=8(m+1-2q)t^2,
\end{equation*}
and comparison the latter relation  with (\ref{Muntzer4}) and (\ref{dimension1}) yields the first claimed inclusion.

Now consider $h_1$. Define $\widehat h_1=\widehat h_0+\frac{1}{4}(u^2-v^2)^2$, so that
$$
|\nabla\widehat h_1|^2=|\nabla\widehat h_0|^2+2(u^2-v^2)(\scal{\nabla_u \widehat h_0}{u}-\scal{\nabla_v \widehat h_0}{v})+(u^2-v^2)^2t^2.
$$
Since $c\in u\otimes v\otimes  w$ we find
$$
\scal{\nabla_u \widehat h_0}{u}=\frac{1}{2}c_w^\TOP c_{wu}u=\frac{1}{2}c_w^\TOP c_{w}=2h_0,
$$
which, in particular,  implies  $\scal{\nabla_u \widehat h_0}{u}=\scal{\nabla_v \widehat h_0}{v}$. Thus, by (\ref{h00})
$$
|\nabla \widehat h_1|^2=|\nabla \widehat h_0|^2+(u^2-v^2)^2t^2=4t^2\widehat h_0+(u^2-v^2)^2t^2=4\widehat h_1t^2.
$$
Now it is easy to see that
\begin{equation*}\label{isop-gamma1}
|\nabla_t h_1|^2=16t^6, \qquad
\Delta_t h_1=8(m+1-2q)t^2,
\end{equation*}
which yields the second claimed inclusion. Thus the theorem is proved.

\subsection{The non-existence results}

First combining Corollary~\ref{cor:normdeg} with Theorem~\ref{thE} we obtain the following

\begin{corollary}
\label{cor:iso}
To any exceptional eigencubic $f$ of type $(n_1,3\ell+2)$ one can associate two isoparametric quartic polynomials $h_0$ and $h_1$ by virtue of Theorem~\ref{thE},  $h_0 \in \mathrm{Isop}(n_1-1,\ell+1)$ and $h_1\in \mathrm{Isop}(n_1,\ell)$.
\end{corollary}

On the other hand, in \cite{OT1} Ozeki and Takeuchi by using representations of Clifford algebras produced two infinite series of isoparametric families with four principal curvatures. In \cite{FKM}, Ferus, Karcher, and M\"nzner generalized the Ozeiki-Takeuchi approach to produce an infinite family of isoparametric hypersurfaces. More precisely, a quartic polynomial $f$ is said to be of OT-FKM-type \cite[\S~4.7]{Cecil1} if it is congruent to
$$
F_{\mathcal{H}}:=t^4-2\sum_{i=0}^{s}(t^\TOP H_i t)^2, \qquad t\in \R{2m},
$$
where $\mathcal{H}=\{H_0,\ldots, H_s\}\in \mathrm{Cliff}(\R{2m},s)$.
It is easily verified that
\begin{equation}\label{Hsys}
F_{\mathcal{H}}\in \mathrm{Isop}(s,m-s-1).
\end{equation}
Since the class $\mathrm{Cliff}(\R{2m},s)$ is non-empty if and only if
\begin{equation}\label{HurEx}
s\le \rho(m),
\end{equation}
the existence of the quartic polynomial (\ref{Hsys}) of OT-FKM-type is equivalent to the inequality
\begin{equation}\label{HurExm1}
\min\{ m_1,m_2\}\le \rho(m_1+m_2+1).
\end{equation}
A culmination of this story is the following deep classification result due to  Cecil-Chi-Jensen \cite{CJ} and Immerwoll \cite{Immer}: if $h\in \mathrm{Isop}(m_1,m_2)$ and either of the inequalities  $m_2\ge 2m_1-1$ or $m_1\ge 2m_2-1$ holds then $h$ is OT-FKM-type.

\begin{proposition}\label{pr:n1=2}
There are no exceptional eigencubics of type $(2,8)$, $(2,14)$, $(2,26)$.
\end{proposition}

\begin{proof}
We argue by contradiction and suppose  that there exist an exceptional eigencubic $f$ of type $(n_1,n_2)\in \{(2,8), (2,14), (2,26)\}$. By Corollary~\ref{cor:iso} we can associate to $f$ an isoparametric quartic polynomial $h_1\in \mathrm{Isop}(m_1,m_2)$ with $m_1=n_1$ and $m_2=\frac{n_2-2}{3}$. Since $n_2\ge 8$ we have
$$
m_2+1-2m_1=\frac{n_2-2}{3}+1-n_1=\frac{n_2-5}{3}\ge 0
$$
thus, by Cecil-Chi-Jensen-Immerwoll theorem, $h_1$ is of OT-FKM-type. This  by virtue of (\ref{HurExm1}) yields $2=\min\{m_1,m_2\}\le \rho(m_1+m_2+1)$. But, for our choice of $(n_1,n_2)$ the number  $m_1+m_2+1=\frac{n_2+7}{3}$ is  odd, hence $\rho(m_1+m_2+1)=1$, a contradiction.

\end{proof}

\begin{proposition}
There are no exceptional eigencubics of type $(3,8)$.
\end{proposition}

\begin{proof}
Again, we argue by contradiction and suppose  that $f$ is an exceptional eigencubic of type $(3,8)$. By Corollary~\ref{cor:normdeg}, $\ell=2$, $q=n_1=3$ and $m=q+\ell+1=6$. This implies by Corollary~\ref{cor:iso} and (\ref{Muntzer100}) that the isoparametric quartic polynomials $h_i$ associated to $f$ satisfy
\begin{equation}\label{isoh}
h_0\in \mathrm{Isop}(2,3),\qquad -h_1\in \mathrm{Isop}(2,3).
\end{equation}
By Cecil-Chi-Jensen-Immerwoll theorem both $h_0$ and $-h_1$ are of OT-FKM-type. Hence there exists a Clifford system $\mathcal{H}=\{H_0,\ldots, H_s\}\in \mathrm{Cliff}(\R{2m},s)$ associated with $h_0$ and $\mathcal{F}=\{F_0,\ldots, F_r\}\in \mathrm{Cliff}(\R{2m},r)$ associated with $h_1$. In virtue of $m=6$ the inequality (\ref{HurEx}) yields $s,r\le 2$. Then we infer from (\ref{Hsys}) and (\ref{isoh}) that $s=r=2$. By using (\ref{hpolys}) we have
$$
h_0\equiv (u^2+v^2)^2-2 c_w^2=t^4-2\sum_{i=0}^{2}(t^\TOP H_i t)^2
$$
and
$$
-h_1\equiv u^4-6u^2v^2+v^4+2 c_w^2=t^4-2\sum_{j=0}^{2}(t^\TOP F_i t)^2,
$$
hence eliminating $c_w^2$ yields
$$
\sum_{i=0}^{2}(t^\TOP H_i t)^2+\sum_{j=0}^{2}(t^\TOP F_i t)^2=4u^2v^2.
$$
The latter identity implies the following block forms associated with the decomposition $t=u\oplus v$:
$$
H_i=\left(
        \begin{array}{cc}
          \textbf{0} & S_i \\
          S_i^\TOP &  \textbf{0}\\
        \end{array}
      \right),
\qquad F_j=\left(
        \begin{array}{cc}
          \textbf{0} & S_{j+3} \\
          S_{j+3}^\TOP &  \textbf{0}\\
        \end{array}
      \right),\quad i,j=0,1,2,
$$
thus
$$
\sum_{i=0}^{5}(u^\TOP S_i v)^2=4u^2v^2,\qquad u,v\in\R{6}.
$$
But the latter identity is a composition formula of size $[6,6,6]$ (see, for instance, \cite{Shapiro}) which contradicts to the Hurwitz theorem stating that a composition formula of size $[k,k,k]$ does exist only if $k=1,2,4,8$. The contradiction proves the proposition.
\end{proof}


\begin{thebibliography}{FKM}

\bibitem[A]{Abresh}
Abresch U., {Isoparametric hypersurfaces with four or six distinct principal curvatures}, \textit{Math. Ann.} \textbf{264}~(1983), 283--302

\bibitem[Ad]{Adams}
Adams J.F., Lectures on Exceptional Lie Groups, Chicago Lectures in Mathematics. University of Chicago Press, Chicago, IL, 1996

\bibitem[ABS]{Atiyah}
Atiyah M.,  Bott R.,  Shapiro M., Clifford Modules, \textit{Topology} \textbf{3}(1964), 3--38

\bibitem[Ba]{Baez}
Baez, John. The Octonions, \textit{Bull. Amer. Math. Soc.}, 39(2002) 145--205.

\bibitem[BW]{Baird}
Baird P., Wood J.C., Harmonic morphisms between Riemannian manifolds. London Mathematical Society Monographs. New Series, 29. The Clarendon Press, Oxford University Press, Oxford, 2003.


\bibitem[BG]{BG}
Bombieri, E.; Giusti, E. Harnack's inequality for elliptic differential equations on minimal surfaces.
Invent. Math. \textbf{15}~ (1972), 24--46.

\bibitem[BGG]{BGG}
Bombieri, E., de Giorgi, E., Giusti, E.: Minimal cones and the Bernstein
problem. {Invent. Math}. \textbf{82}, 243--269 (1968)

\bibitem[Ca1]{Ca1}
Calabi E., Improper affne hypersurfaces of convex type and a generalization of a theorem
by K. J\"{o}rgens, \textit{Michigan Math. J.}, \textbf{5}~(1958), 105--126.



%

\bibitem[Car] {Cartan}
\'{E}. Cartan, {Sur des familles remarquables d'hypersurfaces isoparam\'{e}triques dans les espaces sph\'{e}riques},
\textit{Math. Z.} \textbf{45} (1939), 335--367.

\bibitem[C1]{Cecil}
Cecil,  Th.E., Isoparametric and Dupin hypersurfaces, \textit{SIGMA Symmetry Integrability Geom. Methods Appl}. \textbf{4} (2008), Paper 062.

\bibitem[C2]{Cecil1}
Cecil, Th.E., Lie Sphere Geometry With Applications to Submanifolds, 2nd ed., Universitext, Springer, New York, 2008


\bibitem[CCC]{CJ}
T. Cecil, Q. S. Chi and G. Jensen, {Isoparametric hypersurfaces with four principal curvartures}, \textit{Ann. of Math.}, \textbf{166}~(2007), 1--76.



\bibitem[Fa]{Farina} Farina A., Two results on entire solutions of Ginzburg-Landau system in higher dimensions,
J. Funct. Anal. \textbf{214}~ (2004), 386--395.

\bibitem[Fl]{Fleming}
Fleming W. H., Flat chains over a finite coefficient group, \textit{Trans. Amer. Math. Soc.}, \textbf{121}~(1966),
160--186.

\bibitem[FKM]{FKM}
Ferus D., Karcher H. and M\"{u}nzner H.-F., Cliffordalgebren und neue isoparametrische
Hyperfl\"{a}chen, \textit{Math. Z.}, \textbf{177} (1981), 479--502.

\bibitem[F]{Fomenko}
Fomenko, A.T.: {Variational principles of topology. Multidimensional minimal surface theory}. Mathematics and its Applications, 42. Kluwer Acad. Publ. Group, Dordrecht (1990)


\bibitem[GX]{GeXie}
Ge J.,  Xie Y., {Gradient map of isoparametric polynomial and its application to Ginzburg-Landau system.} \textit{J. Funct. Anal.}, \textbf{258}(2010), 1682--1691.

\bibitem[GN]{GN}
Godlinski M., Nurowski P., $GL(2,\R{})$-geometry of ODE's, J. Geom. Physics, \textbf{60}(2010), 991--1027

\bibitem[GW]{Wallach}
Goodman, R., Wallach, N.R.:
Representations and invariants of the classical groups. Cambridge Univ. Press, Cambridge  (1998)

\bibitem[HL]{HL}
Harvey R., Lawson H. B., Calibrated geometries, \textit{Acta Math.} \textbf{148}(1982), 47--157.


\bibitem[H1]{Hsiang66}
Hsiang W.~Y., On the compact homogeneous minimal submanifolds, \textit{Proc. Nat. Acad. Sci. U.S.A}. \textbf{56}~(1966) 5--6.


\bibitem[H]{Hsiang67}
Hsiang W.~Y.,
{Remarks on closed minimal submanifolds in the standard Riemannian $m$-sphere}, \textit{{J. Diff. Geom.}} \textbf{1}(1967), 257--267.


\bibitem[HL]{Hsiang-Lawson}
Hsiang W.~Y., Lawson H.~B., Minimal submanifolds of low cohomogeneity, J. Diff Geom. \textbf{5}~ (1971), 1-38.

\bibitem[Hu]{Huse}
Husemoller D.,
Fibre bundles. 3rd edition. Graduate Texts in Mathematics, 20. Springer-Verlag, New York, 1994


\bibitem[Im]{Immer}
Immervoll S., {The classification of isoparametric hypersurfaces
with four distinct principal curvatures}, \textit{Ann. Math.}, \textbf{168} (2008), 1011--1024.




\bibitem[J]{Jorg}
J\"orgens K., \"Uber die L\"osungen der Didderentialgleichung $rt-s^2 = 1$, \textit{Math. Ann.}, \textbf{127}~(1954), 130--134.


\bibitem[L]{Lawson}
Lawson H.~Bl., Complete minimal surfaces in $S\sp{3}$. \textit{Ann. of Math.}  \textbf{92} (1970) 335--374.

\bibitem[MMM]{MMM}
Massari U., Miranda M., Miranda M., Jr.: The Bernstein theorem in higher dimensions. \textit{Boll. Unione Mat. Ital.} (9) \textbf{1}~(2008), no. 2, 349--359.

\bibitem[Mi]{Milnor}
Milnor, J., Singular points of complex hypersurfaces. Annals of Mathematics Studies, No. 61. Princeton University Press, Princeton, NJ. 1968

\bibitem[M]{Miranda}
Miranda M., Recollections on a conjecture in mathematics. \textit{Mat. Contemp.} \textbf{35}~(2008), 143--150.


\bibitem[Mu1]{Mun1}
M\"unzner H.-F., \textit{Isoparametrische Hyperfl\"achen in Sph\"aren}. {Math. Ann.} \textbf{251} (1980), 57--71.

\bibitem[Mu2]{Mun2}
M\"unzner H.-F., Isoparametrische Hyperfl\"achen in Sph\"aren, II. Uber die Zerlegung der Sph\"are in Ballb\"undel, \textit{Math. Ann.} \textbf{256} (1981), 215–-232.


\bibitem[NV1]{NV1}
Nadirashvili N., Vl\v{a}du\c{t} S., \textit{Geom. Funct. Anal.}, Nonclassical solutions of fully nonlinear elliptic equations,
17(2007), 1283--1296

\bibitem[NV2]{NV2}
Nadirashvili N.,  Vl\v{a}du\c{t} S.,
Nonclassical Solutions of Fully Nonlinear Elliptic Equations II. Hessian Equations and Octonions, arXiv:0912.3126

\bibitem[Ni]{Nitsche}
Nitsche J.C.C., Lectures on minimal surfaces. Vol. 1. Introduction, fundamentals, geometry and basic boundary value problems. Cambridge University Press, Cambridge, 1989

\bibitem[No]{Nomi}
Nomizu K., Some results in E. Cartan's theory of isoparametric families of hypersurfaces.
\textit{Bull. Am. Math. Soc.} \textbf{79}~ (1973), 1184--1189.

\bibitem[Nu]{Nur}
Nurowski P., Distinguished dimensions for special Riemannian geometries, J. Geom. Phys. \textbf{58}~(2008),  1148--1170.
Preprint: arXiv:math/0601020

\bibitem[Os]{Osserman}
Osserman, R. The minimal surface equation. Seminar on nonlinear partial differential equations (Berkeley, Calif., 1983), 237--259, Math. Sci. Res. Inst. Publ., 2, Springer, New York, 1984.

\bibitem[OT1]{OT1}
H.~Ozeki \& M.~Takeuchi, {On some types of isoparametric hypersurfaces in spheres. I.}
\textit{T\^{o}hoku Math. J.} (2) \textbf{27} (1975), no. 4, 515--559.

\bibitem[OT2]{OT2}
H.~Ozeki \& M.~Takeuchi, {On some types of isoparametric hypersurfaces in spheres. I.}
\textit{T\^{o}hoku Math. J.} (2) \textbf{28} (1976), no. 1, 7--55.


\bibitem[Pe]{Perdomo1}
Perdomo O., Characterization of order 3 algebraic immersed minimal surfaces of $S^3$. {Geom. Dedicata}, \textbf{129}~ (2009), 23--34

\bibitem[P]{Pogor}
Pogorelov A.V., The Multidimensional Minkowski Problem, Wiley, New York, 1978.

\bibitem[Sh]{Shapiro}
Shapiro D. B., \textit{Compositions of quadratic forms},
de Gruyter Expositions in Mathematics, 33. Walter de Gruyter \& Co., Berlin, 2000.

\bibitem[S1]{SimonL}
Simon, L.: The minimal surface equation. Geometry, V, 239--272, Encyclopaedia Math. Sci., \textbf{90}, Springer, Berlin, (1997)

\bibitem[S2]{Simon89}
Simon, L.: Entire solutions of the minimal surface equation, {J. Diff. Geom.}, \textbf{30}, 634--688 (1989)


\bibitem[SS]{SS}
Simon, L., Solomon, B.: Minimal hypersurfaces asymptotic to quadratic cones in $\R{n+1}$. {Invent. math.}, \textbf{86}, 535--551 (1986)



\bibitem[SJ]{Simons}
Simons, J.: Minimal varieties in Riemannian manifolds. {Ann. of Math.},
\textbf{88},  62--105 (1968)

\bibitem[Sm]{Smoczyk}
Smoczyk K., Algebraic selfsimilar solutions of the mean curvature flow. Preprint, 2009.


\bibitem[So1]{BS1}
Solomon B., \textit{The harmonic analysis of cubic isoparametric minimal hypersurfaces. II. Dimensions $3$ and $6$}.
Amer. J. Math. \textbf{112} (1990), no. 2, 157--203.

\bibitem[So2]{BS2}
Solomon B., \textit{The harmonic analysis of cubic isoparametric minimal hypersurfaces. II. Dimensions $12$ and $24$}.
Amer. J. Math. \textbf{112} (1990), no. 2, 205--241.



\bibitem[St]{Stoltz}
Stolz S., \textit{Multiplicities of Dupin hypersurfaces}, Invent. Math. \textbf{138}~(1999), 253--279.


\bibitem[Ta]{Taka}
Takahashi T., Minimal immersions of Riemannian manifolds, \textit{Math. Soc. Japan}, \textbf{18}~(1966), 380--385.


%

\bibitem[T1]{TkCartan}
Tkachev V.G., \textit{A generalization of Cartan's theorem on isoparametric cubics}, Proc. Amer. Math. Soc., \textbf{138} ~(2010), 2889--2895.
reprint arXiv:1001.2387

\bibitem[T2]{TkCliff}
Tkachev V.G., \textit{Minimal cubic cones via Clifford algebras}, {Complex Anal. Oper. Theory}, \textbf{4}~ (2010), no 3, 685--700. Preprint arXiv:1003.0215



\bibitem[Th1]{Thorb}
Thorbergsson G.,
A survey on isoparametric hypersurfaces and their generalizations. Handbook of differential geometry, Vol. I, 963--995, North-Holland, Amsterdam, 2000.

\bibitem[Th2]{Thorb2}
Thorbergsson G., Singular Riemannian foliations and isoparametric submanifolds, \textit{Milan J. of Math.}, \textbf{78}~ (2010), 355--370




\bibitem[Wa]{Waerden}
Waerden, van der: Algebra I, Springer-Verlag, 1966

%


\end{thebibliography}
\end{document}